\documentclass[a4paper,fleqn]{cas-sc}

\makeatletter
\def\ps@pprintTitle{%
  \let\@oddhead\@empty
  \let\@evenhead\@empty
  \let\@oddfoot\@empty
  \let\@evenfoot\@oddfoot
}
\makeatother

\ExplSyntaxOn
\cs_set:Npn \__first_footerline: {
  \group_begin:
  \small \sffamily
  \group_end:
}
\ExplSyntaxOff

\usepackage[numbers]{natbib}
\newcommand{\upcite}[1]{\textsuperscript{\citep{#1}}}
\usepackage{subcaption}
\usepackage{float}
\usepackage{placeins}
\usepackage{hyperref}
\usepackage{graphicx}
\usepackage{soul}
\usepackage{xcolor}
\hypersetup{
    colorlinks=true,
    linkcolor=blue,
    citecolor=blue,
    urlcolor=blue
}
\newdefinition{property}{Property}[section]
\newtheorem{theorem}{Theorem}[section]
\newtheorem{lemma}[theorem]{Lemma}
\newdefinition{definition}{Definition}[section]
\newdefinition{remark}{Remark}[section]
\newproof{proof}{Proof}
\def\tsc#1{\csdef{#1}{\textsc{\lowercase{#1}}\xspace}}
\tsc{WGM}
\tsc{QE}

\begin{document}
\let\WriteBookmarks\relax
\def\floatpagepagefraction{1}
\def\textpagefraction{.001}

\shorttitle{}

\shortauthors{}

\title [mode = title]{Convergence and error analysis of problem-independent machine learning (PIML) with linear completeness}

\tnotemark[1]

\tnotetext[1]{This work was supported by the National Natural Science Foundation of China (Nos. 12271079, 12494552, 12494550).}

%

\author[1,3]{Ke Tan}


\fnmark[1]

\ead{tank@mail.dlut.edu.cn}




\author[1,3]{Chong-Jun Li}
\cormark[1]
\fnmark[1]

\ead{chongjun@dlut.edu.cn}



\affiliation[1]{organization={School of Mathematical Sciences, Dalian University of Technology},
            city={Dalian},
            state={Liaoning},
            country={PR China}}

\author[2]{Xu Guo}
\fnmark[2]

\ead{guoxu@dlut.edu.cn}



\affiliation[2]{organization={State Key Laboratory of Structural Analysis, Optimization and CAE Software for Industrial Equipment, Department of Engineering Mechanics, Dalian University of Technology},
            city={Dalian},
            state={Liaoning},
            country={PR China}}
\affiliation[3]{organization={Key Laboratory for Computational Mathematics and Data Intelligence of Liaoning Province, Dalian University of Technology},
            city={Dalian},
            state={Liaoning},
            country={PR China}}

\cortext[1]{Corresponding author}

\fntext[1]{}


\begin{abstract}
Problem‑independent machine learning (PIML) is an effective computational approach that combines the finite element method with machine learning. PIML uses neural networks to learn the mapping from material distributions or element geometries to the numerical shape functions of substructures. This approach avoids repeated computation while preserving the conformity of the underlying finite element discretization, thereby significantly improving the efficiency of substructure-based finite element analysis. 
However, existing studies on PIML have primarily focused on numerical accuracy and computational efficiency, while rigorous theoretical analyses of its approximation errors are still limited.
In this paper, we establish an error analysis framework for PIML from a classical finite element perspective, in which the numerical error is bounded in terms of the finite element discretization error and the network generalization error. Our analysis shows that, when only rigid-body displacement invariance is enforced, the resulting error estimate contains a learning-induced term that does not necessarily decrease under mesh refinement. To ensure convergence, we propose an improved PIML scheme with linear completeness. By using linear displacement field constraints to reconstruct the corresponding columns of the transition matrix, we ensure that the resulting numerical shape-function space contains the space of linear polynomials. We derive error estimates and establish convergence results of the improved scheme in the \(H^1\) and \(L^2\) norms. Numerical experiments verify the theoretical results and demonstrate that the improved PIML achieves optimal convergence rates, good generalization performance, and numerical stability.
\end{abstract}


\begin{highlights}
\item We develop an error analysis framework for PIML from a classical finite element perspective, in which the error is controlled by both the mesh size \(h\) and the generalization error \(\varepsilon\) of the matrix learned by PIML.
\item We propose an improved PIML scheme that ensures convergence and establish a rigorous error bound for the proposed scheme.
\item We present numerical experiments that verify the theoretical results and demonstrate the effectiveness of the improved scheme.
\end{highlights}

\begin{keywords}
Problem‑independent machine learning\sep Numerical shape functions \sep  Finite element method \sep Linear completeness \sep Error analysis
\end{keywords}

\maketitle

\section{Introduction}\label{1}
In recent years, deep learning (DL)-based research using deep neural networks (DNNs) for solving partial differential equations (PDEs) has increased rapidly \cite{39}. When combined with different loss functions, many neural network (NN) methods have been developed, including the physics-informed neural networks (PINNs)  \cite{6}, the deep Galerkin method (DGM)  \cite{36,7}, the deep Ritz method (DRM)  \cite{37,8}, and others. Some convergence analyses have been developed for these methods, particularly for PINNs  \cite{9, 10, 11, 12}. 

    Although neural-network-based methods have achieved considerable success, several difficulties remain unresolved when applying them to solve PDEs. One major challenge is the high computational cost. Training a neural network requires substantial memory and computational resources, and it is not only strongly dependent on hardware but also requires a large amount of time. For example, the time consumed by PINNs can be more than two orders of magnitude greater than that of traditional numerical methods \cite{40}. Another difficulty lies in accuracy. Since many NN methods are meshless, when the true solution of a PDE is a high-frequency oscillatory function, the neural networks used in PINNs or DRM will naturally tend to the low-frequency part of the true solution, and thus fail to approximate the true solution  \cite{13}. This phenomenon is called spectral bias  \cite{38,14}. Moreover, the architectures of these neural networks typically consist of compositions of many nonlinear functions, making exact analytical integration of the loss function generally difficult. Some studies have shown that the integration error is a significant factor influencing training performance  \cite{32,15,16,17}.

   In addition, the finite element method (FEM)  \cite{18}, as a traditional numerical method for solving PDEs, is extremely efficient. It has a well-established theory of error estimation and mature programming support. Although different problems pose different challenges to the FEM, numerous FEM variants have been developed to address them, such as the spline element method  \cite{19}, VEM  \cite{20}, SBFEM  \cite{21}, etc. However, as a problem-dependent and mesh-dependent method, the FEM relies strongly on the underlying weak formulation of the PDEs and mesh quality. For problems requiring repeated analyses, repeated assembly of the stiffness matrix and repeated solution of the resulting systems can be computationally expensive \cite{22, 23, 24}.

    Motivated by the complementary advantages and limitations of DNNs and FEM, some researchers have proposed combining the two methods. Recently, problem‑independent machine learning (PIML) has achieved remarkable results. PIML was developed by Huang et al. in  \cite{22, 23} to solve topology optimization problems. It uses a DNN to fit the mapping between the numerical shape functions of the substructure method and the material distribution, thereby avoiding repeated computation of the numerical shape function. Because PIML only needs to learn the variable part of the numerical shape functions and the stiffness matrix computed using traditional FEM shape functions, it is naturally conforming. These features make PIML notably efficient and accurate: a model can be trained in advance on a computer and then be deployed on other computers. Experimental results demonstrate that in topology optimization problems, PIML reduces time consumption by \(60\%\) \cite{23}. Meanwhile, some other network structures have also been used for PIML, such as deep operator networks (DeepONet)  \cite{35,30} and PINNs \cite{42}. Moreover, PIML has been extended to more types of elements, such as the isoparametric elements in \cite{24}.

    It should be noted that although PIML has achieved high efficiency and accuracy in various problems, there is still limited theoretical work on error estimation for the numerical solutions obtained by PIML. In this paper, we analyze the error of PIML from a classical finite element perspective. Our analysis does not depend on the specific neural network architecture used to approximate the numerical shape functions. In fact, we only need to consider the influence of the generalization error on the finite element computation. 

The main contributions of this paper are summarized as follows.

\begin{enumerate}
    \item We develop an error estimation framework for PIML from a classical finite element perspective, in which the error is controlled by both the mesh size \(h\) and the generalization error \(\varepsilon\) of the matrix learned by PIML.
    \item We propose an improved PIML scheme that ensures convergence and establish a rigorous error bound for the proposed scheme.
    \item We present numerical experiments that verify the theoretical results and demonstrate the effectiveness of the improved scheme.
\end{enumerate}

    This paper is organized as follows. In Section \ref{Sect. 2}, we briefly introduce the substructure method for finite element analysis, which is related to the extended multi‑scale finite element method and problem‑independent machine learning (PIML). Some concepts and definitions, which are crucial for the subsequent research, will be introduced here. In Section \ref{Sect. 3}, we derive error estimates for PIML in terms of the mesh size and the generalization error of the learned transition matrix. In Section \ref{Sect. 4}, we propose an improved approach for PIML to ensure convergence only with respect to the mesh size. In Section \ref{Sect. 5}, we present numerical experiments to verify the convergence of both traditional PIML and the modified PIML.


\section{Preliminaries and methodology}\label{Sect. 2}
\subsection{Substructure method for finite element analysis}
Consider the finite element method (FEM) applied to a second-order boundary-value problem (BVP) defined on a bounded domain \(\Omega\subset \boldsymbol{R}^2 \) or \(\boldsymbol{R}^3\), whose weak form is as follows:
\begin{align}\label{ep:weakform}
\mbox{Find }\ \boldsymbol{u}\in V\subset H^1(\Omega)\ s.t.\ a(\boldsymbol{u},\boldsymbol{v})=l(\boldsymbol{v}),\ \forall \boldsymbol{v} \in V\cap H_0^1(\Omega),
\end{align}
with some corresponding forces and boundary conditions, where $a(\boldsymbol{u},\boldsymbol{v})=\int_{\Omega}\boldsymbol{\varepsilon}(\boldsymbol{u}):\boldsymbol{D}(\boldsymbol{x}):\boldsymbol{\varepsilon}(\boldsymbol{v})d\boldsymbol{x}$ satisfies
\begin{align}
a(u,v)&=a(v,u),\label{ep:symmetric}\\
|a(u,v)|&\leq M\parallel{u}\parallel_V\parallel v\parallel_V, \label{ep:continuous}\\
a(v,v)&\geq \alpha\parallel v\parallel_V^2\label{ep:coercive},
\end{align}
with two positive constants \(M\) and \(\alpha\). That is, assume that the bilinear form \(a(u,v)\) is symmetric, continuous, and coercive, by C\'ea's lemma \cite{18}, the finite element solution \(u_h\in V_h\subset H^1(\Omega)\) has the following property:
\begin{align*}
\parallel u-u_h\parallel_V\leq\frac{M}{\alpha}\inf_{v\in V_h}\parallel u-v\parallel_V.
\end{align*}
It is clear that the energy norm \(\parallel\bullet\parallel_V\) is equivalent to the \(H^1\) norm \(\parallel\bullet\parallel_{H^1(\Omega)}\). In practical computation, the numerical solution vector \(\boldsymbol{U}\in\mathbb{R}^n\) corresponding to \(u_h\) is obtained by solving the system of linear equations
\begin{align}
\boldsymbol{K}\boldsymbol{U}=\boldsymbol{f}, \label{eq:1}
\end{align}
where $n$ is the number of degrees of freedom (DOFs), \(\boldsymbol{K}\) is the global stiffness matrix, and \(\boldsymbol{f}\) is the equivalent nodal force vector. In many complex applications, the FEM requires a sufficiently fine mesh to achieve the desired accuracy, which may result in a very large number of DOFs. This produces a large-scale matrix $\boldsymbol{K}$ and thus leads to substantial memory and computational cost for solving the system in Eq. \eqref{eq:1}. 

\begin{figure}
    \centering
    \begin{subfigure}{0.45\textwidth}
        \centering
        \includegraphics[width=0.8\textwidth]{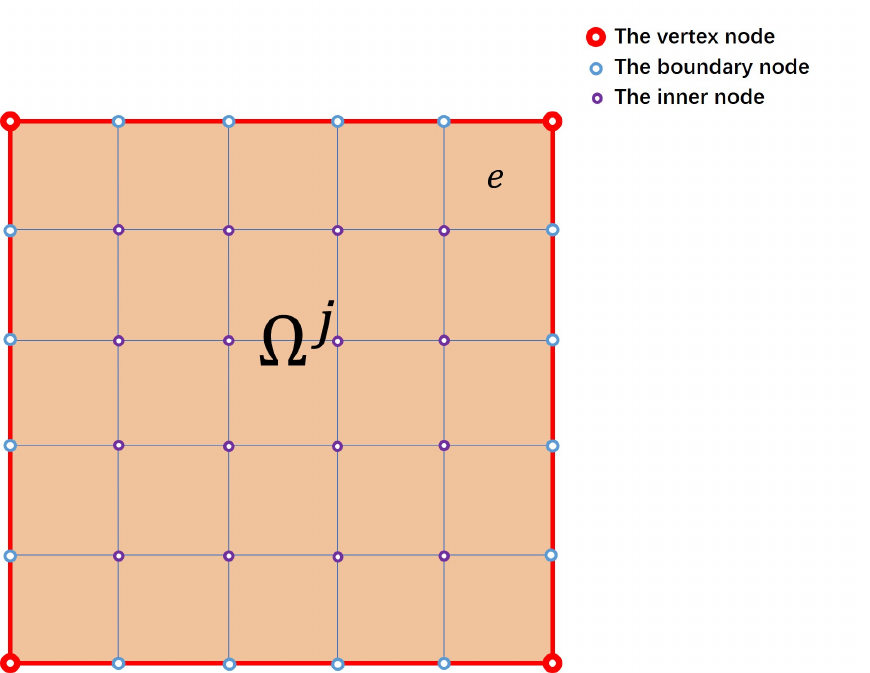}
        \caption{The rectangular substructure}
        \label{fig:sub:a}
    \end{subfigure}
    \hfill
    \centering
    \begin{subfigure}{0.45\textwidth}
        \centering
        \includegraphics[width=0.8\textwidth]{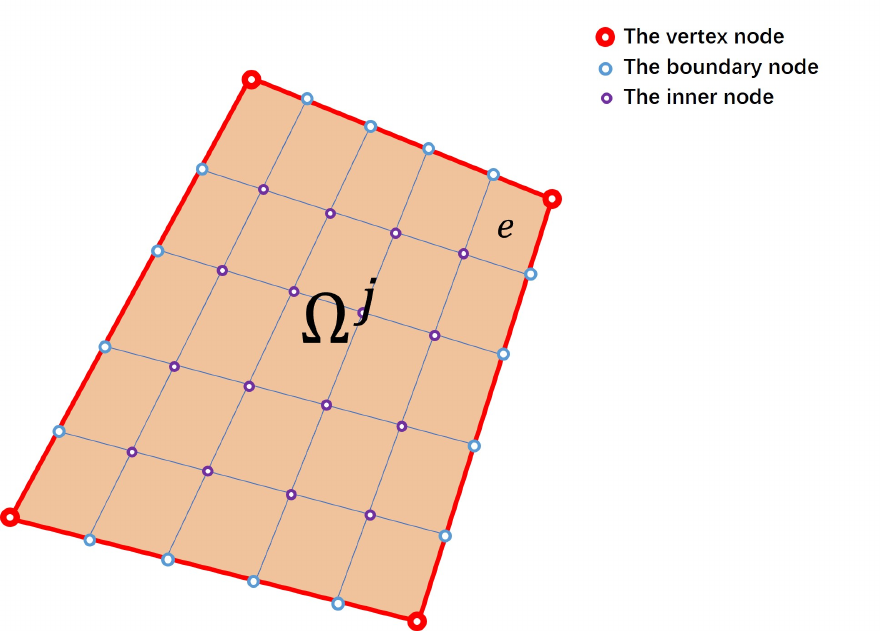}
        \caption{The quadrilateral substructure}
        \label{fig:sub:b}
    \end{subfigure}
    \caption{The substructure of two coarse elements}
    \label{fig:rect}
\end{figure}
To reduce the number of DOFs in the system in Eq. \eqref{eq:1}, a classical approach is the so-called substructure method \cite{23} based on the extended multiscale finite element method \cite{33}. In this method, several fine elements are aggregated to form a coarse element \(\Omega^j\) as shown in \autoref{fig:rect}, where \(\boldsymbol{e}\) denotes a fine element within the coarse element. The local finite element equations \(\boldsymbol{K}^j\boldsymbol{U}^j=\boldsymbol{F}^j\) for the coarse element \(\Omega^j\) can be written in block form as follows \cite{23,24}:
\begin{align*}
\boldsymbol{K}^j\boldsymbol{U}^j=\left(\begin{array}{cc}
		\boldsymbol{K}_{bb}^j&(\boldsymbol{K}^j_{ib})^T\\
		\boldsymbol{K}^j_{ib}&\boldsymbol{K}^j_{ii}
	\end{array}\right)\left(\begin{smallmatrix}
\boldsymbol{U}^j_{b} \\
\boldsymbol{U}^j_{i}
\end{smallmatrix}\right)=\left(\begin{smallmatrix}
\boldsymbol{F}^j_{b} \\
\boldsymbol{F}^j_{i}
\end{smallmatrix}\right),
\end{align*}
where \(\boldsymbol{U}^j_{b}\) denotes the boundary DOFs and \(\boldsymbol{U}^j_{i}\) denotes the internal DOFs in \(\Omega^j\), \(\boldsymbol{F}_b^j\) and \(\boldsymbol{F}_i^j\) denote the equivalent nodal forces of boundary DOFs and interior DOFs. Then \(\boldsymbol{U}^j_{i}\) can be represented by
\begin{align}
\boldsymbol{U}^j_{i}=(\boldsymbol{K}^j_{ii})^{-1}(\boldsymbol{F}^j_{i}-\boldsymbol{K}^j_{ib}\boldsymbol{U}^j_{b})\label{ep:2}.
\end{align}
Thus, the linear system involving only the boundary DOFs \(\boldsymbol{U}^j_b\) in a coarse element \(\Omega^j\) is
\begin{align}
\boldsymbol{K}^{j}_s\boldsymbol{U}^j_{b}=\boldsymbol{{F}}^{j}_s\label{ep:3},
\end{align}
where \(\boldsymbol{K}^{j}_s=\boldsymbol{K}^j_{bb}-(\boldsymbol{K}^j_{ib})^T(\boldsymbol{K}^j_{ii})^{-1}\boldsymbol{K}^j_{ib}\) is the condensed stiffness matrix, and similarly \(\boldsymbol{F}_{s}^j=\boldsymbol{F}_{b}^j-(\boldsymbol{K}^j_{ib})^T(\boldsymbol{K}^j_{ii})^{-1}\boldsymbol{F}^j_{i}\) is the equivalent condensed boundary nodal force vector. Once \(\boldsymbol{U}_{b}^j\) is determined, we can obtain \(\boldsymbol{U}_{i}^j\) by Eq. \eqref{ep:2}.

\begin{figure}[pos=h]
    \centering
    \begin{subfigure}{0.45\textwidth}
        \centering
        \includegraphics[width=0.8\textwidth]{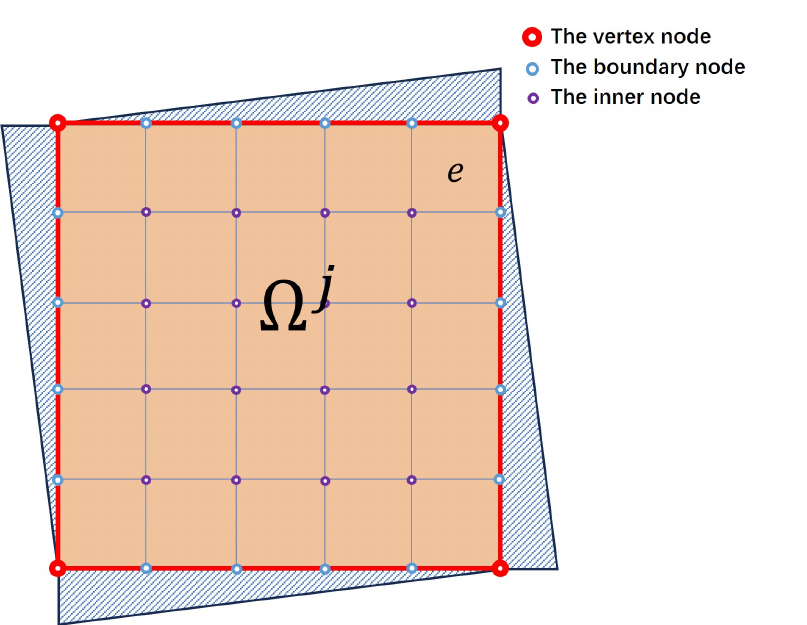}
        \caption{The rectangular substructure}
        \label{fig:sub:a}
    \end{subfigure}
    \hfill
    \centering
    \begin{subfigure}{0.45\textwidth}
        \centering
        \includegraphics[width=0.8\textwidth]{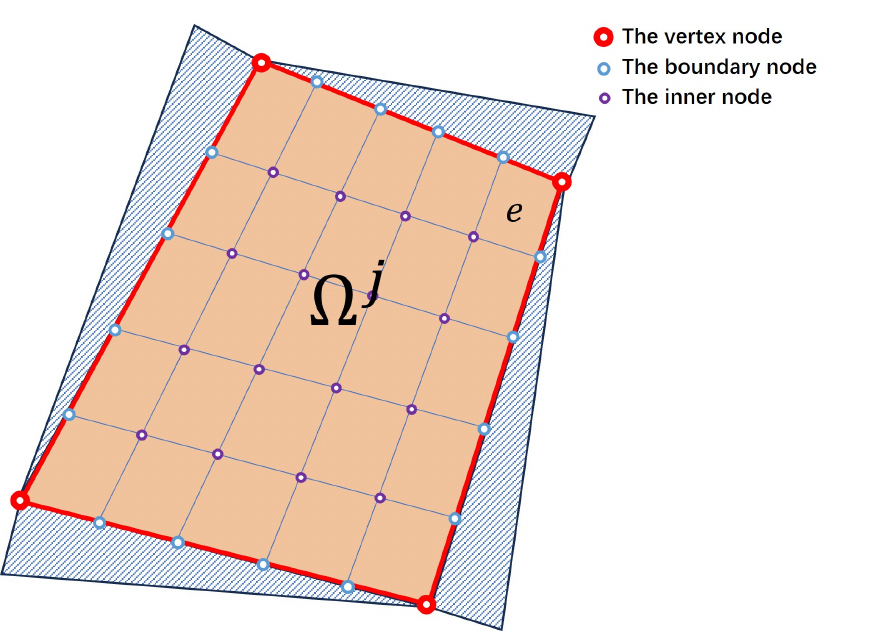}
        \caption{The quadrilateral substructure}
        \label{fig:sub:b}
    \end{subfigure}
    \caption{The substructure with linear boundary constraint}
    \label{fig:1}
\end{figure}

Moreover, to further reduce DOFs, a series of boundary processing methods are used, such as linear boundary constraints as shown in \autoref{fig:1}.
This means that the DOFs at the interior nodes on the boundary can be determined by the two end nodes. Thus, for a two-dimensional elasticity problem, a coarse element containing \(m\times m\) fine elements has $8m$ DOFs without boundary constraints, but only $8$ DOFs when linear boundary constraints are imposed. In order to obtain higher accuracy, we can also use the second-order or third-order polynomial displacement field constraints on the boundary of the coarse element.
\subsection{The numerical shape functions}\label{sect:2S2}
An alternative formulation of the substructure method is based on numerical shape functions \cite{22, 23}. We first consider the case without boundary constraints, which means that the number of DOFs of the condensed coarse element equals the number of boundary DOFs. As shown in \autoref{fig:1}, let \(\Omega^j\) represent a coarse element with \(m\times m\) fine elements, where $m$ is fixed. Denote by \(n_i\) the number of interior nodes and by \(n_b\) the number of boundary nodes, then we have \(n_i=(m-1)^2\) and \(n_b=4m\). We use the bilinear shape functions on each fine element and denote by \(\boldsymbol{N_e}(\boldsymbol{x})\) the row vector composed of all the piecewise bilinear basis functions. The displacement functions on the coarse element \(\Omega^j\) can be represented as \(\boldsymbol{u}^j(\boldsymbol{x})=\boldsymbol{N_e}(\boldsymbol{x})\boldsymbol{U}^j\). By Eq. \eqref{ep:2}, if the nodal force \(\boldsymbol{f}^j\) vanishes in \(\Omega^j\), we define the transition matrix as 
\begin{align}\boldsymbol{N}^j=\left(\begin{smallmatrix}
\boldsymbol{I} \\
\boldsymbol{N}^j_s
\end{smallmatrix}\right)=\left(\begin{smallmatrix}
\boldsymbol{I} \\
-(\boldsymbol{K}_{ii}^j)^{-1}\boldsymbol{K}_{ib}^j
\end{smallmatrix}\right), \label{ep:Nj}
\end{align}
such that \(\boldsymbol{U}^j=\boldsymbol{N}^j\boldsymbol{U}^j_b\). Clearly, \(\boldsymbol{I}\) is a \(2n_b\times 2n_b\) identity matrix and \(\boldsymbol{N}_s^j\) is a \(2n_i\times2n_b\) matrix.  Hence, the displacement functions can be determined by the boundary nodal DOFs as \(\boldsymbol{U}^j_i=\boldsymbol{N}^j_s\boldsymbol{U}^j_b\). We find that \(\boldsymbol{N}_{s}^j\)  has the following structure:
\begin{align*}
    \boldsymbol{N}_{s}^j=\left(
	\begin{array}{ccccccc}
		N^{j1}_{1, xx}&N^{j1}_{1,yx}&N^{j1}_{2,xx}&N^{j1}_{2, yx}&\ldots&N^{j1}_{n_b, xx}&N^{j1}_{n_b, yx}\\
		N^{j1}_{1, xy}&N^{j1}_{1,yy}&N^{j1}_{2,xy}&N^{j1}_{2, yy}&\ldots&N^{j1}_{n_b, xy}&N^{j1}_{n_b, yy}\\
		N^{j2}_{1, xx}&N^{j2}_{1,yx}&N^{j2}_{2,xx}&N^{j2}_{2, yx}&\ldots&N^{j2}_{n_b, xx}&N^{j2}_{n_b, yx}\\
		N^{j2}_{1, xy}&N^{j2}_{1,yy}&N^{j2}_{2,xy}&N^{j2}_{2, yy}&\ldots&N^{j2}_{n_b, xy}&N^{j2}_{n_b, yy}\\
		&&&&\ldots\ldots&&\\
		N^{jn_i}_{1, xx}&N^{jn_i}_{1, yx}&N^{jn_i}_{2, xx}&N^{jn_i}_{2, yx}&\ldots&N^{jn_i}_{n_b, xx}&N^{jn_i}_{n_b, yx}\\
		N^{jn_i}_{1, xy}&N^{jn_i}_{1, yy}&N^{jn_i}_{2, xy}&N^{jn_i}_{2, yy}&\ldots&N^{jn_i}_{n_b, xy}&N^{jn_i}_{n_b, yy}\\
	\end{array}
\right).
\end{align*}
The displacement at a point \(\boldsymbol{x_e}\) in the fine element \(\boldsymbol{e}\) can be represented as \(\boldsymbol{u}^j(\boldsymbol{x_e})=\boldsymbol{N_e}(\boldsymbol{x_e})\boldsymbol{N_e}^j\boldsymbol{U}_b^j\), where \(\boldsymbol{N_e}^j\) is the submatrix of \(\boldsymbol{N}^j\) which consists of the rows corresponding to those DOFs of fine element \(\boldsymbol{e}\). Let \(\boldsymbol{K}^j_e\) denote the contribution of fine element \(\boldsymbol{e}\) to the condensed stiffness matrix. The condensed stiffness matrix and condensed external load  can be written as follows:\\
\begin{align}
\boldsymbol{K}_{s}^j&=\bigwedge_{e}\boldsymbol{K}^j_{e}=\bigwedge_{e}(\boldsymbol{N_e}^j)^T\boldsymbol{B_{e}^T}\boldsymbol{D}\boldsymbol{B_{e}}\boldsymbol{N_e}^j={(\boldsymbol{N}^{j})}^T\boldsymbol{K}_j\boldsymbol{N}^j\label{ep:9},\\
\boldsymbol{F}_{s}^j&=\bigwedge_{e}(\boldsymbol{N_e}^j)^T\boldsymbol{N_e}(\boldsymbol{x})^T\boldsymbol{f_e}=(\boldsymbol{N}^{j})^T\boldsymbol{F}^j,
\end{align}
where \(\bigwedge\) denotes the assembly operator, \(\boldsymbol{B_e}\) is the condensed strain matrix of fine element \(\boldsymbol{e}\), and \(\boldsymbol{f_e}\) is the external load vector. Clearly, the stiffness matrix in Eq. \eqref{ep:9} is equal to that in Eq. \eqref{ep:3} \cite{22,23}. 

Although both methods yield the same condensed matrix, the second approach better illustrates the relationship between the condensed shape functions and the original shape functions, which is important for subsequent analysis. In the following, we use the formulation in Eq. \eqref{ep:9}, and define the numerical shape functions by
\begin{align}
\boldsymbol{N}_l^{j}(\boldsymbol{x})=\boldsymbol{N_e}(\boldsymbol{x})\boldsymbol{N}_{, l}^j, l=1,2,\dots,2n_b,\label{ep:NBF}
\end{align}
where \(\boldsymbol{N}_{, l}^j\) denotes the $l-$th column of $\boldsymbol{N}^j$. The following properties then hold:
\begin{property}Properties of the numerical shape functions:\\
\label{thm:P2p3}
(1) For any $l$, \(\boldsymbol{N}_l^{j}(\boldsymbol{x})=\boldsymbol{N_e}(\boldsymbol{x})\boldsymbol{N}_{, l}^j\) is a continuous piecewise bilinear function and is $C^0$ continuous on \(\Omega^j\). \\
(2) The space spanned by \(\boldsymbol{N}_l^{j}(\boldsymbol{x})\), \(l=1, 2, . . . . , 2n_b\), denoted by $U_h$ is a subspace of $V_h$ spanned by \(\boldsymbol{N_e}(\boldsymbol{x})\). \\
(3) \(\{\boldsymbol{N}_l^{j}(\boldsymbol{x})\}\) forms a set of nodal interpolation basis functions on \(\Omega^j\).
\end{property}
Once the fine mesh is fixed, the properties of the numerical shape functions \(\{\boldsymbol{N}_l^{j}(\boldsymbol{x})\}\) are completely determined by the properties of \(\boldsymbol{N}^j\). Therefore, the study of \(\{\boldsymbol{N}_l^{j}(\boldsymbol{x})\}\) is reduced to the study of \(\boldsymbol{N}^j\). We introduce the following three definitions.
\begin{definition}[\textbf{Partition of unity}]
\label{thm:D2f1}
A \(2(n_i+n_b)\times2n_b\) matrix \(\boldsymbol{N}^j\) is said to satisfy the row-wise partition-of-unity property if the following equations hold:
\begin{align}
\boldsymbol{N}\left[\begin{smallmatrix}
1 & 0 & 1 & 0 & . . . & 1 & 0 \\
0 & 1 & 0 & 1 & . . . & 0 & 1
\end{smallmatrix}\right]_{2\times2n_b}^T=\left[\begin{smallmatrix}
1 & 0 & 1 & 0 & . . . & 1 & 0 \\
0 & 1 & 0 & 1 & . . . & 0 & 1
\end{smallmatrix}\right]_{2\times2(n_i+n_b)}^T.
\end{align}
\end{definition}
\begin{definition}[\textbf{rigid-body displacement invariance}]
\label{thm:D2f2}
Let \(\boldsymbol{X}_b=[x_{b, 1}, x_{b, 2}, x_{b, 3}, . . . x_{b, n_b}]^T\) be the x-coordinates of boundary nodes, \(\boldsymbol{Y}_b=[y_{b, 1}, y_{b, 2}, y_{b, 3}, . . . y_{b, n_b}]^T\) be the y-coordinates. Similarly, we let \(\boldsymbol{X}=[x_{i, 1}, x_{i, 2}, x_{i, 3}, . . . x_{i, n_i+n_b}]^T\) and \(\boldsymbol{Y}=[y_{i, 1}, y_{i, 2}, y_{i, 3}, . . . y_{i, n_i+n_b}]^T\) be the x-coordinates and y-coordinates of the boundary and internal nodes. We say that the matrix \(\boldsymbol{N}^j\) satisfies rigid-body displacement invariance if the following equations hold:
\begin{align}
\boldsymbol{N}^j\left[\begin{smallmatrix}
1 & 0 & 1 & 0 & . . . & 1 & 0 \\
0 & 1 & 0 & 1 & . . . & 0 & 1
\end{smallmatrix}\right]_{2\times2n_b}^T&=\left[\begin{smallmatrix}
1 & 0 & 1 & 0 & . . . & 1 & 0 \\
0 & 1 & 0 & 1 & . . . & 0 & 1
\end{smallmatrix}\right]_{2\times2(n_i+n_b)}^T,\\
\boldsymbol{N}^j[{y}_{b,1}, -{x}_{b,1}, {y}_{b,2}, -{x}_{b,2},\dots, {y}_{b,n_b}, -{x}_{b,n_b}]^T&=[y_{i,1},-x_{i,1}, y_{i,2}, -x_{i,2},\dots, y_{i,n_i+n_b}, -x_{n_i+n_b}]^T.\nonumber
\end{align}
\end{definition}
\begin{definition}[\textbf{Linear completeness}]
\label{thm:D2f3}
Let \(\boldsymbol{X}_b=\left[\begin{smallmatrix}
x_{b, 1} & 0 & x_{b, 2} & 0 & . . . & x_{b, n_b} & 0 \\
0 & x_{b, 1} & 0 & x_{b, 2} & . . . & 0 & x_{b, n_b}
\end{smallmatrix}\right]^T\) be the x-coordinates of boundary nodes, \(\boldsymbol{Y}_b=\left[\begin{smallmatrix}
y_{b, 1} & 0 & y_{b, 2} & 0 & . . . & y_{b, n_b} & 0 \\
0 & y_{b, 1} & 0 & y_{b, 2} & . . . & 0 & y_{b, n_b}
\end{smallmatrix}\right]^T\) be the y-coordinates. Similarly, for internal nodes, let \(\boldsymbol{X}_i=\left[\begin{smallmatrix}
x_{i, 1} & 0 & x_{i, 2} & 0 & . . . & x_{i, n_i} & 0 \\
0 & x_{i, 1} & 0 & x_{i, 2} & . . . & 0 & x_{i, n_i}
\end{smallmatrix}\right]^T\) and \(\boldsymbol{Y}_i=\left[\begin{smallmatrix}
y_{i, 1} & 0 & y_{i, 2} & 0 & . . . & y_{i, n_i} & 0 \\
0 & y_{i, 1} & 0 & y_{i, 2} & . . . & 0 & y_{i, n_i}
\end{smallmatrix}\right]^T\) be their x-coordinates and y-coordinates, and let \(\boldsymbol{1}_{2\times2n_b}=\left[\begin{smallmatrix}
1 & 0 & 1 & 0 & . . . & 1 & 0 \\
0 & 1 & 0 & 1 & . . . & 0 & 1
\end{smallmatrix}\right]_{2\times2n_b}^T\). We say that the matrix \(\boldsymbol{N}^j\) satisfies linear completeness if the following equations hold:
\begin{align}
\boldsymbol{N}^j\boldsymbol{1}_{2\times2n_b}&=\boldsymbol{1}_{2\times2(n_b+n_i)}\nonumber,\\
\boldsymbol{N}^j\boldsymbol{X}_b&=\left[\begin{smallmatrix}\boldsymbol{X}_b\\
\boldsymbol{X}_i\end{smallmatrix}\right]\label{ep:Linearcom},\\
\boldsymbol{N}^j\boldsymbol{Y}_b&=\left[\begin{smallmatrix}\boldsymbol{Y}_b\\
\boldsymbol{Y}_i\end{smallmatrix}\right]\nonumber.
\end{align}
\end{definition}
It is clear that linear completeness implies rigid-body displacement invariance, and rigid-body displacement invariance implies the row-wise partition-of-unity property, but the converses do not hold.

Obviously, the numerical shape functions have the same properties as the transition matrix.
\begin{theorem}
\label{thm:T2m1}
Let \(\boldsymbol{N}^j\) be the transition matrix and let \(\boldsymbol{N_e}(\boldsymbol{x})\boldsymbol{N}^j\) denote the corresponding numerical shape functions. Then:
\begin{enumerate}
  \item[(i)] if \(\boldsymbol{N}^j\) satisfies the row-wise partition-of-unity property, then \(\boldsymbol{N_e}(\boldsymbol{x})\boldsymbol{N}^j\) satisfies the partition-of-unity property, i.e.,
	\begin{align*}
	\boldsymbol{N_e}(\boldsymbol{x})\boldsymbol{N}^j\left[\begin{smallmatrix}
1&0&1&0&\dots&1&0 \\
0&1&0&1&\dots&0&1
\end{smallmatrix}\right]_{2\times2n_b}^T=\left[\begin{smallmatrix}
1&0 \\
0&1
\end{smallmatrix}\right],
	\end{align*}
\item[(ii)] if \(\boldsymbol{N}^j\) has rigid-body displacement invariance, then \(\boldsymbol{N_e}(\boldsymbol{x})\boldsymbol{N}^j\) satisfies rigid-body displacement invariance, i.e.,
	\begin{align*}
	&\boldsymbol{N_e}(\boldsymbol{x})\boldsymbol{N}^j\left[\begin{smallmatrix}
1&0&1&0&\dots&1&0 \\
0&1&0&1&\dots&0&1
\end{smallmatrix}\right]_{2\times2n_b}^T=\left[\begin{smallmatrix}
1&0 \\
0&1
\end{smallmatrix}\right],\\
&\boldsymbol{N_e}(\boldsymbol{x})\boldsymbol{N}^j\left[\begin{smallmatrix}
y_1&-x_1&y_2&-x_2&\dots&y_{n_b}&-x_{n_b}
\end{smallmatrix}\right]^T=\left[\begin{smallmatrix}
y \\
-x
\end{smallmatrix}\right],
	\end{align*}
\item[(iii)] if \(\boldsymbol{N}^j\) has linear completeness, then \(\boldsymbol{N_e}(\boldsymbol{x})\boldsymbol{N}^j\) satisfies linear completeness, i.e.,
	\begin{align*}
	&\boldsymbol{N_e}(\boldsymbol{x})\boldsymbol{N}^j\left[\begin{smallmatrix}
1&0&1&0&\dots&1&0 \\
0&1&0&1&\dots&0&1
\end{smallmatrix}\right]_{2\times2n_b}^T=\left[\begin{smallmatrix}
1&0 \\
0&1
\end{smallmatrix}\right],\\
&\boldsymbol{N_e}(\boldsymbol{x})\boldsymbol{N}^j\left[\begin{smallmatrix}
x_1&0&x_2&0&\dots&x_{n_b}&0 \\
0&x_1&0&x_2&\dots&0&x_{n_b}
\end{smallmatrix}\right]^T=\left[\begin{smallmatrix}
x&0 \\
0&x
\end{smallmatrix}\right],\\
&\boldsymbol{N_e}(\boldsymbol{x})\boldsymbol{N}^j\left[\begin{smallmatrix}
y_1&0&y_2&0&\dots&y_{n_b}&0 \\
0&y_1&0&y_2&\dots&0&y_{n_b}
\end{smallmatrix}\right]^T=\left[\begin{smallmatrix}
y&0 \\
0&y
\end{smallmatrix}\right].
	\end{align*}
\end{enumerate}
\end{theorem}
We also note that \(\boldsymbol{N}_s^j\) is independent of the coarse-element size \(h\).
\begin{theorem}
\label{thm:T2m2}
The entries of the transition matrix \(\boldsymbol{N}_s^j\) are independent of the size of the coarse element.
\end{theorem}
\begin{proof}
Recall that \(\boldsymbol{N}_s^j=-{(\boldsymbol{K}_{ii}^j)}^{-1}\boldsymbol{K}_{ib}^j\). For a second-order linear elasticity problem, the local stiffness matrix \(\boldsymbol{K}^j\) is assembled from the stiffness matrices of the fine elements. The stiffness matrix associated with a fine element \(\boldsymbol{e}\) is given by:
\begin{align*}
\boldsymbol{K_e}=\int_{\boldsymbol{e}} \boldsymbol{B_e}^T(\boldsymbol{x})\boldsymbol{D}(\boldsymbol{x})\boldsymbol{B_e}(\boldsymbol{x})\boldsymbol{dx},
\end{align*}
where \(\boldsymbol{B_e}(\boldsymbol{x})\) is the strain-displacement matrix and \(\boldsymbol{D}(\boldsymbol{x})\) is the constitutive matrix. We assume that the material coefficient distribution remains unchanged in reference coordinates under geometric scaling. When the coarse element is scaled by its size \(h\), it is clear that the fine elements are scaled by the same factor. For two-dimensional problems, we have
\begin{align*}
\bar{\boldsymbol{K_e}}=\int_{\bar{\boldsymbol{e}}} \boldsymbol{B_e}^T(\frac{\boldsymbol{x}}{h})\boldsymbol{D}(\frac{\boldsymbol{x}}{h})\boldsymbol{B_e}(\frac{\boldsymbol{x}}{h})\boldsymbol{dx},
\end{align*}
where \(\bar{\boldsymbol{K}_e}\) denotes the stiffness matrix of the scaled element and \(\bar{\boldsymbol{e}}\) is the scaled fine element. We set \(\boldsymbol{y}=\frac{\boldsymbol{x}}{h}\) which implies
\begin{align*}
\bar{\boldsymbol{K_e}}&=\int_{\boldsymbol{e}}\frac{1}{h}\boldsymbol{B_e}^T(\boldsymbol{y})\boldsymbol{D}(\boldsymbol{y})\frac{1}{h}\boldsymbol{B_e}(\boldsymbol{y})h^2\boldsymbol{dy}\\
&=\int_{\boldsymbol{e}}\boldsymbol{B_e}^T(\boldsymbol{y})\boldsymbol{D}(\boldsymbol{y})\boldsymbol{B_e}(\boldsymbol{y})\boldsymbol{dy},\\
&=\boldsymbol{K_e}.
\end{align*}
Then we obtain \(\bar{\boldsymbol{K}}^j=\boldsymbol{K}^j\), it implies \({(\bar{\boldsymbol{K}}_{ii}^j)}^{-1}{\bar{\boldsymbol{K}}_{ib}^j}={({\boldsymbol{K}}_{ii}^j)}^{-1}{\boldsymbol{K}}_{ib}^j\). For a three-dimensional problem, we have
\begin{align*}
\bar{\boldsymbol{K_e}}&=\int_{{e}}\frac{1}{h}\boldsymbol{B_e}^T(\boldsymbol{y})\boldsymbol{D}(\boldsymbol{y})\frac{1}{h}\boldsymbol{B_e}(\boldsymbol{y})h^3\boldsymbol{dy}\\
&=\int_{e}h\boldsymbol{B_e}^T(\boldsymbol{y})\boldsymbol{D}(\boldsymbol{y})\boldsymbol{B_e}(\boldsymbol{y})\boldsymbol{dy}\\
&=h\boldsymbol{K_e},
\end{align*}
then \(\bar{\boldsymbol{K}}^j=h\boldsymbol{K}^j\), it implies:
\begin{align*}
{(\bar{\boldsymbol{K}}_{ii}^j)}^{-1}\bar{\boldsymbol{K}}_{ib}^j=h^{-1}h{(\boldsymbol{K}_{ii}^j)}^{-1}\boldsymbol{K}_{ib}^j={(\boldsymbol{K}_{ii}^j)}^{-1}\boldsymbol{K}_{ib}^j.
\end{align*}
This completes the proof.
\end{proof}
\begin{remark}Although the derivation above assumes that the equivalent nodal force inside the coarse element vanishes, the same equivalence between the two formulations also holds when the equivalent nodal force is nonzero. This can be verified directly by replacing \((\boldsymbol{K}^j_{ii})^{-1}\boldsymbol{K}^j_{ib}\) with \(\boldsymbol{N}^j_s\). When the body force is nonzero, the displacements at the internal nodes must be obtained via \(\boldsymbol{U}^j_i=\boldsymbol{N}^j_s\boldsymbol{U}^j_b+{(\boldsymbol{K}_{ii}^j)}^{-1}F^j_i\).
\end{remark}
\subsection{Problem-Independent Machine Learning (PIML)}
Although the substructure method indeed reduces DOFs, the computational cost of computing $\boldsymbol{N}^j$ remains high. In some problems that require repeated stiffness matrix calculations and in which material distributions or mesh geometries change repeatedly, substantial computational time is required. From Eq. \eqref{ep:2} and Eq. \eqref{ep:9}, the transition matrix is defined by the local stiffness matrix of the coarse element. Once the transition matrix has been obtained, the condensed stiffness matrix can be obtained from the transition matrix. To obtain the transition matrix efficiently, the PIML approach employs machine learning to predict the transition matrix \cite{22,23,24,30}. The idea of PIML is to construct a mapping from the material distribution of each coarse element to $\boldsymbol{N^j}$. Once the mapping has been constructed, $\boldsymbol{N^j}$ of any coarse element can be directly obtained.

PIML incorporates machine learning (ML) into the computation of $\boldsymbol{N^j}$. This is based on three conclusions \cite{22,23}: (1) There is a unique relationship between $\boldsymbol{N^j}$ and the material distribution of the coarse element. This property ensures the well-posedness and effectiveness of the learning problem. (2) $\boldsymbol{N^j}$ is problem-independent, which means that it is independent of the global boundary conditions for PDEs. Once the mapping has been obtained, the model can be applied to problems with various global boundary conditions for PDEs. (3) Training samples are easily obtained, and the neural network has sufficient approximation capacity to approximate the target mapping. These properties provide the basis for constructing an effective learning model, although the prediction accuracy also depends on the training procedure and the network architecture. 

\begin{figure}
    \centering
    \includegraphics[width=0.55\textwidth]{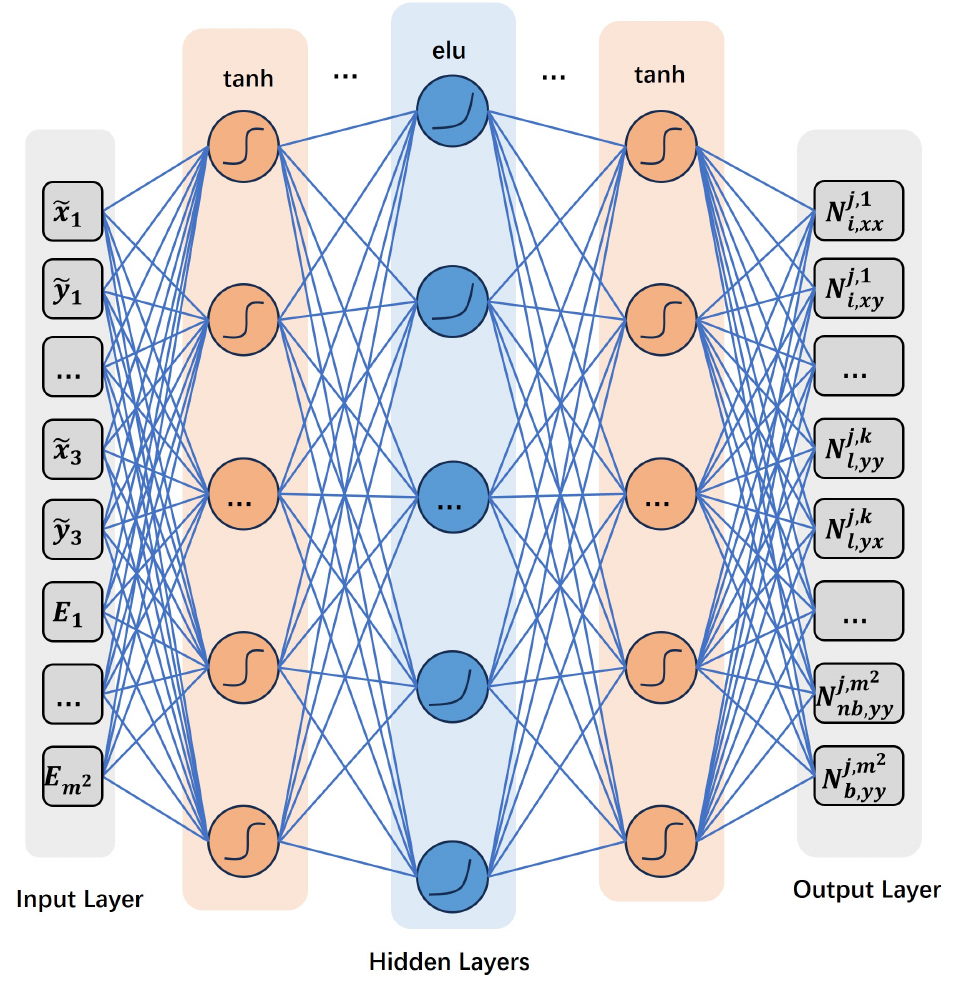}
    \caption{Diagram of the deep neural network \upcite{22,23}}\label{fig:3}
\end{figure}

We first assume that the structure of the neural network is a multilayer perceptron (MLP) \cite{34}, the most common neural network structure. In our research, the inputs of the MLP are the distribution of the material coefficient and the x- and y-coordinates of the vertices of the coarse element, whereas the outputs are the values of the numerical shape functions at the fine mesh nodes, or more directly, \(\boldsymbol{N^j_s}\). The architecture of the network is shown in \autoref{fig:3}.

Thus, PIML essentially acts as a matrix fitter. Although other network architectures may also be used in PIML, their common objective is to fit the mapping between \(\boldsymbol{N}^j_s\) and the material distribution of the coarse element. When the neural network has been trained, the predicted matrix denoted by \(\widetilde{\boldsymbol{N}}_s^j\) cannot be used directly in finite element computation, because it may not satisfy consistency at fine- and coarse-mesh nodes under rigid-body motion. To achieve this, the following rigid-body displacement invariance constraints on \(\widetilde{\boldsymbol{N}}_s^j\) must be satisfied:
\begin{align}
[1, 0, 1, 0, . . . . , 1, 0]^T&=\widetilde{\boldsymbol{N_s}}^j[1, 0, 1, 0, . . . . , 1, 0]^T,\nonumber\\
[0, 1, 0, 1, . . . . , 0, 1]^T&=\widetilde{\boldsymbol{N_s}}^j[0, 1, 0, 1, . . . . , 0, 1]^T,\label{ep:Rigdre}\\
[y_{i,1}, -x_{i,1}, y_{i,2}, -x_{i,2}, . . . . , y_{i,n_i}, -x_{i,n_i}]^T&=\widetilde{\boldsymbol{N_s}}^j[{y}_{b,1}, -{x}_{b,1}, {y}_{b,2}, -{x}_{b,2}, . . . . , {y}_{b,n_b}, -{x}_{b,n_b}]^T,\nonumber
\end{align}
To satisfy these constraints, we only need to learn \(2n_b-3\) columns of \(\boldsymbol{N^j_s}\) (in a 3D problem the number would be \(3n_b-6\)), and the remaining three columns (six columns for a 3D problem) can be determined by Eq. \eqref{ep:Rigdre}.\\
According to \(\boldsymbol{N}^j=\left(\begin{smallmatrix}
\boldsymbol{I} \\
\boldsymbol{N^j_s}
\end{smallmatrix}\right)\), for the case without boundary constraints, the learned transition matrix is \(\widetilde{\boldsymbol{N}}^j=\left(\begin{smallmatrix}
\boldsymbol{I} \\
\boldsymbol{\widetilde{\boldsymbol{N}}^j_s}
\end{smallmatrix}\right)\). Similar to \cite{22,23}, we use a submatrix \(L\) to represent the linear boundary constraints on the coarse element \(\Omega^j\). Then \(\widetilde{\boldsymbol{N}}^j=\left(\begin{smallmatrix}
\boldsymbol{L} \\
\boldsymbol{\widetilde{\boldsymbol{N}}^j_s}
\end{smallmatrix}\right)\). Since matrix \(L\) depends only on the boundary constraints of the coarse element and is identical for all coarse elements, it does not need to be learned. Denote \(\delta\boldsymbol{N}^j=\boldsymbol{N}^j-\widetilde{\boldsymbol{N}}^j\). From the neural network matrix learning algorithm, \(\lvert\delta\boldsymbol{N}^j\lvert_\infty\) can be used to represent an upper bound on the generalization error for matrix learning.
\begin{remark}For rectangular elements, the inputs contain only the material distribution over the fine elements. For quadrilateral elements, owing to translational invariance, one vertex can be placed at the origin; therefore, the coordinates of the remaining three vertices are additionally required as network inputs, as shown in \autoref{fig:3}.
\end{remark}
\begin{remark}
For the three-dimensional case, the above definitions and theorems can be extended similarly by adding the \(z-\)coordinates.
\end{remark}
\section{Error analysis}\label{Sect. 3}
In this section, we analyze the errors at the coarse-element nodes using standard finite element error estimates and the properties of the numerical shape functions. We first establish the boundedness of the interpolation operator defined by the numerical shape functions with boundary DOFs. \\
\subsection{Interpolation operator on the reference element \(\widehat{\Omega}^j\)}
\begin{figure}
    \centering
	\centering
    \includegraphics[width=0.6\textwidth]{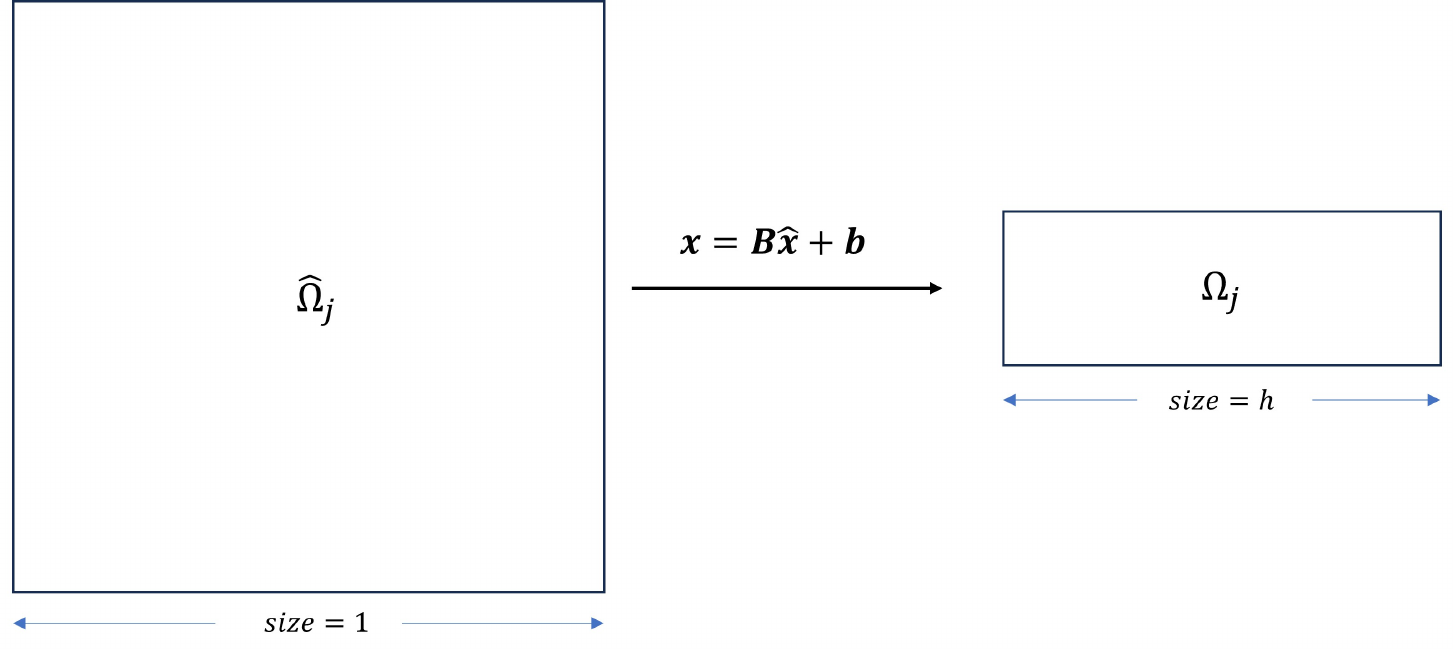}
    \caption{The affine transformation}\label{fig:affine}
\end{figure}

We consider the boundedness of the interpolation operator on the reference rectangular element \(\widehat{\Omega}^j\), which is obtained through an affine mapping \(\mathfrak{B}^j\):
\begin{align*}
\widehat{\Omega}^j\rightarrow\Omega^j:\boldsymbol{x}=\boldsymbol{B}^j\hat{\boldsymbol{x}}+\boldsymbol{b},
\end{align*}
where \(\boldsymbol{B}^j\) is an invertible matrix, the affine transformation maps the reference element \(\widehat{\Omega}^j\) with unit size onto the coarse element \(\Omega^j\) of size \(h\) as shown in \autoref{fig:affine}. Because the properties of this affine mapping \(\mathfrak{B}^j\) depend on those of \(\boldsymbol{B}^j\), we will use \(\boldsymbol{B}^j\) to denote the operator \(\mathfrak{B}^j\).
Without loss of generality, let \(\boldsymbol{N}^j\) be an \(n_1\times n_2\) matrix, where \(n_1=2(n_i+n_b),\ n_2=2n_b\) for two-dimensional cases and \(n_1=3(n_i+n_b),\ n_2=3n_b\) for three-dimensional cases. Let \(\widehat{\boldsymbol{N}}_{\boldsymbol{e}}(\boldsymbol{x})\) denote the vector of fine-element shape functions on the reference element \(\widehat{\Omega}^j\)with unit size. Then \(\widehat{\boldsymbol{N}}_{\boldsymbol{e}}(\boldsymbol{x})\boldsymbol{N}^j\) defines an interpolation operator \(\widehat{\Pi}:H^2\rightarrow H^1\):
\begin{align}
\widehat{\Pi}{f}=\sum_{l=1}^{n_2}\widehat{\boldsymbol{N}}_{l}^j(\boldsymbol{x})f(\boldsymbol{x}_l)=\sum_{l=1}^{n_2}\widehat{\boldsymbol{N}}_{\boldsymbol{e}}(\boldsymbol{x})\boldsymbol{N}^j_{,l}f(\boldsymbol{x_l})=\widehat{\boldsymbol{N}}_{\boldsymbol{e}}(\boldsymbol{x})\sum_{l=1}^{n_2}\boldsymbol{N}^j_{,l}f(\boldsymbol{x_l}),\label{ep:referop}
\end{align}
where \(\boldsymbol{x}_l\) denotes the nodal points associated with the boundary DOFs. Then we have the following lemma.
\begin{lemma} The  $H^1$-seminorm and $L^2$-norm of shape functions \(\widehat{N}_{l}^j(\boldsymbol{x})=\widehat{\boldsymbol{N}}_{\boldsymbol{e}}(\boldsymbol{x})\boldsymbol{N}^j_{,l}, l=1,2,\dots,n_2\) on the reference element \(\widehat{\Omega}^j\) are uniformly bounded.
\label{thm:bound1}
\end{lemma}
\begin{proof}
Note that \(\widehat{\boldsymbol{N}}_{\boldsymbol{e}}(\boldsymbol{x})\) consists of the bilinear basis functions $Q_4$ in two dimensions or the trilinear basis functions $H_{8}$ in three dimensions. Since the reference element \(\widehat{\Omega}^j\) has unit diameter, there exist constants $C_1$ and $C_2$ such that the following inequalities hold :
\begin{align}
\lvert \widehat{N}_l^j(\boldsymbol{x})\lvert_{H^1(\widehat{\Omega})}\leq C_1\lvert\boldsymbol{N}^j\lvert_{\infty}, l=1,2,\dots,n_2,\label{ep:14}\\
\lvert \widehat{N}_l^j(\boldsymbol{x})\lvert_{L^2(\widehat{\Omega})}\leq C_2\lvert\boldsymbol{N}^j\lvert_{\infty}, l=1,2,\dots,n_2,\label{ep:15}
\end{align}
where \(\lvert \boldsymbol{{N}}^j \rvert_{\infty}\) denotes the maximum absolute value of the entries of \(\boldsymbol{N}^j\).
\end{proof}

We define \(\widehat{\widetilde{\Pi}}\) as the interpolation operator obtained by PIML on the reference coarse element \(\widehat{\Omega}^j\),i.e.,
\begin{align}
\widehat{\widetilde{\Pi}}=\sum_{l=1}^{n_2}\widehat{\widetilde{\boldsymbol{N}}}_{l}^j(\boldsymbol{x})f(\boldsymbol{x}_l)=\sum_{l=1}^{n_2}\widehat{\boldsymbol{N}}_{\boldsymbol{e}}(\boldsymbol{x})\widetilde{\boldsymbol{N}}^j_{,l}f(\boldsymbol{x_l})=\widehat{\boldsymbol{N}}_{\boldsymbol{e}}(\boldsymbol{x})\sum_{l=1}^{n_2}\widetilde{\boldsymbol{N}}^j_{,l}f(\boldsymbol{x_l}),\label{ep:referoppiml}
\end{align}
Then we have the difference operator between \(\widehat{\Pi}\) in Eq.\eqref{ep:referop} and \(\widehat{\widetilde{\Pi}}\) in Eq.\eqref{ep:referoppiml} on the reference coarse element \(\widehat{\Omega}^j\):
\begin{align*}
\widehat{\delta{\Pi}}=\widehat{\Pi}-\widehat{\widetilde{\Pi}}=\widehat{\boldsymbol{N}}_{\boldsymbol{e}}(\boldsymbol{x})(\boldsymbol{N}^j-\widetilde{\boldsymbol{N}}^j)=\widehat{\boldsymbol{N}}_{\boldsymbol{e}}(\boldsymbol{x})\delta{\boldsymbol{N}^j}.
\end{align*}
We use \(\varepsilon=\lvert \delta{\boldsymbol{N}}^j \rvert_{\infty}\) to denote the generalization error for \({\boldsymbol{N}}^j\), then, \(\widehat{\delta{\Pi}}\) can be bounded in terms of the generalization error of the learned transition matrix, as established below.
\begin{lemma} Suppose \(\delta\boldsymbol{N^j}\) is bounded, i.e., \(\varepsilon=\lvert\delta\boldsymbol{N^j}\lvert_{\infty}\), and \(\boldsymbol{x}\in\mathbb{R}^d,d=2,3\).  Then the operator \(\widehat{\delta\Pi}\) is bounded,
\begin{align}
\parallel\widehat{\delta\Pi} u\parallel_{H^1(\widehat{\Omega}^j)}\leq\varepsilon C_3\parallel u\parallel_{H^2(\widehat{\Omega}^j)}.
\end{align}
\label{thm:bound2}
\end{lemma}
\begin{proof}
Let \(u(\boldsymbol{x})\in H^2\), and let \(\delta\boldsymbol{N}^j_{,l}\) denote the $l-$th column of \(\delta\boldsymbol{N}^j\) as in Eq. \eqref{ep:NBF}. We obtain
\begin{align*}
\parallel\widehat{\delta\Pi}u\parallel_{H^1(\widehat{\Omega}^j)}&=\parallel\widehat{\boldsymbol{N}}_{\boldsymbol{e}}(\boldsymbol{x})\sum_{l}^{n_2}\delta\boldsymbol{N}^j_{, l}u(\boldsymbol{x}_l)\parallel_{H^1(\widehat{\Omega}^j)}\\
&\leq\sum_{k}^{n_1}\big(\parallel \widehat{N}^k(\boldsymbol{x})\sum_{l=1}^{n_2}\delta N^j_{k, l}u(\boldsymbol{x}_l)\parallel_{H^1(\widehat{\Omega}^j)}\big)\\
&\leq\sum_{k}^{n_1}\big(\parallel \widehat{N}^k(\boldsymbol{x})\parallel_{H^1(\widehat{\Omega}^j)}\big(\sum_{l=1}^{n_2}\lvert \delta N^j_{k, l}\lvert\lvert{u(\boldsymbol{x}_l})\lvert\big)\big)\\
&\leq\ n_2\lvert \delta\boldsymbol{N}^j\lvert_{\infty}\parallel u\parallel_{C^0(\widehat{\Omega}^j)}\sum_k^{n_1}\parallel \widehat{N}^k(\boldsymbol{x})\parallel_{H^1(\widehat{\Omega}^j)},
\end{align*}
where \(\parallel u\parallel_{C^0(\widehat{\Omega}^j)}=\sup_{\boldsymbol{x}\in \widehat{\Omega}^j} |u(\boldsymbol{x})|\) and \(\widehat{N}^k(\boldsymbol{x})\) is the $k-$th basis function of \(\widehat{\boldsymbol{N}}_{\boldsymbol{e}}(\boldsymbol{x})\). By Eq. \eqref{ep:14} and Eq. \eqref{ep:15},
\begin{align*}
\ n_2\lvert\delta{\boldsymbol{N}^j}\lvert_{\infty}\parallel u\parallel_{C^0(\widehat{\Omega}^j)}\sum_{k=1}^{n_1}\parallel\widehat{N}^k(\boldsymbol{x})\parallel_{H^1(\widehat{\Omega}^j)}&\leq n_1n_2\sqrt{C_1^2+C_2^2}\lvert\boldsymbol{\delta N^j}\lvert_{\infty}\parallel u\parallel_{C^0(\widehat{\Omega}^j)}\\
&\leq n_1n_2\sqrt{C_1^2+C_2^2}K\lvert\boldsymbol{\delta N^j}\lvert_{\infty}\parallel u\parallel_{H^2(\widehat{\Omega}^j)},
\end{align*}
The last inequality follows from the Sobolev embedding theorem, where \(K\) is a constant independent of \(u\), and we set \(C_3=n_1n_2\sqrt{C_1^2+C_2^2}K\). Thus, we obtain the boundedness of \(\widehat{\delta\Pi}\):\\
\begin{align*}
\parallel\widehat{\delta\Pi} u\parallel_{H^1(\widehat{\Omega}^j)}\leq\varepsilon C_3\parallel u\parallel_{H^2(\widehat{\Omega}^j)}.
\end{align*}
\end{proof}
\FloatBarrier
\subsection{Interpolation operator on the coarse element \(\Omega^j\)}
For a coarse element \(\Omega^j\), let \(\Pi\) denote the interpolation operator defined by the numerical shape functions \(\{\boldsymbol{N}^j_l(\boldsymbol{x})\}\) as in Eq. \eqref{ep:NBF} , i.e.,
\begin{align}
\Pi{f}=\sum_{l=1}^{n_2}\boldsymbol{N_e}(\boldsymbol{x}){\boldsymbol{N}}_{,l}^jf(\boldsymbol{x}_l)=\boldsymbol{N_e}(\boldsymbol{x})\sum_{l=1}^{n_2}{\boldsymbol{N}}_{,l}^jf(\boldsymbol{x}_l).\label{ep:Toperator}
\end{align}
Correspondingly,  let \(\tilde{\Pi}\) denote the interpolation operator defined by the numerical shape functions \({\boldsymbol{N_e}(\boldsymbol{x})\widetilde{\boldsymbol{N}}^j_{,l}}\) obtained by PIML:
\begin{align}
\widetilde{\Pi}f=\sum_{l=1}^{n_2}\boldsymbol{N_e}(\boldsymbol{x})\widetilde{\boldsymbol{N}}_{,l}^j f(\boldsymbol{x}_l)=\boldsymbol{N_e}(\boldsymbol{x})\sum_{l=1}^{n_2}\widetilde{\boldsymbol{N}}_{,l}^j f(\boldsymbol{x}_l)\label{ep:RPIMLop}.
\end{align}
The difference operator between $\Pi$ and $\widetilde{\Pi}$ is defined by
\begin{align}
\delta\Pi=\Pi-\widetilde{\Pi}=\boldsymbol{N_e}(\boldsymbol{x})\delta{\boldsymbol{N}^j}=\boldsymbol{N_e}(\boldsymbol{x})(\boldsymbol{N}^j-\widetilde{\boldsymbol{N}}^j)=\boldsymbol{N_e}(\boldsymbol{x})\delta\boldsymbol{N}^j.
\end{align}
It is clear that \(\delta{\Pi}\) is still a linear operator. 

From Theorem \ref{thm:T2m2}, \(\boldsymbol{N}_s^j\) is independent of the size of coarse elements, so the generalization error \(\varepsilon\) of \(\boldsymbol{N}_s^j\) is also independent of the coarse-element size. Note that \(\delta\Pi\) is the difference between two interpolation operators constructed from continuous piecewise bilinear functions, so it is affine-invariant, which means that \(\forall\widehat{v}\in H^2(\widehat\Omega)\) and \(\forall v\in H^2(\Omega)\), \(\widehat{\delta\Pi v}=\widehat{\delta\Pi}\widehat{v}\). Let \(\mathbb{P}_k\) be the space of all polynomials of degree at most \(k\). In particular, \(\mathbb{P}_0\) denotes the set of constants.
\begin{lemma}For a coarse element \(\Omega^j\) of size \(h\) and \(u\in H^2(\Omega^j)\), suppose \(\varepsilon=\lvert\delta\boldsymbol{N^j}\lvert_{\infty}\), we have
\label{thm:L3m1}
\begin{align}
\lvert\delta\Pi u\lvert_{H^1(\Omega^j)}&\leq\varepsilon C_5\parallel u\parallel_{H^2(\Omega^j)},\\
\parallel\delta\Pi u\parallel_{L^2(\Omega^j)}&\leq\varepsilon C_6h\parallel u\parallel_{H^2(\Omega^j)}.
\end{align}
\end{lemma}
\begin{proof}
Since both \(\boldsymbol{{ N^j}}\) and \(\widetilde{\boldsymbol{N}}^j\) satisfy the row-wise partition-of-unity condition, for any \(p\in \mathbb{P}_0\), let \(P_l\) denote the vector composed of function values of $p$ at all nodes in \(\Omega^j\). We have
\begin{align*}
\delta\Pi{p}=\boldsymbol{N_e}(\boldsymbol{x})\sum_{l=1}^{n_2}(\boldsymbol{N}_{,l}^jP_l-\widetilde{\boldsymbol{N}}_{,l}^jP_l)=0,
\end{align*}
which implies that
\begin{align}
\delta\Pi u=\delta\Pi u-\delta\Pi p=\delta\Pi(u-p).\nonumber
\end{align}
On the reference element \(\widehat{\Omega}^j\) with unit size, under the affine mapping \(\boldsymbol{B}^j\) and by the  Quotient Mapping Theorem \cite{18}, we have
\begin{align}
\parallel\widehat{\delta\Pi u}\parallel_{H^1(\widehat{\Omega}^j)}&=\parallel\widehat{\delta\Pi(u-p)}\parallel_{H^1(\widehat{\Omega}^j)}\nonumber \\
&\leq\parallel\widehat{\delta\Pi}\parallel \inf_{\hat{p}\in\widehat{\mathbb{P}}_0}\parallel\widehat{u}-\widehat{p}\parallel_{H^2(\widehat{\Omega}^j)}\nonumber \\
&\leq\varepsilon C_3\inf_{\hat{p}\in\widehat{\mathbb{P}}_0}\parallel\widehat{u}-\widehat{p}\parallel_{H^2(\widehat{\Omega}^j)}.\label{ep:23}
\end{align}
Note \(\parallel\widehat{u}-\widehat{p}\parallel^2_{H^2(\widehat{\Omega}^j)}=\parallel\widehat{u}-\widehat{p}\parallel^2_{H^1(\widehat{\Omega}^j)}+\lvert\widehat{u}-\widehat{p}\lvert^2_{H^2(\widehat{\Omega}^j)}\). \\
Because \(\widehat{p}\in\widehat{\mathbb{P}}_0\), we have \(\lvert\widehat{u}-\widehat{p}\lvert^2_{H^2(\widehat{\Omega}^j)}=\lvert\widehat{u}\lvert^2_{H^2(\widehat{\Omega}^j)}\). \\
By the Quotient Mapping Theorem \cite{18}, \(\parallel\widehat{u}-\widehat{p}\parallel^2_{H^1(\widehat{\Omega}^j)}\leq (D_1)^2\lvert\widehat{u}\lvert^2_{H^1(\widehat{\Omega}^j)}\). \\
Combining the above estimates, we have
\begin{align*}
\inf_{\widehat{p}\in\widehat{\mathbb{P}}_0}\parallel\widehat{u}-\widehat{p}\parallel^2_{H^2(\widehat{\Omega}^j)}&\leq (D_1)^2\lvert\widehat{u}\lvert^2_{H^1(\widehat{\Omega}^j)}+\lvert\widehat{u}\lvert^2_{H^2(\widehat{\Omega}^j)}\\
&\leq\max\{(D_1)^2, 1\}(\lvert\widehat{u}\lvert^2_{H^1(\widehat{\Omega}^j)}+\lvert\widehat{u}\lvert^2_{H^2(\widehat{\Omega}^j)}),\nonumber
\end{align*}
combining it with Eq. \eqref{ep:23}:
\begin{align}
\parallel\widehat{\delta\Pi u}\parallel^2_{H^1(\widehat{\Omega}^j)}\leq\varepsilon^2 (C_3)^2\max\{(D_1)^2, 1\}(\lvert\widehat{u}\lvert^2_{H^1(\widehat{\Omega}^j)}+\lvert\widehat{u}\lvert^2_{H^2(\widehat{\Omega}^j)})\label{ep:100}
\end{align}
Using the standard scaling relations for Sobolev norms, we obtain
\begin{align}
\lvert\widehat{u}\lvert_{H^1(\widehat{\Omega}^j)}&\leq {B_1}\parallel\boldsymbol{B}^j\parallel^1det(\boldsymbol{B}^j)^{-\frac{1}{2}}\lvert u\lvert_{H^1({\Omega}^j)}\label{ep:25},\\
\lvert\widehat{u}\lvert_{H^2(\widehat{\Omega}^j)}&\leq {B_2}\parallel\boldsymbol{B}^j\parallel^2det(\boldsymbol{B}^j)^{-\frac{1}{2}}\lvert u\lvert_{H^2({\Omega}^j)}\label{ep:26},\\
\lvert u\lvert_{H^1({\Omega}^j)}&\leq {B_3}\parallel(\boldsymbol{B}^j)^{-1}\parallel^1 det(\boldsymbol{B}^j)^{\frac{1}{2}}\lvert\widehat{u}\lvert_{H^1(\widehat{\Omega}^j)}\label{ep:27},\\
\lvert u\lvert_{H^0({\Omega}^j)}&\leq {B_4}\parallel(\boldsymbol{B}^j)^{-1}\parallel^0 det(\boldsymbol{B}^j)^{\frac{1}{2}}\lvert\widehat{u}\lvert_{H^0(\widehat{\Omega}^j)}\label{ep:28}.
\end{align}
Combining Eqs. \eqref{ep:100}-\eqref{ep:27},
\begin{align}
\lvert\delta\Pi u\lvert^2_{H^1({\Omega}^j)}&\leq ({B_3})^2\parallel(\boldsymbol{B}^j)^{-1}\parallel^2 det(\boldsymbol{B}^j)^{1}\lvert\widehat{\delta\Pi{u}}\lvert^2_{H^1(\widehat{\Omega}^j)}\nonumber\\
&\leq ({B_3})^2\parallel(\boldsymbol{B}^j)^{-1}\parallel^2 det(\boldsymbol{B}^j)^{1}\parallel\widehat{\delta\Pi{u}}\parallel^2_{H^1(\widehat{\Omega}^j)}\nonumber\\
&\leq \varepsilon^2(C_3)^2\max\{(D_1)^2, 1\}({B_3})^2\parallel(\boldsymbol{B}^j)^{-1}\parallel^2 det(\boldsymbol{B}^j)^{1}(\lvert\widehat{u}\lvert^2_{H^1(\widehat{\Omega}^j)}+\lvert\widehat{u}\lvert^2_{H^2(\widehat{\Omega}^j)}).\label{ep:29}
\end{align}
For \(\lvert\widehat{u}\lvert^2_{H^1(\widehat{\Omega}^j)}+\lvert\widehat{u}\lvert^2_{H^2(\widehat{\Omega}^j)}\), we have
\begin{align}
\lvert\widehat{u}\lvert^2_{H^1(\widehat{\Omega}^j)}+\lvert\widehat{u}\lvert^2_{H^2(\widehat{\Omega}^j)}&\leq ({B_1})^2\parallel\boldsymbol{B}^j\parallel^2det(\boldsymbol{B}^j)^{-1}\lvert u\lvert^2_{H^1({\Omega}^j)}+
({B_2})^2\parallel\boldsymbol{B}^j\parallel^4det(\boldsymbol{B}^j)^{-1}\lvert u\lvert^2_{H^2({\Omega}^j)}.\label{ep:30}
\end{align}
For a rectangular mesh, there exists a constant \(C_4\) such that
\begin{align}
\parallel\boldsymbol{B}^j\parallel\parallel(\boldsymbol{B}^j)^{-1}\parallel\leq C_4,\label{ep:31}
\end{align}
and there exists a constant \(D_2\) satisfying
\begin{align}
\parallel\boldsymbol{B}^j\parallel\leq D_2h.\label{ep:32}
\end{align}
Substituting Eqs. \eqref{ep:30}-\eqref{ep:32} into Eq. \eqref{ep:29}, we get
\begin{align}
\lvert\delta\Pi u\lvert^2_{H^1({\Omega}^j)}\leq\varepsilon^2(C_3)^2\max\{(D_1)^2, 1\}({B_3})^2[({B_1})^2(C_4)^2\lvert u\lvert^2_{H^1({\Omega}^j)}+({B_2})^2(C_4)^2(D_2)^2h^2\lvert u\lvert^2_{H^2({\Omega}^j)}],\nonumber
\end{align}
Because \(\lvert u\lvert_{H^1({\Omega}^j)}\leq\parallel u\parallel_{H^2({\Omega}^j)}\) and \(\lvert u\lvert_{H^2({\Omega}^j)}\leq\parallel u\parallel_{H^2({\Omega}^j)}\), we have:
\begin{align*}
\lvert\delta\Pi u\lvert^2_{H^1({\Omega}^j)}\leq\varepsilon^2(C_3)^2\max\{(D_1)^2, 1\}({B_3})^2[({B_1})^2(C_4)^2\parallel u\parallel^2_{H^2({\Omega}^j)}+({B_2})^2(C_4)^2(D_2)^2h^2\parallel u\parallel^2_{H^2({\Omega}^j)}].\nonumber
\end{align*}
For sufficiently small $h$, the second term is smaller than the first, leading to
\begin{align}
\lvert\delta\Pi u\lvert^2_{H^1({\Omega}^j)}\leq2\varepsilon^2(C_3)^2\max\{(D_1)^2, 1\}({B_3})^2({B_1})^2(C_4)^2\parallel u\parallel^2_{H^2({\Omega}^j)}.\nonumber
\end{align}
Let \(C_5=C_3\sqrt{2\max\{(D_1)^2, 1\}}{B_3}{B_1}C_4\). Hence, we obtain:
\begin{align}
\lvert\delta\Pi u\lvert_{H^1({\Omega}^j)}\leq\varepsilon C_5\parallel u\parallel_{H^2({\Omega}^j)}.\label{ep:35}
\end{align}
For \(\parallel\delta\Pi u\parallel_{L^2}\), we replace the left hand side of Eq. \eqref{ep:100} by \(\lvert \widehat{\delta\Pi{u}}\lvert_{H^0}\) and then combine Eq. \eqref{ep:28}, which gives:
\begin{align}
\lvert\delta\Pi u\lvert_{H^0({\Omega}^j)}\leq\varepsilon C_6h\parallel u\parallel_{H^2({\Omega}^j)},\label{ep:36}
\end{align}
where \(C_6=2C_3D_2B_2B_4C_4\), noting that \(\parallel\delta\Pi u\parallel_{L^2({\Omega}^j)}=\lvert\delta\Pi u\lvert_{H^0({\Omega}^j)}\). This completes the proof. 
\end{proof}
\FloatBarrier
\subsection{Error estimates for the PIML interpolation operator}
Note that Lemma \ref{thm:L3m1} alone is insufficient to establish the error estimate for PIML approximation, because the interpolation operator \(\Pi\) in Eq. \eqref{ep:Toperator}  defined by the theoretical substructure method lacks linear completeness for the heterogeneous problem. We must first prove the convergence of its interpolation error. 
\begin{lemma}Suppose that \(\Omega\) is a two- or three-dimensional domain equipped with a rectangular or cuboid mesh \(\mathcal{T}\) and that \(u\in H^2(\Omega)\) is the exact solution of a second-order linear elastic problem, where the governing equations are required to satisfy the properties given in Eqs. \eqref{ep:symmetric} - \eqref{ep:coercive}, and \(f\in L^2(\Omega)\) is the body force, we have\label{thm:L3m4}
\begin{equation}
\parallel u-\Pi u\parallel_{H^1(\Omega)}\leq C_{0}h\parallel u\parallel_{H^2(\Omega)}+C_{7}h\parallel f\parallel_{L^2(\Omega)}.
\end{equation}
\end{lemma}
\begin{proof}
On each coarse element \(\Omega^j\), let \(\boldsymbol{u}^j_b\) denote the vector of nodal values of $u$ on the boundary, then consider the FEM solution \({u}_h^j\) on \(\Omega^j\), we have 
\begin{align}
\parallel u-{u}_h^j\parallel_{H^1(\Omega^j)}\leq C_0h\parallel u\parallel_{H^2(\Omega^j)},\nonumber
\end{align}
where $C_0$ is a constant. The local finite element solution \({u}_h^j\) can be obtained by:
\begin{align}
\boldsymbol{K}^j\boldsymbol{u}^j=\left(\begin{array}{cc}
		\boldsymbol{K}_{bb}^j&(\boldsymbol{K}^j_{ib})^T\\
		\boldsymbol{K}^j_{ib}&\boldsymbol{K}^j_{ii}
	\end{array}\right)\left(\begin{matrix}
\boldsymbol{u}^j_{b} \\
\boldsymbol{u}^j_{i}
\end{matrix}\right)=\left(\begin{matrix}
\boldsymbol{F}^j_{b} \\
\boldsymbol{F}^j_{i}
\end{matrix}\right),\nonumber
\end{align}
where \(\boldsymbol{u}^j\) is the vector of nodal values of \({u}^j_h\). Then \({u}^j_h\) can be represented as:
\begin{align*}
{u}^j_h=\boldsymbol{N_e}(\boldsymbol{x})\boldsymbol{u}^j=\boldsymbol{N_e}(\boldsymbol{x})[\boldsymbol{N^j}\boldsymbol{u}^j_b+\left(\begin{matrix}\boldsymbol{0}\\
(\boldsymbol{K}^{j}_{i, i})^{-1}\boldsymbol{F}_i^j\end{matrix}\right)],
\end{align*} 
where \(\boldsymbol{0}\) denotes the zero vector with the same dimension as \(\boldsymbol{u}^j_b\). Now note that \(u-{u}_h^j\) is close to the restriction of \(u-\Pi u\) on \(\Omega^j\), it has the following form:
\begin{align}
u-{u}_h^j=u-\boldsymbol{N_e}(\boldsymbol{x})[\boldsymbol{N^j}\boldsymbol{u}^j_b+\left(\begin{matrix}\boldsymbol{0}\\
(\boldsymbol{K}^{j}_{i, i})^{-1}\boldsymbol{F}_i^j\end{matrix}\right)]=u-\Pi u-\boldsymbol{N_e}^i(\boldsymbol{x})(\boldsymbol{K}^{j}_{i, i})^{-1}\boldsymbol{F}_i^j,\nonumber
\end{align}
where \(\boldsymbol{N_e}^i(\boldsymbol{x})\) denotes the row vector of shape functions of the fine element associated with the interior nodes of the coarse element \(\Omega^j\). We next establish an estimate for the load-dependent correction term \((\boldsymbol{K}^{j}_{i, i})^{-1}\boldsymbol{F}_i^j\)  as \(h\rightarrow0\):
\begin{align}
\lvert F_k^j\mid&=\lvert\int_{\Omega^j_k}N_e^k(\boldsymbol{x})f(\boldsymbol{x})d\Omega^j\mid\nonumber\\
&\leq\int_{\Omega^j}\lvert N_e^k(\boldsymbol{x})f(\boldsymbol{x})\lvert d\Omega^j\nonumber\\
&\leq\parallel N_e^k(\boldsymbol{x})\parallel_{L^2(\Omega^j)}\parallel f(\boldsymbol{x})\parallel_{L^2(\Omega^j)}.
\end{align}
Because the diameter of \(\Omega^j\) is $h$, 
\begin{align}
\parallel N_e^k(\boldsymbol{x})\parallel_{L^2(\Omega^j)}\leq \widehat{C}_1h\quad\Omega\subset\mathbb{R}^2,\nonumber
\end{align}
or
\begin{align}
\parallel N_e^k(\boldsymbol{x})\parallel_{L^2(\Omega^j)}\leq\widehat{C}_1h^{\frac{3}{2}}\quad\Omega\subset\mathbb{R}^3.\nonumber
\end{align}
Furthermore, every entry of \((\boldsymbol{K}^{j}_{i, i})^{-1}\) satisfies the following bound:
\begin{align}
\lvert(\boldsymbol{K}^{j}_{i, i})^{-1}(s, t)\lvert\leq \widehat{C}_2\quad\Omega\subset\mathbb{R}^2,\nonumber
\end{align}
or
\begin{align}
\lvert(\boldsymbol{K}^{j}_{i, i})^{-1}(s, t)\lvert\leq \widehat{C}_2h^{-1}\quad\Omega\subset\mathbb{R}^3,\nonumber
\end{align}
where constants $\widehat{C}_1$ and $\widehat{C}_2$ depend only on the shape of \(\Omega^j\). Thus, for the vector \((\boldsymbol{K}^{j}_{i, i})^{-1}\boldsymbol{F}_i^j\), each entry satisfies the following:
\begin{align}
\lvert{((\boldsymbol{K}^{j}_{i, i})^{-1}\boldsymbol{F}_i^j)_k}\lvert&\leq m^2\widehat{C}_1\widehat{C}_2h\parallel f(\boldsymbol{x})\parallel_{L^2(\Omega^j)}\quad\Omega\subset\mathbb{R}^2,\nonumber\\
\lvert{((\boldsymbol{K}^{j}_{i, i})^{-1}\boldsymbol{F}_i^j)_k}\lvert&\leq m^3\widehat{C}_1\widehat{C}_2h^{\frac{1}{2}}\parallel f(\boldsymbol{x})\parallel_{L^2(\Omega^j)}\quad\Omega\subset\mathbb{R}^3.\nonumber
\end{align}
Consequently, we obtain:
\begin{align}
\parallel\boldsymbol{N_e}^i(\boldsymbol{x})(\boldsymbol{K}^{j}_{i, i})^{-1}\boldsymbol{F}_i^j\parallel_{H^1(\Omega^j)}&\leq\sum_k\parallel N_e^k(\boldsymbol{x}){((\boldsymbol{K}^{j}_{i, i})^{-1}\boldsymbol{F}_i^j)_k}\parallel_{H^1(\Omega^j)}\nonumber\\
&\leq m^4\widehat{C}_1\widehat{C}_2\widehat{C}_3h\parallel f(\boldsymbol{x})\parallel_{L^2(\Omega^j)}\quad\Omega\subset\mathbb{R}^2.\nonumber
\end{align}
For \(\Omega\subset\mathbb{R}^3\), \(m^4\) is replaced by \(m^6\), the last inequality follows from:
\begin{align*}
\parallel N_e^k(\boldsymbol{x})\parallel_{H^1(\Omega^j)}&\leq \widehat{C}_3\quad\Omega\subset\mathbb{R}^2,\\
\parallel N_e^k(\boldsymbol{x})\parallel_{H^1(\Omega^j)}&\leq \widehat{C}_3h^{\frac{1}{2}}\quad\Omega\subset\mathbb{R}^3.
\end{align*}
Setting \(C_7=m^4\widehat{C}_1\widehat{C}_2\widehat{C}_3\) for \(\Omega\subset\mathbb{R}^2\) or \(C_7=m^6\widehat{C}_1\widehat{C}_2\widehat{C}_3\) for \(\Omega\subset\mathbb{R}^3\),  we obtain the estimate :
\begin{align}
\parallel u-\Pi u\parallel_{H^1(\Omega^j)}&=\parallel u-{u}_h^j+\boldsymbol{N_e}^i(\boldsymbol{K}^{j}_{i, i})^{-1}\boldsymbol{F}_i^j\parallel_{H^1(\Omega^j)}\nonumber\\
&\leq\parallel u-{u}_h^j\parallel_{H^1(\Omega^j)}+\parallel\boldsymbol{N_e}^i(\boldsymbol{x})(\boldsymbol{K}^{j}_{i, i})^{-1}\boldsymbol{F}_i^j\parallel_{H^1(\Omega^j)}\nonumber\\
&\leq C_0h\parallel u\parallel_{H^2(\Omega^j)}+C_7h\parallel f\parallel_{L^2(\Omega^j)}.\nonumber
\end{align}
Summing over all coarse elements yields:
\begin{align}
\parallel u-\Pi u\parallel^2_{H^1(\Omega)}&=\sum_j\parallel u-\Pi u\parallel^2_{H^1(\Omega^j)}\nonumber\\
&\leq2\sum_j(\parallel u-{u}_h^j\parallel^2_{H^1(\Omega^j)}+\parallel\boldsymbol{N_e}^i(\boldsymbol{x})(\boldsymbol{K}^{j}_{i, i})^{-1}\boldsymbol{F}_i^j\parallel^2_{H^1(\Omega^j)})\nonumber\\
&\leq2\sum_j(C^2_0h^2\parallel u\parallel^2_{H^2(\Omega^j)}+C^2_7h^2\parallel f\parallel^2_{L^2(\Omega^j)}).\nonumber
\end{align}
For rectangular meshes, $C_0$ and $C_7$ have uniform upper bounds in the whole domain $\Omega$. Thus, absorbing the uniform bounds into $C_0$ and $C_7$ and also using the same symbols, we have
\begin{equation*}
\parallel u-\Pi u\parallel_{H^1(\Omega)}\leq C_{0}h\parallel u\parallel_{H^2(\Omega)}+C_{7}h\parallel f\parallel_{L^2(\Omega)}.
\end{equation*}
\end{proof}

By Lemma \ref{thm:L3m4}, we can obtain the error estimation of the solution by the PIML as follows.

\begin{theorem}Suppose that \(\Omega\) is a two- or three-dimensional domain equipped with a rectangular or cuboid mesh \(\mathcal{T}\) and that \(u\in H^2(\Omega)\) is the exact solution of a second-order linear elastic problem, where the governing equations are required to satisfy the properties given in Eqs. \eqref{ep:symmetric} - \eqref{ep:coercive}. Let \(\widetilde{u}_h\in H^1(\Omega)\) be its approximation by the PIML, and \(\varepsilon=\max_{\Omega_j\subset \mathcal{T}}\lvert\delta\boldsymbol{N^j}\lvert_{\infty}\) be the generalization error for matrix learning. Then the following estimates hold:
\label{thm:T3m4}
\begin{align}
&\parallel u-\widetilde{u}_h\parallel_{H^1(\Omega)}\leq (M_1h+M_2\varepsilon)\parallel u\parallel_{H^2(\Omega)}+M_3h\parallel f\parallel_{L^2(\Omega)},\\
&\parallel u-\widetilde{u}_h\parallel_{L^2(\Omega)}\leq (M_4h^2+M_5h\varepsilon+M_6\varepsilon^2)\parallel u\parallel_{H^2(\Omega)}+(M_7h^2+M_8h\varepsilon)\parallel f\parallel_{L^2(\Omega)}.
\end{align}
\end{theorem}
\begin{proof}
Denote \(\Omega^j\) by the $j-$th element of \(\mathcal{T}\),  \(\Pi\) by the corresponding theoretical interpolation operator in Eq.\eqref{ep:Toperator} and \(\widetilde{\Pi}\) by the PIML interpolation operator in Eq.\eqref{ep:RPIMLop}. We have\\
\begin{align}
\parallel u-\widetilde{\Pi}u\parallel_{H^1(\Omega)}\leq\parallel u-\Pi u\parallel_{H^1(\Omega)}+\parallel\delta\Pi u\parallel_{H^1(\Omega)},\nonumber
\end{align}
where \(\delta\Pi=\Pi-\widetilde\Pi\), and note that the first term satisfies
\begin{align}
\parallel u-\Pi u\parallel_{H^1(\Omega)}\leq C_{0}h\parallel u\parallel_{H^2(\Omega)}+C_7h\parallel f\parallel_{L^2(\Omega)},\nonumber
\end{align}
where $h$ is the maximum element diameter in \(\mathcal{T}\). For the second term, on each element \(\Omega^j\), Eq. \eqref{ep:35} and Eq. \eqref{ep:36} give
\begin{align}
\parallel\delta\Pi u\parallel_{H^1(\Omega^j)}\leq\sqrt{\varepsilon^2 (C_{5}^j)^2\parallel u\parallel^2_{H^2(\Omega^j)}+\varepsilon^2 (C^j_6)^2h^2\parallel u\parallel^2_{H^2(\Omega^j)}}.\nonumber
\end{align}
Thus,
\begin{align}
\parallel\delta\Pi u\parallel^2_{H^1(\Omega)}\leq\sum_{j=1}\varepsilon^2 (C_{5}^j)^2\parallel u\parallel^2_{H^2(\Omega^j)}+\varepsilon^2 (C^j_6)^2h^2\parallel u\parallel^2_{H^2(\Omega^j)}.\label{ep:56}
\end{align}
For sufficiently small $h$, the second term is smaller than the first, and since \(\mathcal{T}\) is rectangular, all constants in Eq. \eqref{ep:56} are uniformly bounded, thus:
\begin{align}
\parallel\delta\Pi u\parallel_{H^1(\Omega)}\leq\varepsilon C_{5}\parallel u\parallel_{H^2(\Omega)}.\nonumber
\end{align}
We obtain the \(H^1-\)norm estimate:
\begin{align}
\parallel u-\widetilde{u}_h\parallel_{H^1(\Omega)}&\leq C_c\parallel u-\widetilde{\Pi}u\parallel_{H^1(\Omega)}\nonumber\\
&\leq C_c(C_{0}h+\varepsilon C_{5})\parallel u\parallel_{H^2(\Omega)}+C_cC_7h\parallel f\parallel_{L^2(\Omega)},\nonumber
\end{align}
where \(C_c=\frac{M}{\alpha}\) is the control constant of C\'ea's lemma. For the \(L^2-\)norm estimate, using the Aubin-Nitsche argument, we have
\begin{align*}
\parallel u-\widetilde{u}_h\parallel_{L^2(\Omega)}\leq C_{8}\parallel u-\widetilde{u}_h\parallel_{H^1(\Omega)}\sup_{g\in L^2(\Omega)}\{\frac{1}{\parallel g\parallel_{L^2(\Omega)}}\inf_{v_h\in\widetilde{U}_h}\parallel \Phi-v_h\parallel_{H^1(\Omega)}\},
\end{align*}
where \(\widetilde{U}_h\) is the space spanned by the numerical shape functions obtained by PIML, and \(\Phi\) is the solution of the dual problem of Eq. \eqref{ep:weakform}:
\begin{align*}
B(v,\Phi)=<g,v>, \forall v\in H^1(\Omega), g\in L^2(\Omega)
\end{align*}
Due to the regularity estimate of the second-order elliptic problem \cite{18}, we have \(\Phi\in H^2(\Omega)\) and
\begin{align*}
\parallel\Phi\parallel_{H^2(\Omega)}\leq C_9\parallel{g}\parallel_{L^2(\Omega)}.
\end{align*}
Combining our previous analysis, we obtain
\begin{align}
&\parallel u-\widetilde{u}_h\parallel_{L^2(\Omega)}\leq C_8C_c(C_0C_9h+\varepsilon C_5C_9+C_7h)((C_0h+\varepsilon C_5)\parallel u\parallel_{H^2(\Omega)}+C_7h\parallel f\parallel_{L^2(\Omega)}).\nonumber
\end{align}
Let \(M_1=C_cC_0, M_2=C_cC_5, M_3=C_cC_7\) and \(M_4=C_8C_9C_cC_0^2+C_0C_7C_8C_c, M_5=2C_0C_5C_8C_9C_c+C_5C_7C_8C_c, M_6=C_8C_9C_cC_5^2, M_7=C_7C_8C_c(C_0C_9+C_7), M_8=C_5C_7C_8C_9C_c\). This completes the proof.
\end{proof}

Theorem \ref{thm:T3m4} shows that the error of PIML is controlled by both the mesh size \(h\) and the generalization error \(\varepsilon\) of the matrix learned by PIML. Consequently, for a fixed trained network, mesh refinement alone may lead to a nonzero error upper bound. Although the error floor can be reduced by decreasing the network generalization error \(\varepsilon\), it is desirable to construct a method whose convergence is governed primarily by the mesh size, because reducing the generalization error may be considerably more difficult than refining the mesh. This problem can be addressed by improving the PIML by imposing linear completeness constraints on the trained numerical shape functions, which is presented in the next section.

\section{Improved PIML with linear completeness}\label{Sect. 4}
Let \(\widetilde{U_h}\) denote the space spanned by the learned numerical shape functions \({\boldsymbol{N_e}(\boldsymbol{x})\widetilde{\boldsymbol{N}}^j}\) obtained by PIML, which is a subspace of \(V_h\) spanned by the piecewise bilinear shape functions \(\boldsymbol{N_e}(\boldsymbol{x})\) defined on the fine elements. Note that \(\mathbb{P}_1\subset V_h\). 
However, the solution in \(\widetilde{U_h}\) is not the best approximation in \(V_h\), but rather in a space of piecewise polynomials that may not contain \(\mathbb{P}_1\). Without linear completeness, the standard optimal-order approximation estimates cannot be obtained directly from polynomial reproduction arguments. Therefore, if \(\widetilde{U_h}\) contains piecewise linear polynomial fields and satisfies suitable stability conditions, optimal order approximation estimates can be established.
\subsection{The linear completeness constraints on the numerical shape functions}
To achieve this goal, we introduce an improved PIML formulation with linear completeness. In the original formulation, \(\widetilde{U_h}\) was required to satisfy the partition-of-unity condition. We now require it to satisfy the linear-completeness condition introduced in Definition \ref{thm:D2f3} of Section \ref{sect:2S2}.
To impose linear completeness, consider a coarse element: when the displacements of the boundary nodes follow a linear displacement field, the internal nodes should also conform to the same displacement field. 

For a coarse element divided into \(m\times m\) fine elements without boundary constraints, there are \(n_b=4m\) boundary nodes and \(n_i=(m-1)^2\) interior nodes. Let \([x_{b_1}, x_{b_2},\dots,x_{b_{n_b}}]\) and \([y_{b_1},y_{b_2},\dots,y_{b_{n_b}}]\) be the \(x-\) and \(y-\)coordinates of the boundary nodes respectively. Let \([x_{i_1}, x_{i_2},\dots,x_{i_{n_i}}]\) and \([y_{i_1},y_{i_2},\dots,y_{i_{n_i}}]\) be the \(x-\) and \(y-\)coordinates of the interior nodes respectively. Denote by \(\boldsymbol{T}^j_b\) the \(2n_b\times6\) matrix as follows:
\begin{align}
\boldsymbol{T}^j_b=\left(
	\begin{array}{cccccc}
		1&0&x_{b_1}&0&y_{b_1}&0\\
		0&1&0&x_{b_1}&0&y_{b_1}\\
		1&0&x_{b_2}&0&y_{b_2}&0\\
		0&1&0&x_{b_2}&0&y_{b_2}\\
		&&&\ldots\ldots&&\\
		1&0&x_{b_{n_b}}&0&y_{b_{n_b}}&0\\
		0&1&0&x_{b_{n_b}}&0&y_{b_{n_b}}\\
	\end{array}
\right).\nonumber
\end{align}
Similarly, denote by \(\boldsymbol{T}^j_i\) the \(2n_i\times6\) matrix as:
\begin{align}
\boldsymbol{T}^j_i=\left(
	\begin{array}{cccccc}
		1&0&x_{i_1}&0&y_{i_1}&0\\
		0&1&0&x_{i_1}&0&y_{i_1}\\
		1&0&x_{i_2}&0&y_{i_2}&0\\
		0&1&0&x_{i_2}&0&y_{i_2}\\
		&&&\ldots\ldots&&\\
		1&0&x_{i_{n_i}}&0&y_{i_{n_i}}&0\\
		0&1&0&x_{i_{n_i}}&0&y_{i_{n_i}}\\
	\end{array}
\right).\nonumber
\end{align}
According to the linear completeness conditions in Eq.\eqref{ep:Linearcom}, the improved PIML transition matrix $\widetilde{\boldsymbol{N}}^j$ must then satisfy
\begin{align}
\boldsymbol{T}^j_i=\widetilde{\boldsymbol{N}}^j\boldsymbol{T}^j_b=\widetilde{\boldsymbol{N}}^j_{1}\boldsymbol{T}^j_{b_1}+\widetilde{\boldsymbol{N}}^j_2\boldsymbol{T}^j_{b_2}.\label{ep:65}
\end{align}
We order the boundary nodes so that the first three selected nodes are non-collinear. Denote by \(\boldsymbol{T}^j_{b1}\) the \(6\times6\) submatrix composed of the first 6 rows of \(\boldsymbol{T}^j_b\), and by \(\boldsymbol{T}^j_{b2}\) the \((2n_b-6)\times6\) submatrix composed of other rows of \(\boldsymbol{T}^j_b\). The matrix \(\widetilde{\boldsymbol{N}}^j\) is divided into two submatrices \(\widetilde{\boldsymbol{N}}^j_1\) and \(\widetilde{\boldsymbol{N}}^j_2\) according to \(\boldsymbol{T}^j_{b1}\) and \(\boldsymbol{T}^j_{b2}\). That is, \(\widetilde{\boldsymbol{N}}^j_1\) is composed of the first 6 columns of \(\widetilde{\boldsymbol{N}}^j\), and \(\widetilde{\boldsymbol{N}}^j_2\) is composed of the other columns of \(\widetilde{\boldsymbol{N}}^j\). We note that the submatrix \(\boldsymbol{T}^j_{b_1}\) associated with three non-collinear boundary nodes is invertible. Thus, we have
\begin{align}
\widetilde{\boldsymbol{N}}^j_{1}=(\boldsymbol{T}^j_i-\widetilde{\boldsymbol{N}}^j_{2}\boldsymbol{T}^j_{b_2}){(\boldsymbol{T}^j_{b_1})}^{-1}\label{ep:processL}.
\end{align}
This means that we only need to learn the \(8m-6\) columns of \(\widetilde{\boldsymbol{N}}_2^j\) for the two-dimensional case. Similarly, for the three-dimensional case, we only need to learn \(18m^2-6\) columns of \(\widetilde{\boldsymbol{N}}_2^j\). Then, for the improved PIML, denote by 
\(\widetilde{\widetilde{\boldsymbol{N}}}^j\) the new transition matrix obtained by Eq. \eqref{ep:processL} and the corresponding learned numerical shape functions are defined by 
\begin{align}
\widetilde{\widetilde{\boldsymbol{N}}}_l^{j}(\boldsymbol{x})=\boldsymbol{N_e}(\boldsymbol{x})\widetilde{\widetilde{\boldsymbol{N}}}_{, l}^j, l=1,2,\dots,2n_b,\label{ep:NBF2}
\end{align}
where \(\widetilde{\widetilde{\boldsymbol{N}}}_{, l}^j\) denotes the $l-$th column of $\widetilde{\widetilde{\boldsymbol{N}}}^j$. 
Denote by \(\widetilde{\widetilde{\Pi}}\) the interpolation operator defined by \(\widetilde{\widetilde{\boldsymbol{N}}}^j\):
\begin{align}
\widetilde{\widetilde{\Pi}}=\sum_{l=1}^{n_2}\boldsymbol{N_e}(\boldsymbol{x})\widetilde{\widetilde{\boldsymbol{N}}}_{,l}^jf(\boldsymbol{x}_l)\label{ep:LPIMLop}.
\end{align}

Furthermore, if the numerical shape functions are subject to boundary constraints, we modify \(\boldsymbol{T}^j_b\) so that its entries correspond to the linear displacement field at these constrained nodes. For example, the \(\boldsymbol{T}^j_b\) with linear boundary constraints only has 8 rows corresponding to the four nodes of a coarse element. 

\begin{figure}[pos=h]
    \centering
    \begin{subfigure}{0.45\textwidth}
        \centering
        \includegraphics[trim=2cm 9cm 2cm 9cm, clip, width=0.8\textwidth]{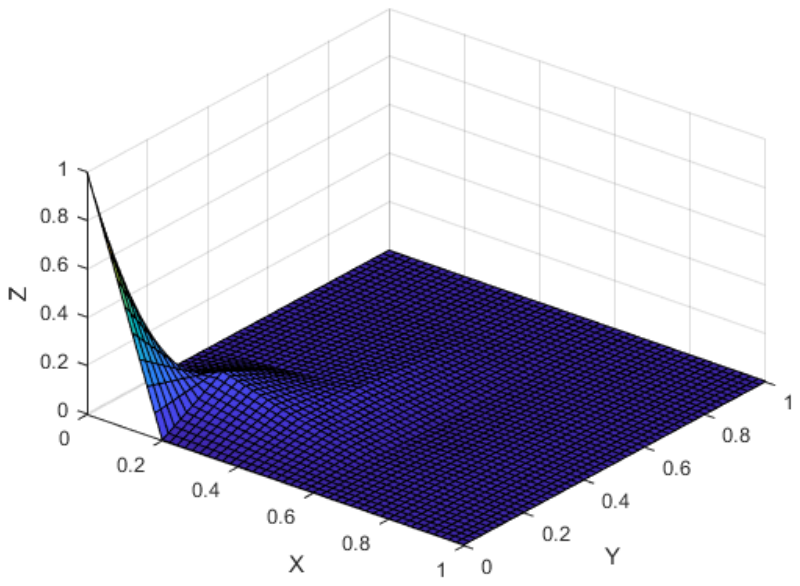}
        \caption{The x-component of the first shape function}
        \label{fig:sub:a}
    \end{subfigure}
    \hfill
    \begin{subfigure}{0.45\textwidth}
        \centering
        \includegraphics[trim=2cm 9cm 2cm 9cm, clip, width=0.8\textwidth]{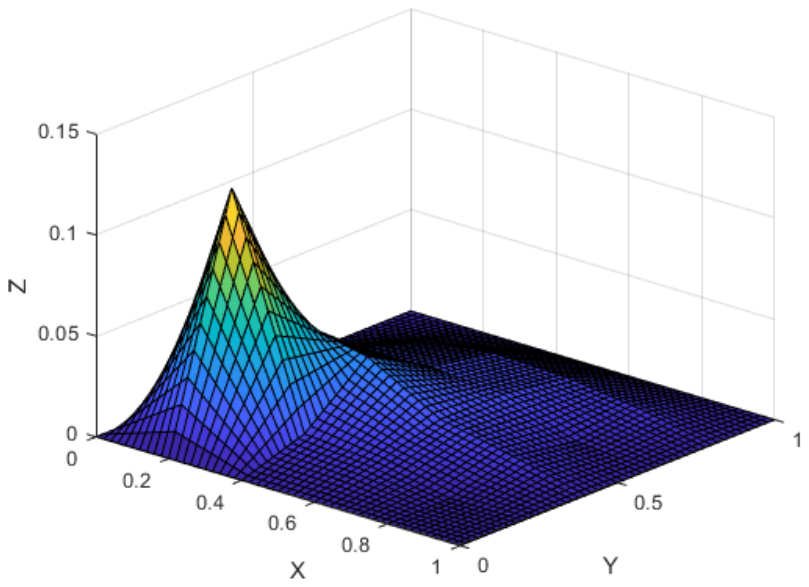}
        \caption{The y-component of the first shape function}
        \label{fig:sub:b}
    \end{subfigure}
    \caption{Numerical shape functions with linear completeness and no boundary constraint}
    \label{fig:shape1}
\end{figure}

\begin{figure}[pos=h]
    \centering
    \begin{subfigure}{0.45\textwidth}
        \centering
        \includegraphics[trim=2cm 9cm 2cm 9cm, clip, width=0.8\textwidth]{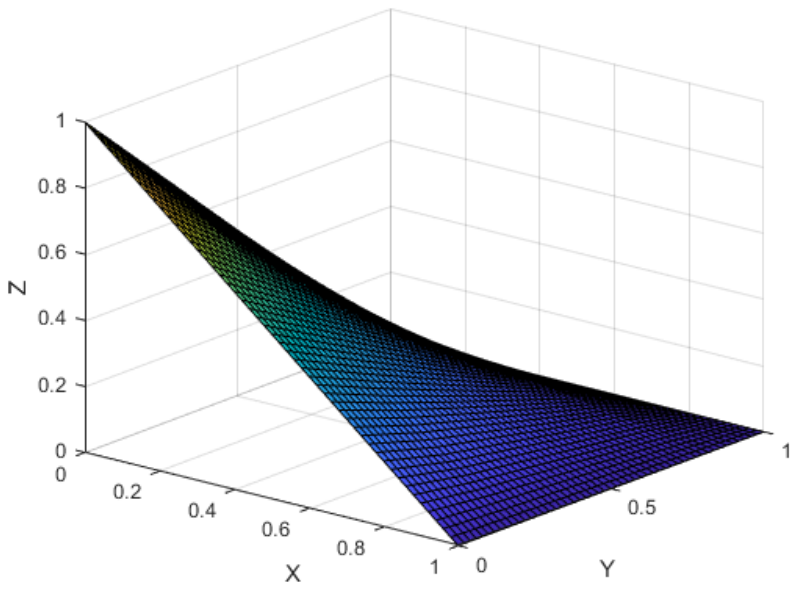}
        \caption{The x-component of the first shape function}
        \label{fig:sub:a}
    \end{subfigure}
    \hfill
    \begin{subfigure}{0.45\textwidth}
        \centering
        \includegraphics[trim=2cm 9cm 2cm 9cm, clip, width=0.8\textwidth]{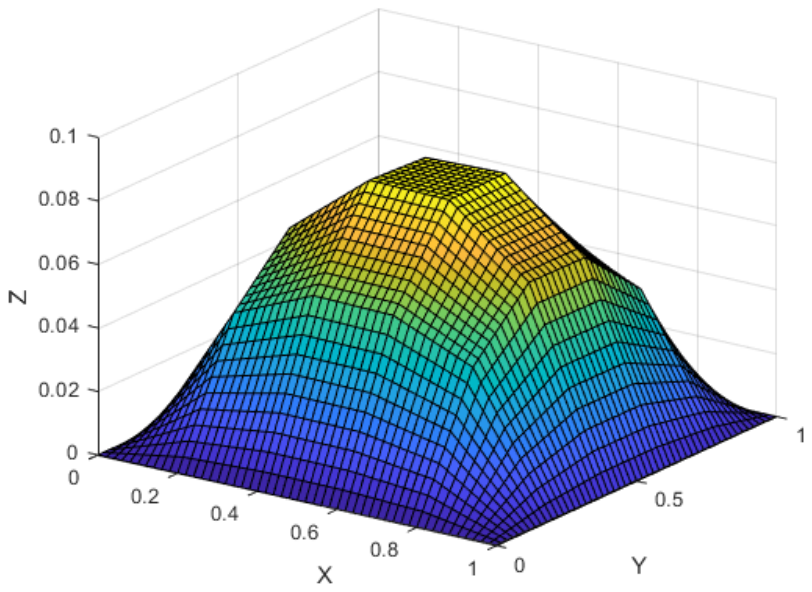}
        \caption{The y-component of the first shape function}
        \label{fig:sub:b}
    \end{subfigure}
    \caption{Numerical shape function with linear completeness and a linear boundary constraint}
    \label{fig:shape2}
\end{figure}

Based on Theorem \ref{thm:T2m1}, it is clear that the learned numerical shape functions defined by Eq. \eqref{ep:NBF2} have linear completeness. Examples of the resulting numerical shape functions are shown in \autoref{fig:shape1} and \autoref{fig:shape2}.

To estimate the interpolation error by the \(\widetilde{\widetilde{\Pi}}\),  
let \(\bar{\Pi}\) be the interpolation operator obtained by the transition matrix $\bar{\boldsymbol{N}}^j$ satisfying the linear completeness conditions through the reconstruction formula \eqref{ep:processL} as follows:
\begin{align}\label{ep:barNj}
\bar{\boldsymbol{N}}^j_1&=(\boldsymbol{T}^j_i-\boldsymbol{N}^j_2\boldsymbol{T}^j_{b_2})(\boldsymbol{T}^j_{b_1})^{-1},\\
\bar{\boldsymbol{N}}^j_2&=\boldsymbol{N}^j_2,
\end{align}
where $\boldsymbol{N}^j_2$ is the submatrix of $\boldsymbol{N}^j$ in Eq. \eqref{ep:Nj} by the substructure method. Obviously, by Lemmas \ref{thm:bound1} and  \ref{thm:bound2}, the operator \(\bar{\Pi}\) is uniformly bounded. Under the uniform boundedness assumption on the exact transition matrices and the geometric regularity condition imposed on the reconstruction nodes, the reconstructed transition matrices are uniformly bounded. Consequently, the associated interpolation operators are uniformly bounded from \(H^2(\widehat\Omega)\) to \(H^1(\widehat\Omega)\).
Note that if the problem involves a homogeneous material, then \(\bar{\Pi}\equiv\Pi\) in Eq. \eqref{ep:Toperator}. 

Suppose the generalization error of learning for every entry in \(\widetilde{\boldsymbol{N}}^j_{2}\) is \(\delta_2\), then the error \(\delta_1\) for the entries in \(\widetilde{\boldsymbol{N}}^j_{1}\) by Eq. \eqref{ep:processL} satisfies
\begin{align*}
\delta_1\leq C_T\delta_2,
\end{align*}
where \(C_T=\parallel{\boldsymbol{T}^j_{b_2}{(\boldsymbol{T}^j_{b_1})}^{-1}}\parallel_1\). 
For any coarse element \(\Omega^j\) with size \(h\), it is clear that \(\boldsymbol{T}^j_{b_2}{(\boldsymbol{T}^j_{b_1})}^{-1}\) is independent of \(h\). For rectangular reference element \(\widehat{\Omega}=[0,1]\times[0,1]\), we choose the three nodes \((0,0), (1,0), (1,1)\) for \(\boldsymbol{T}_{b_1}^j\), then  we can obtain that the absolute value of the entries of \(\{(T^j_{b_1})^{-1}\}\) is bounded by 1, thus, \(\parallel{{(\boldsymbol{T}^j_{b_1})}^{-1}}\parallel_1\leq6\). Since other boundary nodes for \(\boldsymbol{T}^j_{b_2}\) are located in \([0,1]\times[0,1]\), then the \(l_1\) norm \(\parallel\boldsymbol{T}^j_{b_2}\parallel_1\leq 2n_b-6\).
So \(C_{T}\leq 6(2n_b-6)\), it implies \(\delta_1=O(\delta_2)\). Hence, we can also suppose the generalization error for matrix learning as \(\varepsilon=\lvert\delta\bar{\boldsymbol{N}}^j\lvert_{\infty}=\lvert\bar{\boldsymbol{N}}^j-\widetilde{\widetilde{\boldsymbol{N}}}^j\lvert_{\infty}\).

\subsection{The error estimates of the improved PIML}
We now establish error estimates and convergence results for the improved numerical shape functions in Eq. \eqref{ep:NBF2}. By the definitions, both interpolation operators \(\bar{\Pi}\) and \(\widetilde{\widetilde{\Pi}}\) have linear completeness. Denote $\delta\bar{\Pi}=\bar{\Pi}-\widetilde{\widetilde{\Pi}}$, then we have,
\begin{lemma} 
For a coarse element \(\Omega^j\) of size \(h\) and \(u\in H^2(\Omega^j)\), suppose \(\varepsilon=\lvert\delta\bar{\boldsymbol{N}}^j\lvert_{\infty}\), we have
\label{thm:L4m1}
\begin{align}
\lvert\delta\bar{\Pi} u\lvert_{H^1({\Omega}^j)}&\leq\varepsilon C_{11}h\lvert u\lvert_{H^2({\Omega}^j)},\label{ep:DL1}\\
\parallel\delta\bar{\Pi} u\parallel_{L^2({\Omega}^j)}&\leq\varepsilon C_{12}h^2\lvert u\lvert_{H^2({\Omega}^j)}.
\end{align}
\end{lemma}
\begin{proof}
Since the two transition matrices satisfy the linear completeness condition, for all \(p\in \mathbb{P}_1\), let \(P_l\) denote the vector composed of function values of \(p\) at all nodes in \(\Omega^j\). We have
\begin{align}
\delta\bar{\Pi}{p}=\boldsymbol{N_e}(\boldsymbol{x})\sum_{l=1}^{n_2}(\bar{\boldsymbol{N}}_{,l}^jP_l-\widetilde{\boldsymbol{N}}_{,l}^jP_l)=0,
\end{align}
which implies that
\begin{align}
\delta\bar{\Pi} u=\delta\bar{\Pi} u-\delta\bar{\Pi} p=\delta\bar{\Pi}(u-p).\nonumber
\end{align}
On the reference element \(\widehat{\Omega^j}\) with unit size, under the affine mapping \(\boldsymbol{B}^j\) and by the Quotient Mapping Theorem \cite{18}, we have
\begin{align}
\parallel\widehat{\delta\bar{\Pi} u}\parallel_{H^1(\widehat{\Omega}^j)}&=\parallel\widehat{\delta\bar{\Pi} u-p}\parallel_{H^1(\widehat{\Omega}^j)}\nonumber \\
&\leq\parallel\widehat{\delta\bar{\Pi}}\parallel\inf_{\widehat{p}\in\widehat{\mathbb{P}}_1}\parallel\widehat{u}-\widehat{p}\parallel_{H^2(\widehat{\Omega}^j)}\nonumber \\
&\leq\varepsilon C_3\inf_{\widehat{p}\in\widehat{\mathbb{P}}_1}\parallel\widehat{u}-\widehat{p}\parallel_{H^2(\widehat{\Omega}^j)},\label{ep:64}
\end{align}
the last inequality is obtained by Lemma \ref{thm:bound2}. \\
Because \(\widehat{p}\in\widehat{\mathbb{P}}_1\), we have \(\lvert\widehat{u}-\widehat{p}\lvert^2_{H^2(\widehat{\Omega}^j)}=\lvert\widehat{u}\lvert^2_{H^2(\widehat{\Omega}^j)}\). \\
By the Quotient Mapping Theorem \cite{18}, we have \(\parallel\widehat{u}-\widehat{p}\parallel^2_{H^2(\widehat{\Omega}^j)}\leq (D_1)^2\lvert\widehat{u}\lvert^2_{H^2(\widehat{\Omega}^j)}\). \\
Combining this with Eq. \eqref{ep:64}:
\begin{align}
\parallel\widehat{\delta\bar{\Pi} u}\parallel^2_{H^1(\widehat{\Omega}^j)}\leq\varepsilon^2 C^2_3(D_1)^2\lvert\widehat{u}\lvert^2_{H^2(\widehat{\Omega}^j)}.\label{ep:68}
\end{align}
Combining this estimate with Eqs. \eqref{ep:25}-\eqref{ep:28},
\begin{align}
\lvert\delta\bar{\Pi} u\mid^2_{H^1(\widehat{\Omega}^j)}&\leq ({B_3})^2\parallel{(\mathbf{B}^j)}^{-1}\parallel^2 det(\boldsymbol{B}^j)^{1}\lvert\widehat{u}\lvert^2_{H^1(\widehat{\Omega}^j)}\nonumber\\
&\leq ({B_3})^2\parallel{(\mathbf{B}^j)}^{-1}\parallel^2 det(\boldsymbol{B}^j)^{1}\parallel\widehat{u}\parallel^2_{H^1(\widehat{\Omega}^j)}\nonumber\\
&\leq \varepsilon^2(C_3)^2(D_1)^2({B_3})^2\parallel{(\mathbf{B}^j)}^{-1}\parallel^2 det(\boldsymbol{B}^j)^{1}(\lvert\widehat{u}\lvert^2_{H^2(\widehat{\Omega}^j)}).\label{ep:69}
\end{align}
For \(\lvert\widehat{u}\lvert^2_{H^2(\widehat{\Omega}^j)}\), we have
\begin{align}
\lvert\widehat{u}\lvert^2_{H^2(\widehat{\Omega}^j)}&\leq ({B_2})^2\parallel\mathbf{B}^j\parallel^4det(\boldsymbol{B}^j)^{-1}\lvert u\lvert^2_{H^2({\Omega}^j)}.\label{ep:70}
\end{align}
For a rectangular mesh, there exists a constant \(C_4\) such that
\begin{align}
\parallel\mathbf{B}^j\parallel\parallel{(\mathbf{B}^j)}^{-1}\parallel\leq C_4,\label{ep:71}
\end{align}
and a constant \(D_2\) such that
\begin{align}
\parallel\mathbf{B}^j\parallel\leq D_2h.\label{ep:72}
\end{align}
Substituting Eqs. \eqref{ep:70}-\eqref{ep:72} into Eq. \eqref{ep:69}, we get
\begin{align}
\lvert\delta\bar{\Pi} u\lvert^2_{H^1({\Omega}^j)}\leq\varepsilon^2(C_3)^2(D_1)^2({B_3})^2({B_2})^2(C_4)^2(D_2)^2h^2\lvert u\lvert^2_{H^2({\Omega}^j)}.\nonumber
\end{align}
Defining \(C_{11}=C_3D_1{B_3}{B_2}C_4D_2\), we obtain
\begin{align}
\lvert\delta\bar{\Pi} u\lvert_{H^1({\Omega}^j)}\leq\varepsilon C_{11}h\lvert u\lvert_{H^2({\Omega}^j)}.\label{ep:74}
\end{align}
For \(\lvert \delta\bar{\Pi} u\lvert_{L^2({\Omega}^j)}\), we set the left hand side of Eq. \eqref{ep:68} to be \(\lvert \delta\bar{\Pi} u\lvert_{L^2({\Omega}^j)}\) and combine it with Eq. \eqref{ep:28}, which gives
\begin{align}
\lvert\delta\bar{\Pi} u\lvert_{L^2({\Omega}^j)}\leq\varepsilon C_{12}h^2\lvert u\lvert_{H^2({\Omega}^j)}.\nonumber
\end{align}
\end{proof}

Analogously to Theorem \ref{thm:T3m4}, we have the following theorem.
\begin{theorem} Let \(\Omega\subset\mathbb{R}^d, d=2,3\) be a domain equipped with a rectangular or cuboid mesh \(\mathcal{T}\). Let \(u\in H^2(\Omega)\) be the exact solution of a second-order linear elastic problem, where the governing equations are required to satisfy the properties given in Eqs. \eqref{ep:symmetric} - \eqref{ep:coercive}, and \(\widetilde{\widetilde{u}}_h\in H^1(\Omega)\) be its approximation by the improved PIML with interpolation operator \(\widetilde{\widetilde{\Pi}}\). Suppose that \(\varepsilon=\max_{\Omega_j\subset \mathcal{T}}\lvert\delta\bar{\boldsymbol{N}}^j\lvert_{\infty}\) is the generalization error for matrix learning. Then the following estimates hold:
\label{thm:T4m2}
\begin{align}
&\parallel u-\widetilde{\widetilde{u}}_h\parallel_{H^1(\Omega)}\leq (K_1+K_2\varepsilon)h\parallel u\parallel_{H^2(\Omega)},\label{ep:76}\\
&\parallel u-\widetilde{\widetilde{u}}_h\parallel_{L^2(\Omega)}\leq (K_3+K_4\varepsilon+K_5\varepsilon^2)h^2\parallel u\parallel_{H^2(\Omega)}.\label{ep:77}
\end{align}
\end{theorem}
\begin{proof}
Let \(\Omega^j\) denote the $j-$th element of \(\mathcal{T}\), by the definitions of \(\widetilde{\widetilde{\Pi}}\) in Eq.\eqref{ep:LPIMLop} and \(\bar{\Pi}\) in Eq. \eqref{ep:barNj}, both interpolation operators satisfy the linear completeness condition. Consider a subspace \(\bar{U}_h\) of \(V_h\) on each coarse element that contains \(P_1\). Then we have
\begin{align}
\parallel u-\widetilde{\widetilde{\Pi}}u\parallel_{H^1(\Omega)}\leq\parallel u-\bar{\Pi} u\parallel_{H^1(\Omega)}+\parallel\delta\bar{\Pi} u\parallel_{H^1(\Omega)},\nonumber
\end{align}
where \(\delta\bar{\Pi}=\bar{\Pi}-\widetilde{\widetilde{\Pi}}\), and by the Bramble–Hilbert Lemma, the first term satisfies
\begin{align}
\parallel u-\bar{\Pi} u\parallel_{H^1(\Omega)}\leq C_{0}h\parallel u\parallel_{H^2(\Omega)},\nonumber
\end{align}
where $h$ is the maximum element diameter in \(\mathcal{T}\). For the second term, by Eq. \eqref{ep:DL1} on each element \(\Omega^j\), we have
\begin{align}
\parallel\delta\bar{\Pi} u\parallel_{H^1(\Omega^j)}\leq\sqrt{\varepsilon^2 (C_{11}^j)^2h^2\lvert u\lvert^2_{H^2(\Omega^j)}+\varepsilon^2 (C^j_{12})^2h^4\lvert u\lvert^2_{H^2}}.\nonumber
\end{align}
Thus,
\begin{align}
\parallel\delta\bar{\Pi} u\parallel^2_{H^1(\Omega)}\leq\sum_{j}\varepsilon^2(C_{11}^j)^2h^2\lvert u\lvert^2_{H^2(\Omega^j)}+\varepsilon^2 (C^j_{12})^2h^4\lvert u\lvert^2_{H^2}.\label{ep:errorL1}
\end{align}
For sufficiently small $h$, the second term is smaller than the first, and since \(\mathcal{T}\) is rectangular, all constants in Eq. \eqref{ep:errorL1} are uniformly bounded, thus,
\begin{align}
\parallel\delta\bar{\Pi} u\parallel_{H^1(\Omega)}&\leq\varepsilon C_{11}h\lvert u\mid_{H^2(\Omega)}\nonumber\\
&\leq\varepsilon C_{11}h\parallel u\parallel_{H^2(\Omega)}.\nonumber
\end{align}
We obtain the \(H^1-\)norm estimate:
\begin{align}
\parallel u-\widetilde{\widetilde{u}}_h\parallel_{H^1(\Omega)}&\leq C_c\parallel u-\widetilde{\widetilde{\Pi}}u\parallel_{H^1(\Omega)}\nonumber\\
&\leq C_c(C_{0}+\varepsilon C_{11})h\parallel u\parallel_{H^2(\Omega)}.\nonumber
\end{align}
For the \(L^2-\)norm estimate, using the Aubin-Nitsche argument as in Theorem \ref{thm:T3m4}, we have
\begin{align}
\parallel u-\widetilde{\widetilde{u}}_h\parallel_{L^2(\Omega)}&\leq C_cC_7C_8(C_{0}+\varepsilon C_{11})^2h^2\parallel u\parallel_{H^2(\Omega)}.\nonumber
\end{align}
Let \(K_1=C_cC_0,\ K_2=C_cC_{11},\ K_3=C_cC_7C_8C_0^2,\ K_4=2C_cC_7C_8C_0C_{11}\) and \(K_5=C_cC_7C_8C_{11}^2\). This completes the proof.
\end{proof}
\begin{remark}Although the analysis above is presented for the case without boundary constraints, it extends naturally to linear or higher-order polynomial boundary constraints. Since the fixed part \(\boldsymbol{L}\) of \(\widetilde{\widetilde{\boldsymbol{N}}}^j_s\) already satisfies linear completeness, the analysis of \(\widetilde{\widetilde{\boldsymbol{N}}}^j_s\) is analogous to that in the case without boundary constraints.
\end{remark}
\begin{remark}All the results used in the preceding analysis hold in both two and three dimensions. Therefore, the analysis extends directly to three-dimensional problems, as illustrated by the numerical experiments below.
\end{remark}
\begin{remark}Although the above analysis is based on rectangular elements, it can be readily extended to shape-regular quadrilateral elements. A similar analysis is possible, provided that the shape functions of the fine elements on the quadrilateral mesh satisfy linear completeness, e.g., the bilinear basis functions Q4 or the spline basis functions L4 in \cite{chen2010spline}.
\end{remark}

\section{Numerical examples}\label{Sect. 5}
In this section, several examples of linear elasticity problems are presented to validate the proposed error analysis. We denote by MsFEM the theoretical substructure method, by RPIML the PIML with rigid-body displacement invariance, and by LPIML the PIML with linear completeness. For methods without boundary constraints, no suffix is added after their names, and the suffix "-l" indicates that linear boundary constraints are imposed. That is:
\begin{itemize}
  \item \textbf{MsFEM}: The theoretical substructure method without boundary constraints.
  \item \textbf{MsFEM-l}: The theoretical substructure method with linear boundary constraints.
  \item \textbf{LPIML}: The PIML with linear completeness without boundary constraints.
  \item \textbf{LPIML-l}: The PIML with linear completeness and linear boundary constraints.
  \item \textbf{RPIML}: The PIML with rigid-body displacement invariance without boundary constraints.
  \item \textbf{RPIML-l}: The PIML with rigid-body displacement invariance and linear boundary constraints.
\end{itemize}
In all of our examples, we refine the mesh by halving the element size, so the convergence order can be computed as
\begin{align}
\beta=\frac{\log{e_i}-\log{e_{i+1}}}{\log{2}},\nonumber
\end{align}
where \(e_i\) is the root-mean-square error on the nodes of the \(i-\)th refined mesh. The network architecture is shown in \autoref{fig:3}. The neural network consists of 12 layers with the following layer widths: [300, 500, 500, 500, 500, 500, 500, 500, 500, 500, 500, 500]. 200,000 samples were used to train the PIML for rectangles, whether with no boundary constraint or with a linear boundary constraint; 600,000 samples were used for quadrilaterals, with MSE as the loss function. All samples were obtained through uniform random sampling. Young's moduli of the coarse elements were sampled between \(1.0\times10^{-3}\) and \(1\). For the quadrilateral elements, the vertex coordinates are distributed in a unit disk and the interior angles must satisfy \(\frac{\pi}{6}\leq\theta\leq\frac{5\pi}{6}\). The model was trained for 20,000 epochs in PyTorch using the Adam optimizer with a learning rate of \(lr=0.001\) and used for FEM computation in MATLAB.

\subsection{Patch test}
\begin{figure}
    \centering
    \begin{subfigure}{0.45\textwidth}
        \centering
    		\includegraphics[trim=2cm 9cm 2cm 9cm, clip, width=0.8\textwidth]{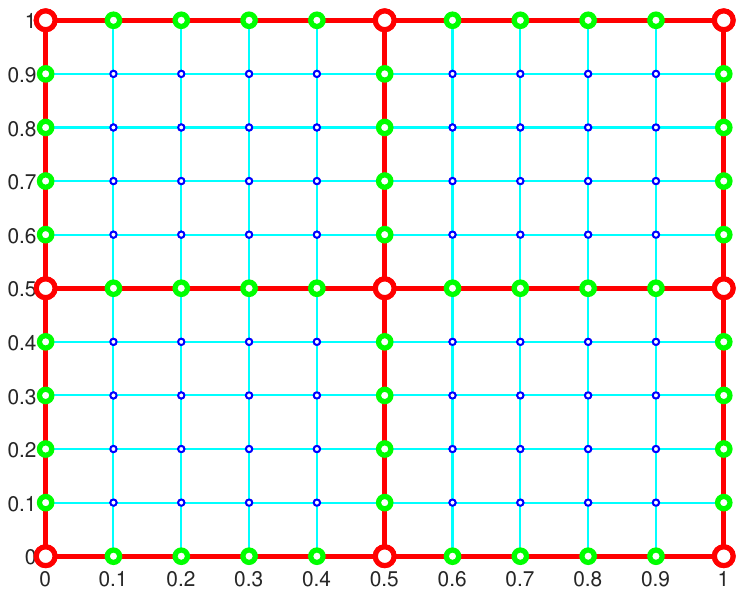}
    		\caption{The rectangular mesh for patch test.}
        \label{fig:sub:a}
    \end{subfigure}
    \hfill
    \begin{subfigure}{0.45\textwidth}
        \centering
    		\includegraphics[trim=2cm 9cm 2cm 9cm, clip, width=0.8\textwidth]{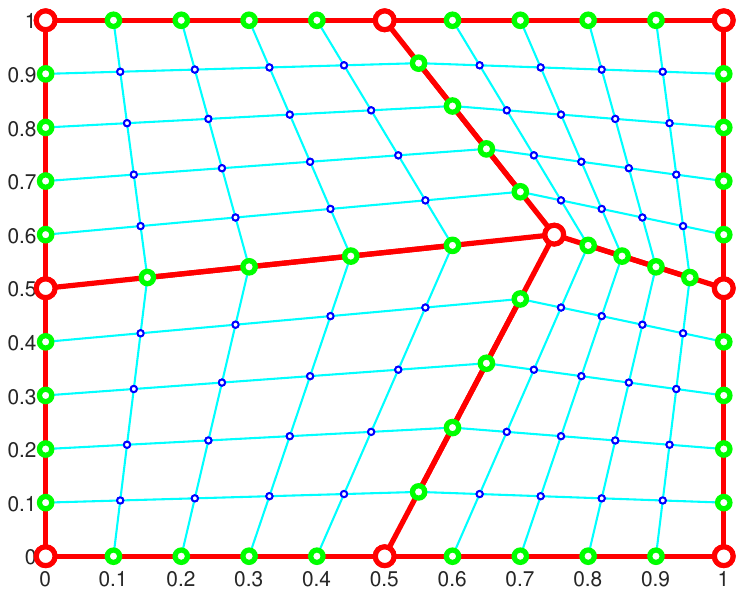}
    		\caption{The quadrilateral mesh for patch test.}
        \label{fig:sub:b}
    \end{subfigure}
    \caption{The mesh for patch test.}
    \label{fig:patchtest}
\end{figure}

This test examines the polynomial reproduction properties of RPIML and LPIML. The square domain is partitioned into \(2\times2\) coarse elements, and each contains \(5\times5\) fine elements, as shown in \autoref{fig:patchtest}.
The red dots denote the coarse element nodes and the green dots denote edge-interior nodes, and the blue dots denote interior nodes. For the case of linear boundary constraints, the fine element nodes are composed of both the blue dots and the green dots; for the case of no boundary constraints, the fine element nodes are composed only of the blue dots. In this experiment, Poisson's ratio is $\upsilon=0. 3$ and Young's modulus is 1000. Three displacement fields are defined as follows:
\begin{align}
&(Translation)\left\{
\begin{aligned}
    &u = 1, \\
    &v = 2.
\end{aligned}
\right. \nonumber\\
&(Rotation)\left\{
\begin{aligned}
    &u = y, \\
    &v = -x.
\end{aligned}
\right. \nonumber\\
&(Linear)\left\{
\begin{aligned}
    &u = x+y, \\
    &v = x+2y+\frac{1}{4}.
\end{aligned}
\right. \nonumber
\end{align}

The first two fields represent rigid-body motions, represented as Translation and Rotation, and the third is a general linear displacement field, denoted by “Linear”. The results are evaluated at both coarse-element nodes (denoted by “c”) and fine-element nodes (denoted by “f”). The method is regarded as passing the test only if the numerical errors on both sets of nodes satisfy the patch-test conditions. The results for rectangular meshes are reported in \autoref{tab:1}, where "P" denotes a pass and "F" denotes a failure. 

We can see that all methods pass the first two fields, because every method reproduces rigid-body motions. RPIML fails the linear field in both parts, but the coarse part of RPIML-l successfully passes the test. This behavior may be attributed to the imposed linear boundary constraints and the regularity of the underlying mesh. LPIML, LPIML-l, MsFEM, and MsFEM-l reproduce all three displacement fields, which shows that they have linear completeness.

In \autoref{tab:quatab}, we see similar results for the case with a quadrilateral mesh; the only difference is that, in this case, RPIML-l fails the patch test at the coarse element nodes.
\begin{table}[width=.9\linewidth,cols=4,pos=h]
  		\centering
  		\caption{Patch test for rectangular mesh: P for pass, F for failure}
  		\label{tab:example}
			\begin{tabular}{ccccccc}
  				\toprule
  				exact solution& $\mathrm{Translation{(c)}}$ & $\mathrm{Rotation{(c)}}$ & $\mathrm{Linear{(c)}}$& $\mathrm{Translation{(f)}}$ & $\mathrm{Rotation{(f)}}$ & $\mathrm{Linear{(f)}}$\\ 
 				 MsFEM & P & P & P& P & P & P \\
 				 MsFEM-l& P & P & P& P & P & P \\
 				 RPIML & P & P & \textcolor{blue}{F} & P & P & \textcolor{blue}{F}\\
 				 RPIML-l& P & P & P & P & P & \textcolor{blue}{F} \\
 				 LPIML& P & P & P & P & P & P \\
 				 LPIML-l& P & P & P & P & P & P \\
 				 \bottomrule
			\end{tabular}\label{tab:1}
\end{table}
\begin{table}[width=.9\linewidth,cols=4,pos=h]
  		\centering
  		\caption{Patch test for quadrilateral mesh: P for pass, F for failure}
  		\label{tab:example}
			\begin{tabular}{ccccccc}
  				\toprule
  				exact solution& $\mathrm{Translation{(c)}}$ & $\mathrm{Rotation{(c)}}$ & $\mathrm{Linear{(c)}}$& $\mathrm{Translation{(f)}}$ & $\mathrm{Rotation{(f)}}$ & $\mathrm{Linear{(f)}}$\\ 
 				 MsFEM-l& P & P & P& P & P & P \\
 				 RPIML-l& P & P & \textcolor{blue}{F} & P & P & \textcolor{blue}{F} \\
 				 LPIML-l& P & P & P & P & P & P \\
 				 \bottomrule
			\end{tabular}\label{tab:quatab}
\end{table}
\subsection{The homogeneous plane elasticity problem}
This example investigates whether enforcing only rigid-body displacement invariance may lead to error stagnation under mesh refinement. The solution domain \(\Omega=[0,1]\times[0,1]\) is a square domain discretized by a square coarse mesh, with each coarse element containing \(5\times5\) fine elements. The exact solution is
\begin{align}
\left\{
\begin{aligned}
    &u = e^x\cos y, \\
    &v = -e^x\sin y.
\end{aligned}
\right. \nonumber
\end{align}

\begin{figure}
    \centering
    \includegraphics[trim=0cm 6.5cm 0cm 7cm, clip, width=0.6\textwidth]{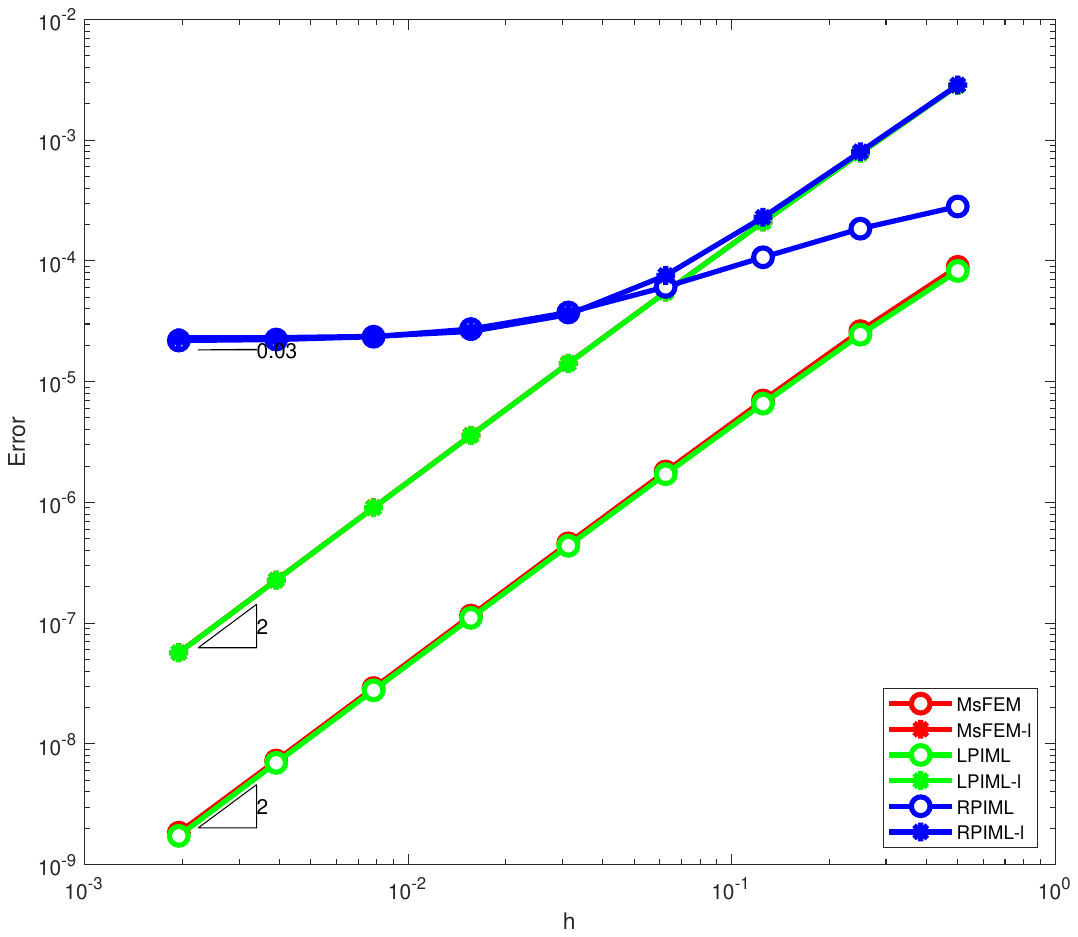}
    \caption{Convergence rates of the different methods in example 5.2}\label{fig:7}
\end{figure}

The errors over all coarse element nodes are computed. The results are presented in \autoref{fig:7} and \autoref{tab:2}. In \autoref{fig:7}, we plot reference lines to illustrate error trends. We observe that the slopes of the curves for RPIML and RPIML-l slowly decrease after several mesh refinements (blue curves). Similar results can be found in \autoref{tab:2}. The errors of the methods with rigid-body displacement invariance decrease much more slowly, and the RPIML error eventually stagnates as the mesh is refined. We also observe that LPIML and LPIML-l exhibit slopes comparable to those of MsFEM and MsFEM-l, indicating optimal convergence rates for this homogeneous problem.
\begin{table}[width=.9\linewidth,cols=4]
  		\centering
  		\caption{The homogeneous plane elasticity problem}
  		\label{tab:example}
			\begin{tabular}{l*{12}{c}}   
				\toprule
					& \multicolumn{12}{c}{Homogeneous Young’s modulus distribution}\\
				\cmidrule(lr){2-13}
				Mesh& MsFEM &$\beta$& MsFEM-l&$\beta$ &LPIML&$\beta$ & LPIML-l&$\beta$ & RPIML&$\beta$ & RPIML-l&$\beta$ \\
				\midrule
				$2\times2$ &8.92E-05&&2.84E-03&&8.25E-05&&2.83E-03&&2.82E-04&&2.86E-03&  \\
				$4\times4$ &2.62E-05&1.8&7.82E-04&1.9&2.44E-05&1.8&7.80E-04&1.9&1.85E-04&0.6&8.00E-04&1.8\\
				$8\times8$ &6.98E-06&1.9&2.10E-04&1.9&6.58E-06&1.9&2.10E-04&1.9&1.07E-04&0.8&2.30E-04&1.8  \\
				$16\times16$ &1.80E-06&2.0&5.51E-05&1.9&1.71E-06&1.9&5.50E-05&1.9&6.09E-05&0.8&7.58E-05&1.6   \\
				$32\times32$ &4.58E-07&2.0&1.42E-05&2.0&4.37E-07&2.0&1.41E-05&2.0&3.74E-05&0.7&3.59E-05&1.1   \\
				$64\times64$ &1.15E-07&2.0&3.59E-06&2.0&1.10E-07&2.0&3.58E-06&2.0&2.72E-05&0.5&2.60E-05&0.5 \\
				$128\times128$ &2.90E-08&2.0&9.05E-07&2.0&2.78E-08&2.0&9.03E-07&2.0&2.35E-05&0.2&2.37E-05&0.1  \\
				$256\times256$ &7.28E-09&2.0&2.27E-07&2.0&6.97E-09&2.0&2.27E-07&2.0&2.23E-05&0.1&2.31E-05&0.0  \\
				$512\times512$ &1.84E-09&2.0&5.69E-08&2.0&1.73E-09&2.0&5.68E-08&2.0&2.19E-05&0.0&2.30E-05&0.0  \\
\bottomrule
\end{tabular}\label{tab:2}
\end{table}
\FloatBarrier
\subsection{The heterogeneous materials problem}
One advantage of the PIML approach is its computational efficiency for heterogeneous-material problems. In this experiment we use the same mesh as the previous experiment to solve two heterogeneous material problems. The first case involves a low-frequency material distribution, whereas the second involves a high-frequency distribution. For the low-frequency case, we prescribe the exact solution:
\begin{align}
\left\{
\begin{aligned}
    &u = e^{(x-y)}, \\
    &v = x\sin{y},
\end{aligned}
\right. \nonumber
\end{align}
with Young’s modulus given by
\begin{align}
E(x, y)=\frac{e^{(x+y)}}{e^2}.\label{ep:low}
\end{align}
For the high-frequency distribution, we conduct another experiment with an exact solution:
\begin{align}
\left\{
\begin{aligned}
    &u = \sin{x}, \\
    &v = e^y,
\end{aligned}
\right. \nonumber
\end{align}
with Young's modulus
\begin{align}
E(x, y)=\frac{\sin(100\pi{x})\cos(80\pi{y})+4}{5}.\label{ep:high}
\end{align}
For both problems, we set Poisson’s ratio to $\upsilon=0.3$. The boundary conditions and body force are derived from the prescribed exact solution and Young’s modulus.

\begin{figure}[h]
    \centering
    \begin{subfigure}{0.45\textwidth}
        \centering
        \includegraphics[trim=2cm 9cm 2cm 9cm, clip, width=1\textwidth]{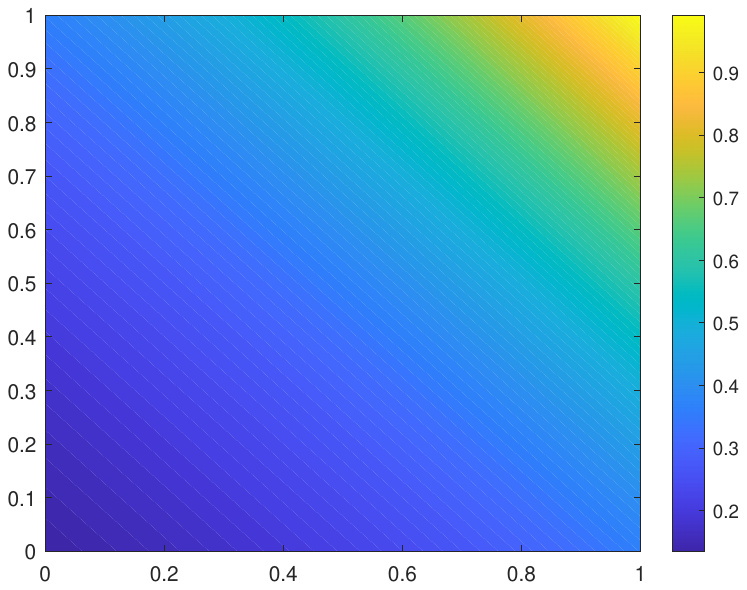}
        \caption{Young's modulus with low-frequency Eq. \eqref{ep:low}}
        \label{fig:sub:a}
    \end{subfigure}
    \hfill
    \begin{subfigure}{0.45\textwidth}
        \centering
        \includegraphics[trim=2cm 9cm 2cm 9cm, clip, width=1\textwidth]{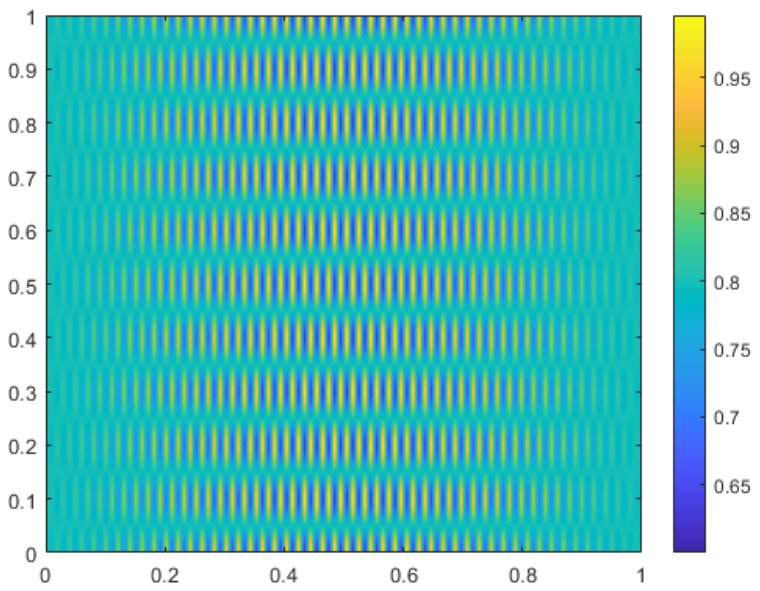}
        \caption{Young's modulus with high-frequency Eq. \eqref{ep:high}}
        \label{fig:sub:b}
    \end{subfigure}
    \caption{The heatmaps of two cases of Young's modulus}
    \label{fig:frquency}
\end{figure}
\begin{figure}
    \centering
    \begin{subfigure}{0.45\textwidth}
        \centering
        \includegraphics[trim=2cm 9cm 2cm 9cm, clip, width=1\textwidth]{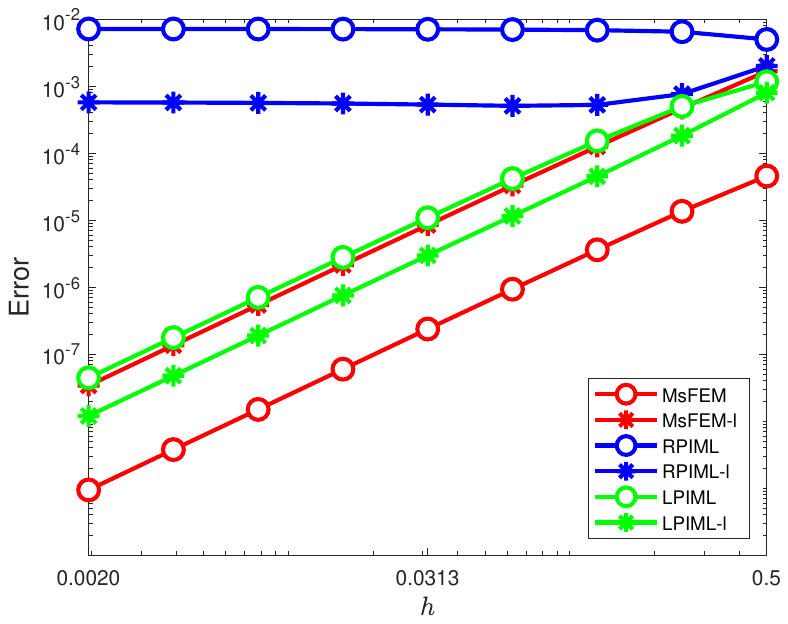}
        \caption{Errors for the case with low-frequency Young’s modulus distribution}
        \label{fig:sub:a}
    \end{subfigure}
    \hfill
    \begin{subfigure}{0.45\textwidth}
        \centering
        \includegraphics[trim=2cm 9cm 2cm 9cm, clip, width=1\textwidth]{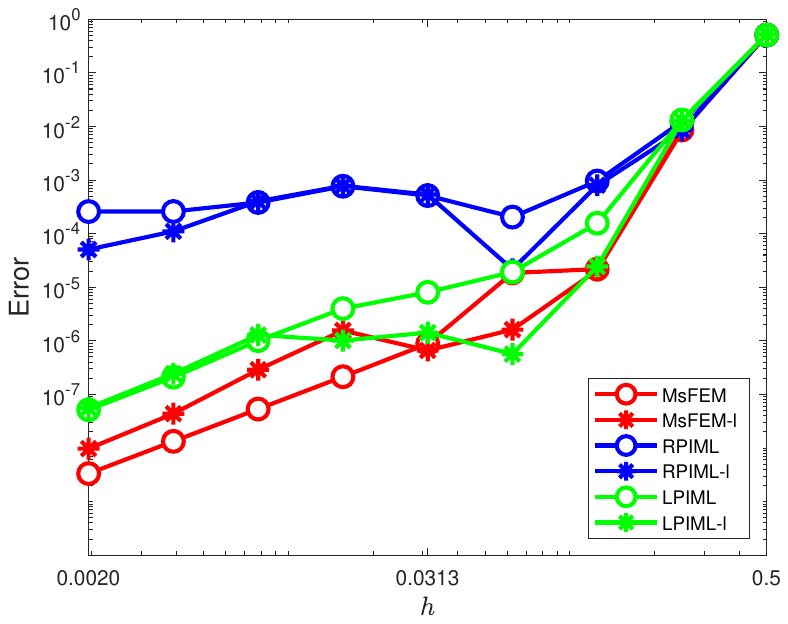}
        \caption{Errors for the case with high-frequency Young’s modulus distribution}
        \label{fig:sub:b}
    \end{subfigure}
    \caption{Errors for the two heterogeneous-material problems}
    \label{fig:frequenconvergence}
\end{figure}

The heatmaps of the two cases of Young's modulus with various frequencies are shown in \autoref{fig:frquency}. The root-mean-square errors at the coarse-element nodes by the six methods are reported in \autoref{tab:3} and \autoref{tab:4}, and the relationship between numerical errors and the element size is shown in \autoref{fig:frequenconvergence}.

The results indicate that PIML with linear completeness achieves higher accuracy than the PIML with rigid-body displacement invariance constraints for both low- and high-frequency material distributions. From \autoref{fig:frequenconvergence}, we can see that both LPIML and LPIML-l exhibit observed convergence rates comparable to those of MsFEM, whereas both RPIML and RPIML-l show a reduced convergence rate, and RPIML even ceases to converge. Meanwhile, we also observe that throughout the entire computation process, the errors of LPIML and LPIML-l are smaller than those of RPIML and RPIML-l. In the low-frequency case, the error of LPIML-l is even less than that of MsFEM-l.
\begin{table}[width=.9\linewidth,cols=4,pos=h]
\centering
\caption{Errors and convergence orders for the low-frequency Young’s modulus distribution}
\label{tab:deflection}
\begin{tabular}{l*{12}{c}}   
\toprule
& \multicolumn{12}{c}{Lower frequency}\\
\cmidrule(lr){2-13}
Mesh& MsFEM & $\beta$ & MsFEM-l &$\beta$ &LPIML &$\beta$ & LPIML-l &$\beta$ & RPIML & $\beta$& RPIML-l &$\beta$  \\
\midrule
$2\times2$ &4.57E-05& & 1.68E-03& &1.17E-03 & & 8.01E-04 & &5.00E-03 & &2.01E-03 &   \\
$4\times4$ &1.36E-05& 1.2 & 4.68E-04 & 1.8 & 5.01E-04 & 1.2 & 1.83E-04 & 2.1 &6.52E-03&-0.4&7.65E-04&1.4    \\
$8\times8$ &3.62E-06& 1.9 & 1.26E-04 & 1.9 & 1.53E-04 & 1.7 & 4.53E-05 & 2.0 &6.88E-03&-0.1&5.27E-04&0.5    \\
$16\times16$ &9.34E-07& 2.0 & 3.30E-05 & 1.9 & 4.19E-05 & 1.9 & 1.16E-05 & 2.0 &7.02E-03 &-0.0&5.10E-04&0.1    \\
$32\times32$ &2.38E-07& 2.0 & 8.49E-06 & 2.0 & 1.09E-05 & 1.9 & 2.98E-06 & 2.0 &7.09E-03&-0.0&5.33E-04&-0.1    \\
$64\times64$ &5.99E-08& 2.0 & 2.15E-06 & 2.0 & 2.79E-06 & 2.0 & 7.55E-07 & 2.0 &7.13E-03&-0.0&5.53E-04&-0.1    \\
$128\times128$ &1.50E-08& 2.0 & 5.43E-07 & 2.0 & 7.06E-07 & 2.0 & 1.90E-07 & 2.0 &7.15E-03&0.0&5.65E-04&-0.0    \\
$256\times256$ &3.77E-09& 2.0 & 1.36E-07 & 2.0 & 1.77E-07 & 2.0 & 4.77E-08 & 2.0 &7.16E-03&0.0&5.72E-04&-0.0    \\
$512\times512$ &9.43E-10& 2.0 &3.41E-08 & 2.0 & 4.45E-08 & 2.0 & 1.20E-08 & 2.0 &7.16E-03&0.0&5.75E-04&-0.0    \\
\bottomrule
\end{tabular}\label{tab:3}
\end{table}
\begin{table}[width=.9\linewidth,cols=4,pos=h]
\centering
\caption{Errors and convergence orders for the high-frequency Young’s modulus distribution}
\label{tab:deflection}
\begin{tabular}{l*{12}{c}}   
\toprule
& \multicolumn{12}{c}{Higher frequency}\\
\cmidrule(lr){2-13}
Mesh& MsFEM & $\beta$ & MsFEM-l &$\beta$ &LPIML &$\beta$ & LPIML-l &$\beta$ & RPIML &$\beta$ & RPIML-l &$\beta$  \\
\midrule
$2\times2$ &5.08E-01& &5.21E-01& &5.08E-01& &5.21E-01& &5.08E-01& &5.21E-01 &   \\
$4\times4$ &8.60E-03&5.9&8.21E-03&6.0&1.32E-02&5.3& 1.28E-02 &5.3&1.24E-02&5.4&8.45E-03&6.0   \\
$8\times8$ &2.18E-05 &8.6&2.07E-05&8.6&1.57E-04&6.4& 2.44E-05 &9.0&9.60E-04&3.7&8.00E-04&3.4 \\
$16\times16$ &1.85E-05&0.2&1.62E-06&3.7&1.93E-05&3.0& 5.68E-07 &5.4&2.04E-04&2.2&2.16E-05&5.2    \\
$32\times32$ &8.83E-07&4.4&6.61E-07&1.3&8.05E-06&1.3 &1.42E-06 &-1.3&5.14E-04&-1.3&5.36E-04&-4.6    \\
$64\times64$ &2.14E-07&2.0&1.56E-06&-1.2&3.99E-06&1.0&1.01E-06 &0.5&7.67E-04&-0.6&7.82E-04&-0.5 \\
$128\times128$ &5.34E-08&2.0&2.89E-07&2.4&1.01E-06&2.0 &1.29E-06&-0.4&3.82E-04&1.0&4.02E-04&1.0  \\
$256\times256$ &1.33E-08&2.0&4.38E-08&2.7&2.10E-07&2.3&2.44E-07&2.4&2.59E-04&0.6&1.10E-04&1.9   \\
$512\times512$ &3.34E-09&2.0&9.75E-09&2.2&5.17E-08&2.0&5.63E-08&2.1&2.59E-04&0.0&5.06E-05&1.1   \\
\bottomrule
\end{tabular}\label{tab:4}
\end{table}
\FloatBarrier
\subsection{Cook's problem with heterogeneous materials}
\begin{figure}[pos=h]
    \centering
    \includegraphics[width=0.45\textwidth]{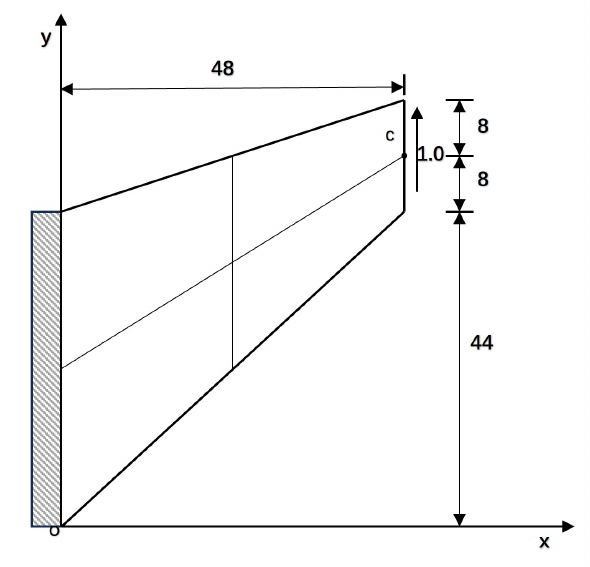}
    \caption{Geometry and mesh for Cook’s problem}
    \label{fig:8}
\end{figure}

For an arbitrary quadrilateral mesh, the numerical shape functions depend on both material distribution and geometry of elements. Therefore, we trained a new neural network including the coordinates of the coarse element vertices as input variables for quadrilateral substructures to approximate the numerical shape functions. In this example, we consider Cook’s problem and report experimental results only for substructures with linear boundary constraints.  The experimental domain, shown in \autoref{fig:8}, is discretized using a sequence of quadrilateral meshes obtained by successively halving the element size. The Young’s modulus distribution of the problem is
\begin{align}
    E(x, y)=\frac{e^{(x-y)}}{e^2}.\nonumber
\end{align}
The left boundary is clamped and a uniformly distributed shear traction is applied to the right boundary. We use the finite element solution on a $2560\times2560$ mesh as the reference solution and compare the values of the displacement component \(v\) at point C obtained by three methods with linear boundary constraints. The errors and convergence rates are reported in \autoref{tab:5}. We can see that as the mesh is refined, the values obtained by MsFEM-l and LPIML-l approach the reference value, whereas RPIML-l fails to provide a convergent approximation.
\begin{table}[width=.9\linewidth,cols=4,pos=h]
\centering
\caption{Errors for Cook’s problem}
\begin{tabular}{l*{6}{c}}   
\toprule
& \multicolumn{6}{c}{The errors of $v_c$}\\
\cmidrule(lr){2-7}
Mesh& MsFEM-l&$\beta$&LPIML-l&$\beta$&RPIML-l&$\beta$\\
\midrule
$2\times2$ &1.3897E+01&&1.4983E+01&&1.4546E+01&\\
$4\times4$ &5.3917E+00&1.37&7.1529E+00&1.07&9.3057E+00&0.64\\
$8\times8$ &1.6934E+00&1.87&2.8897E+00&1.31&9.5357E+00&-0.04\\
$16\times16$ &5.0021E-01&1.76&9.9750E-01&1.53&1.0866E+01&-0.19\\
$32\times32$ &1.5293E-01&1.71&3.1400E-01&1.67&1.1821E+01&-0.12\\
$64\times64$ &4.9513E-02&1.63&9.9500E-02&1.66&1.2300E+01&-0.06\\
$128\times128$ &1.6820E-02&1.56&3.2700E-02&1.61&1.2511E+01&-0.02\\
$256\times256$ &5.8551E-03&1.52&1.1200E-02&1.55&1.2603E+01&-0.01\\
\midrule
Reference & \multicolumn{6}{c}{$35.91$}\\
\bottomrule
\end{tabular}\label{tab:5}
\end{table}
\FloatBarrier
\subsection{Annular domain with a clamped inner boundary}
\begin{figure}[pos=h]
    \centering
    \includegraphics[width=0.4\textwidth]{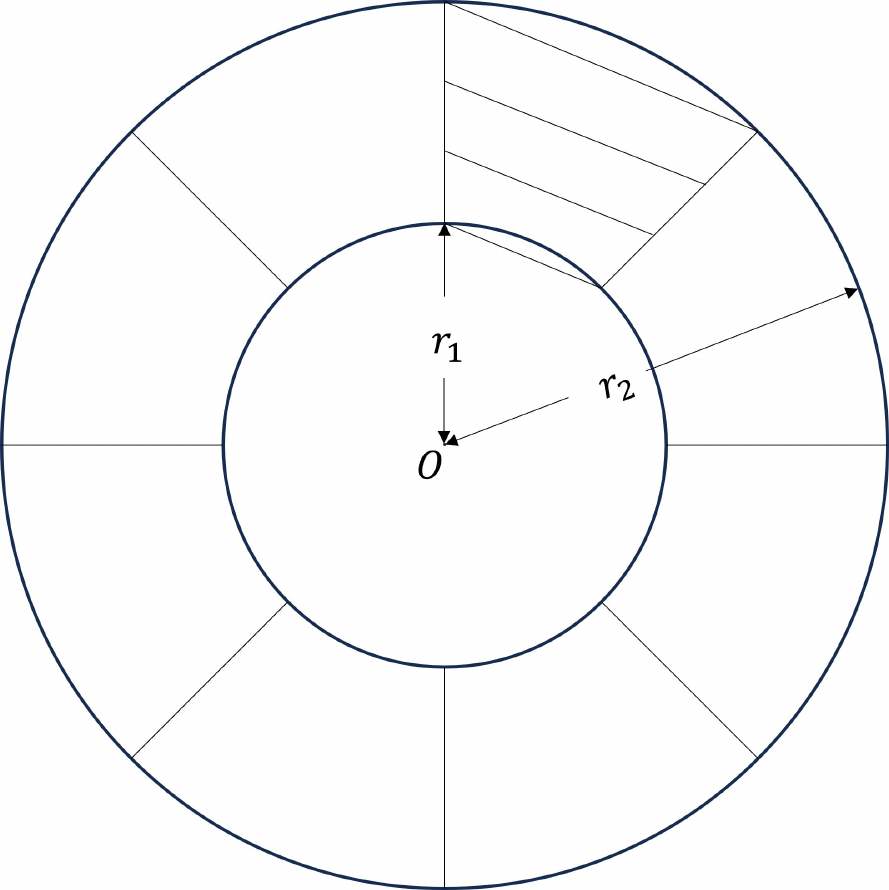}
    \caption{Mesh of the annular domain}
    \label{fig:ring}
\end{figure}

\begin{figure}[pos=h]
    \centering
    \begin{subfigure}{0.45\textwidth}
        \centering
        \includegraphics[trim=4cm 10cm 4cm 10cm, clip, width=1\textwidth]{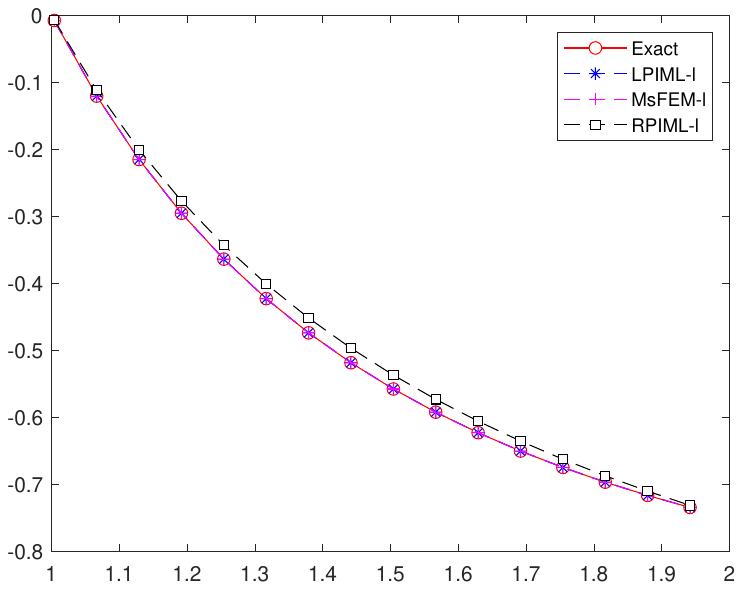}
        \caption{The displacement $u$ at $y=0$ }
        \label{fig:sub:a}
    \end{subfigure}
    \hfill
    \begin{subfigure}{0.45\textwidth}
        \centering
        \includegraphics[trim=4cm 10cm 4cm 10cm, clip, width=1\textwidth]{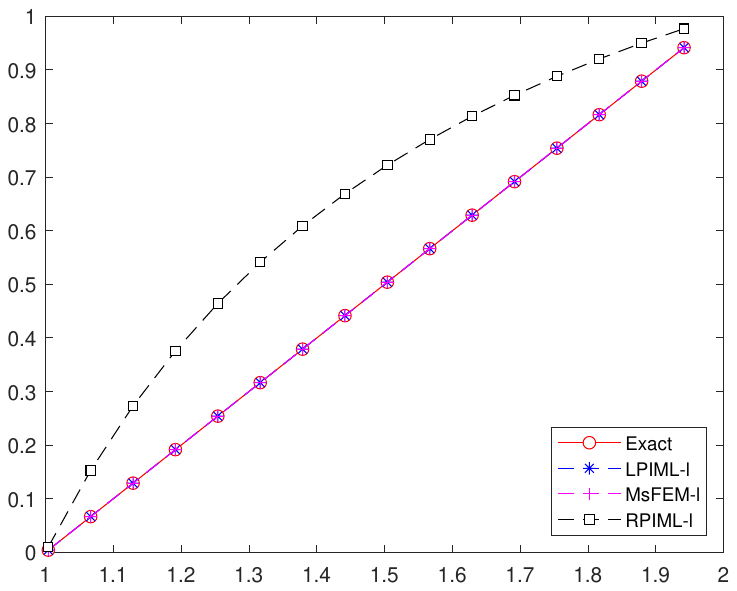}
        \caption{The displacement $v$ at $y=0$}
        \label{fig:sub:b}
    \end{subfigure}
    \caption{Exact and numerical solutions along \(y\equiv0\) in the annular-domain}
    \label{fig:chart}
\end{figure}

In this experiment, we consider an annular domain as shown in \autoref{fig:ring}. The radii of the inner and outer rings are denoted by \(r_1\) and \(r_2\). A displacement field that vanishes on the inner ring is given by
\begin{align}
    u=\frac{1}{r}-\frac{1}{r_1^2},\  v=r-r_1^2.\nonumber
\end{align}
Young's modulus is $E=\frac{r}{r_2}$,
where \(r=\sqrt{x^2+y^2}\) and \(r_1=1,r_2=2\) in our problem. The boundary conditions and body force are derived from the prescribed exact solution and Young’s modulus. The subdivision of a part of the circular domain is illustrated in \autoref{fig:ring}. Because of symmetry, the remaining angular sectors are discretized in the same manner. The mesh is refined by halving the elements along the radial and circumferential directions. The root-mean-square error over all coarse-element nodes is reported in \autoref{tab:ring}. The results show that LPIML-l and MsFEM-l both exhibit the expected convergence behavior,  while RPIML-l fails. \autoref{fig:chart} compares the exact solution and the numerical solution of \(u\) and \(v\) along the positive half of the \(x-\)axis.

\begin{table}[width=.9\linewidth,cols=4,pos=h]
\centering
\caption{Errors and convergence orders for the annular-domain problem}
\label{tab:ring}
\begin{tabular}{l*{6}{c}}   
\toprule
& \multicolumn{6}{c}{The errors of the ring}\\
\cmidrule(lr){2-7}
Mesh& MsFEM-l & \(\beta\) & LPIML-l &\(\beta\) &RPIML-l &\(\beta\)   \\
\midrule
$4\times8$ &2.4534E-02&  & 2.1152E-02& & 2.1422E-01 &      \\
$8\times16$ &5.7070E-03& 2.10 & 5.4165E-03& 1.97 & 1.5100E-01& 0.50     \\
$16\times32$ &1.3767E-03& 2.05 & 1.3441E-03 & 2.01 & 1.4604E-01& 0.05    \\
$32\times64$ &3.3843E-04& 2.02 & 3.3282E-04 & 2.01 & 1.5009E-01& -0.04    \\
$64\times128$ &8.3923E-05& 2.01 & 8.2669E-05 & 2.01 & 1.5298E-01& -0.03    \\
$128\times256$ &2.0897E-05& 2.01 & 2.0591E-05 & 2.01 & 1.5462E-01& -0.02   \\
$256\times512$ &5.2140E-06& 2.00 & 5.1376E-06 & 2.00 & 1.5554E-01& -0.01   \\
\bottomrule
\end{tabular}
\end{table}
\FloatBarrier
\subsection{Heart-shaped domain with a star-shaped mesh}
\begin{figure}[pos=h]
    \centering
    \includegraphics[trim=4cm 9cm 4cm 10cm, clip, width=0.6\textwidth]{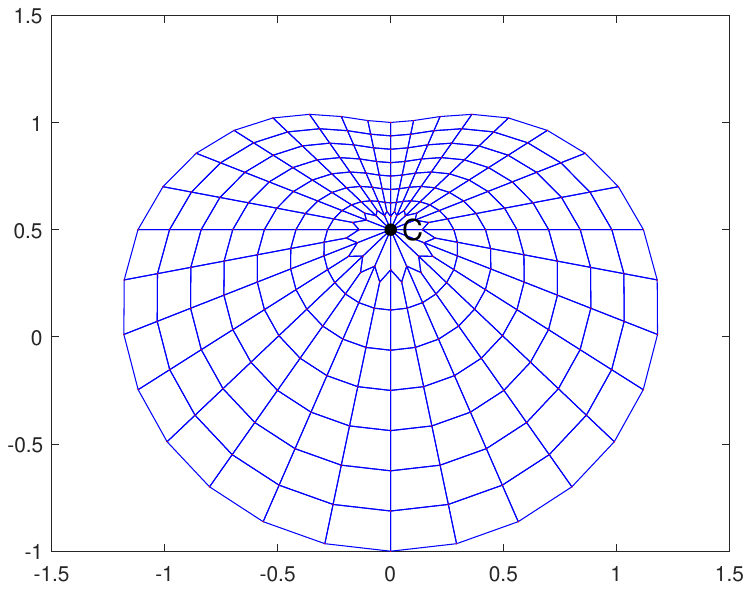}
    \caption{The discretization of the heart‑shaped domain}
    \label{fig:heart}
\end{figure}

\begin{figure}[pos=h]
    \centering
    \begin{subfigure}{0.45\textwidth}
        \centering
        \includegraphics[trim=4cm 10cm 4cm 10cm, clip, width=0.8\textwidth]{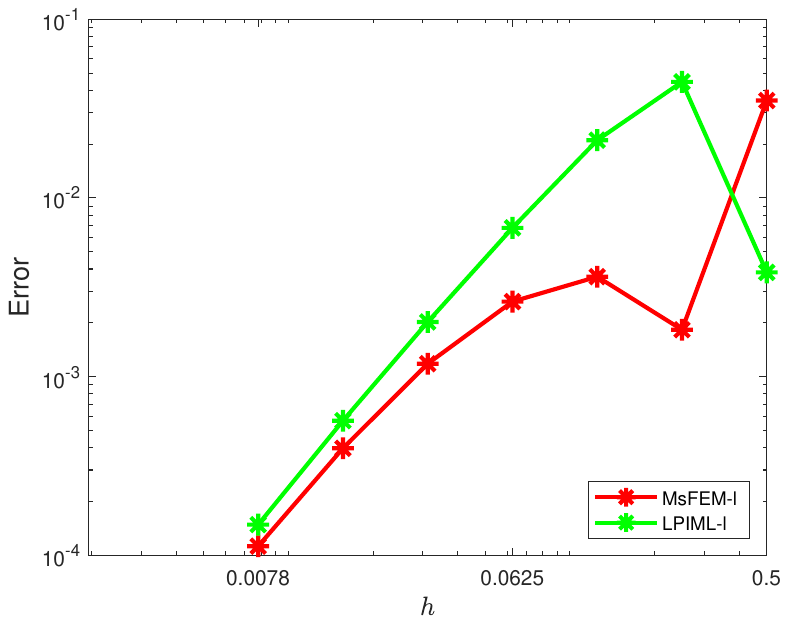}
        \caption{The errors of $u$ at point \(C\) }
        \label{fig:sub:a}
    \end{subfigure}
    \hfill
    \begin{subfigure}{0.45\textwidth}
        \centering
        \includegraphics[trim=4cm 10cm 4cm 10cm, clip, width=0.8\textwidth]{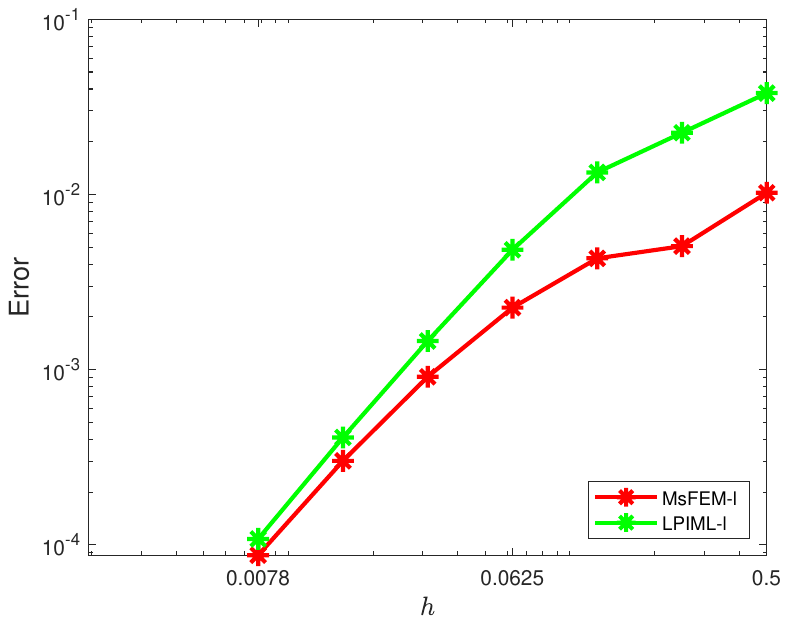}
        \caption{The errors of $v$ at point \(C\)}
        \label{fig:sub:b}
    \end{subfigure}
    \caption{Errors at point \(C\)}
    \label{fig:UVconvergence}
\end{figure}

The PIML training meshes satisfy certain shape-regularity conditions: the interior angles of each element must be greater than \(\frac{\pi}{6}\) and less than \(\frac{5\pi}{6}\). In this example, we consider a heart-shaped domain discretized by a star-shaped mesh radiating outward from point \(C\) with coordinates \((0,\frac{1}{2})\) as shown in \autoref{fig:heart}. We observe that, as the mesh is refined, the interior angles of all the elements near point \(C\) decrease. This example is therefore used to test the generalization ability and numerical stability of the PIML model. We impose clamped boundary conditions. The distribution of Young's modulus and body force are
\begin{align}
    E(x,y)=\frac{\sqrt{x^2+y^2+1}}{3},\ \boldsymbol{f}(x,y)=[1,-1]^T.\nonumber
\end{align}
Since no exact solution is available, we use a reference solution obtained by a \(512\times2048\) mesh where each element is further subdivided into two triangular elements. We report only the displacement at point \(C\). The errors of the two methods relative to the reference solution are shown in \autoref{tab:heart} and \autoref{fig:UVconvergence}. 

\begin{figure}[pos=h]
    \centering
    \begin{subfigure}[b]{0.45\textwidth}
        \centering
        \includegraphics[trim=4cm 9cm 4cm 9cm, clip,width=\textwidth]{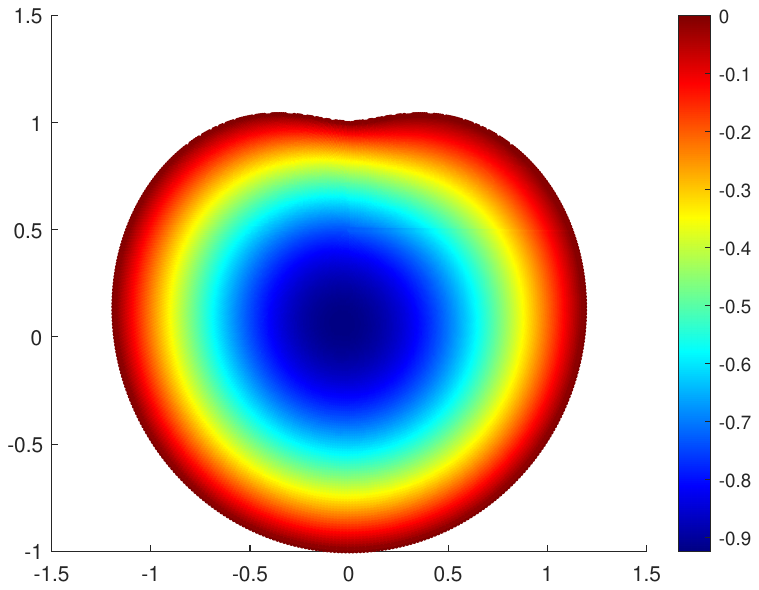}
        \caption{Displacement component \(u\) computed by MsFEM-l}
        \label{fig:sub1}
    \end{subfigure}
    \hfill  
    \begin{subfigure}[b]{0.45\textwidth}
        \centering
        \includegraphics[trim=4cm 9cm 4cm 9cm, clip,width=\textwidth]{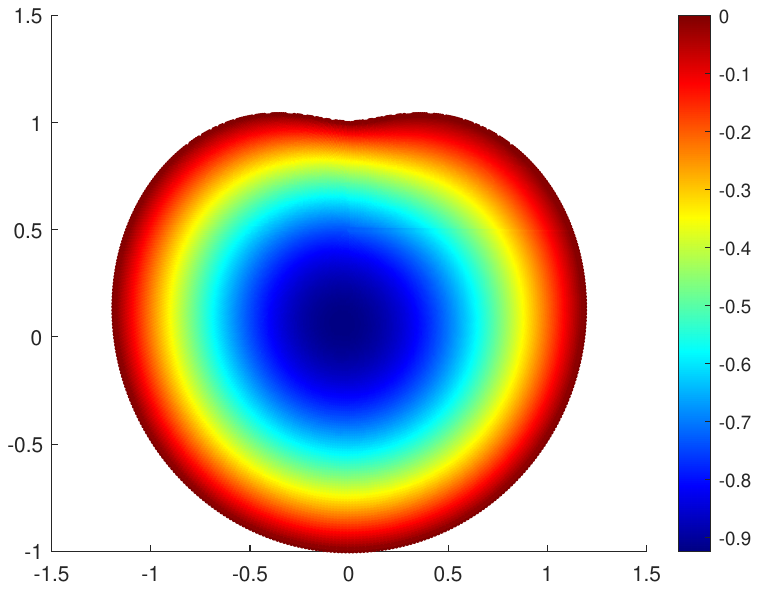}
        \caption{Displacement component \(u\) computed by LPIML-l}
        \label{fig:sub2}
    \end{subfigure}
    \begin{subfigure}[b]{0.45\textwidth}
        \centering
        \includegraphics[trim=4cm 9cm 4cm 9cm, clip,width=\textwidth]{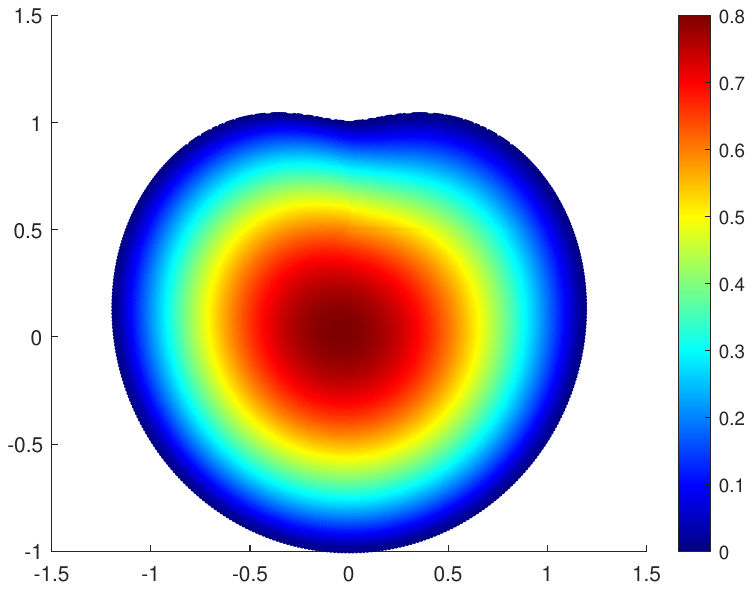}
        \caption{Displacement component \(v\) computed by MsFEM-l}
        \label{fig:sub3}
    \end{subfigure}
    \hfill
    \begin{subfigure}[b]{0.45\textwidth}
        \centering
        \includegraphics[trim=4cm 9cm 4cm 9cm, clip,width=\textwidth]{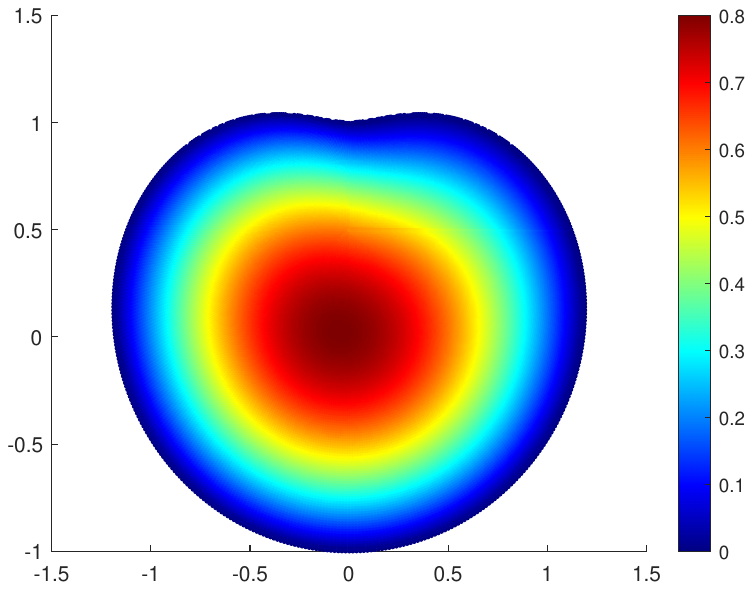}
        \caption{Displacement component \(v\) computed by LPIML-l}
        \label{fig:sub4}
    \end{subfigure}
    \caption{Numerical displacement fields in example 5.6}
    \label{fig:fourgrid}
\end{figure}

The results indicate that LPIML-l maintains good generalization performance and numerical stability even for severely distorted meshes. The displacement fields \(u\) and \(v\) obtained by MsFEM-l and LPIML-l are shown in \autoref{fig:fourgrid}.

\begin{table}[width=.9\linewidth,cols=4,pos=h]
\centering
\caption{Displacement errors at point \(C\)}
\label{tab:deflection}
\begin{tabular}{l*{8}{c}}   
\toprule
& \multicolumn{4}{c}{The errors of $u$}& \multicolumn{4}{c}{The errors of $v$}\\
\cmidrule(lr){2-5}
\cmidrule(lr){6-9}
Mesh& \(u_{MsFEM-l}\)&\(\beta\)&\(u_{LPIML-l}\)&\(\beta\)& \(v_{MsFEM-l}\)&\(\beta\)&\(v_{LPIML-l}\)&\(\beta\)\\
\midrule
$2\times8$ &3.5092E-02&&3.8313E-03&&1.0227E-02&&3.7904E-02&\\
$4\times16$ &1.8285E-03&4.26&4.4632E-02&-3.54&5.0699E-03&1.01&2.2498E-02&0.75\\
$8\times32$ &3.6224E-03&-0.99&2.1049E-02&1.08&4.3185E-03&0.23&1.3357E-02&0.75\\
$16\times64$ &2.6282E-03&0.46&6.7991E-03&1.63&2.2583E-03&0.94&4.8323E-03&1.47\\
$32\times128$ &1.1797E-03&1.16&2.0219E-03&1.75&9.1017E-04&1.31&1.4567E-03&1.73\\
$64\times256$ &3.9677E-04&1.57&5.6613E-04&1.84&3.0147E-04&1.59&4.0999E-04&1.83\\
$128\times512$ &1.1233E-04&1.82&1.4834E-04&1.93&8.7199E-05&1.79&1.0802E-04&1.92\\
\midrule
Reference & \multicolumn{8}{c}{$u=-0.7280,v=0.5900$}\\
\bottomrule
\end{tabular}\label{tab:heart}
\end{table}
We observe that the numerical solutions do not exhibit pronounced error concentration and retain reasonable convergence rates at point \(C\).
\FloatBarrier
\subsection{The three-dimensional cantilever}
\begin{figure}[pos=h]
\centering
    \includegraphics[width=0.8\textwidth]{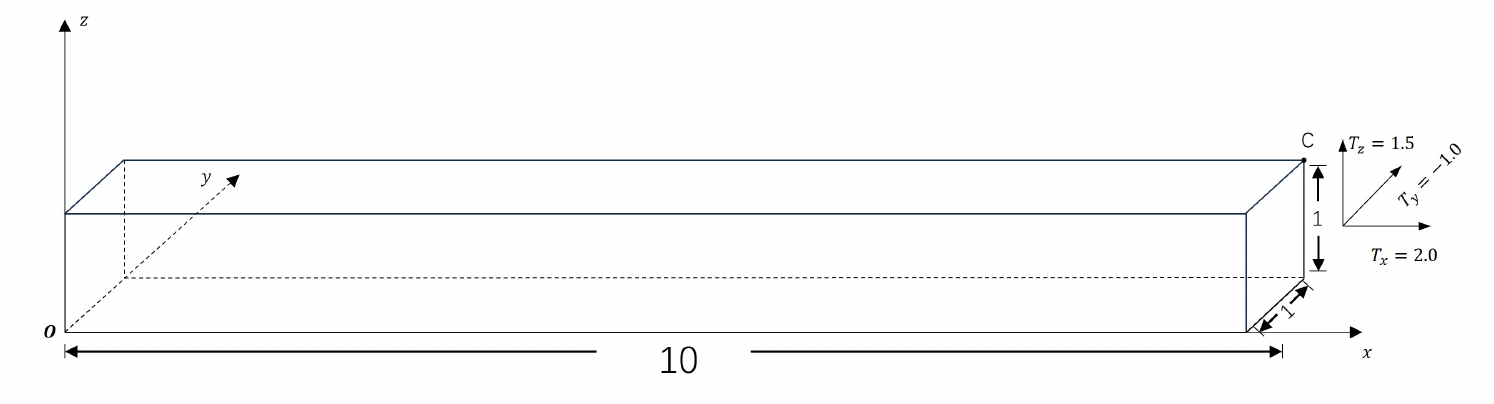}
    \caption{Geometry of the three‑dimensional cantilever}
    \label{fig:cantilever}
\end{figure}

In order to compare the computational costs of the three methods, we consider the three-dimensional cantilever shown in \autoref{fig:cantilever}. The left end of the cantilever is clamped, whereas the right end is free. Three traction components are applied to the right end:
\begin{align*}
T_x&=2.0,\\
T_y&=-1.0,\\
T_z&=1.5.
\end{align*}
We set Poisson’s ratio to 0.3, and Young's modulus is given by
\begin{align}
E(x,y,z)=\frac{e^{z+y}(\sin({10\pi{x}})+1.5)}{e^2+3}.\nonumber
\end{align}
Since no exact solution is available, we use Richardson extrapolation based on the solutions of the three-dimensional linear element \(H_8\) computed on \(200\times20\times20\) and \(400\times40\times40\) meshes. This yields a reference solution with an accuracy comparable to that obtained on an  \(800\times80\times80\) mesh. Every coarse element is divided into \(5\times5\times5\) fine elements. We compare the errors of MsFEM-l, LPIML-l, and RPIML-l for different numbers of coarse elements. The errors in the three displacement components \(u,v,w\) at point \(C\) are reported in \autoref{tab:cantilever}. The results show that LPIML-l and RPIML-l achieve accuracy comparable to that of MsFEM-l in this test.\\
\begin{table}[width=.9\linewidth,cols=4,pos=h]
\centering
\caption{Errors for the three-dimensional cantilever problem}
\begin{tabular}{l*{6}{c}}   
\toprule
& \multicolumn{6}{c}{The errors of $u_c$}\\
\cmidrule(lr){2-7}
Mesh& MsFEM-l&$\beta$&LPIML-l&$\beta$&RPIML-l&$\beta$\\
\midrule
$20\times2\times2$ &6.521E-01&&6.195E-01&&6.505E-01&\\
$40\times4\times4$ &2.604E-01&1.32&2.547E-01&1.28&2.598E-01&1.32\\
$80\times8\times8$ &7.291E-02&1.84&7.027E-02&1.86&7.281E-02&1.84\\
$160\times16\times16$ &1.995E-02&1.87&1.920E-02&1.87&1.996E-02&1.87\\
\midrule
Reference & \multicolumn{6}{c}{$1.63$}\\
\bottomrule
& \multicolumn{6}{c}{The errors of $v_c$}\\
\cmidrule(lr){2-7}
Mesh& MsFEM-l&$\beta$&LPIML-l&$\beta$&RPIML-l&$\beta$\\
\midrule
$20\times2\times2$ &1.049E+00&&1.146E+00&&1.052E+00&\\
$40\times4\times4$ &6.090E-01&0.78&6.290E-01&0.86&6.112E-01&0.78\\
$80\times8\times8$ &2.311E-01&1.40&2.367E-01&1.41&2.324E-01&1.40\\
$160\times16\times16$ &8.194E-02&1.50&8.351E-02&1.50&8.291E-02&1.49\\
\midrule
Reference & \multicolumn{6}{c}{$-8.39$}\\
\bottomrule
& \multicolumn{6}{c}{The errors of $w_c$}\\
\cmidrule(lr){2-7}
Mesh& MsFEM-l&$\beta$&LPIML-l&$\beta$&RPIML-l&$\beta$\\
\midrule
$20\times2\times2$ &4.037E+00&&4.279E+00&&4.042E+00&\\
$40\times4\times4$ &1.774E+00&1.19&1.857E+00&1.20&1.779E+00&1.18\\
$80\times8\times8$ &5.819E-01&1.61&6.009E-01&1.63&5.850E-01&1.60\\
$160\times16\times16$ &1.828E-01&1.67&1.876E-01&1.68&1.854E-01&1.66\\
\midrule
Reference & \multicolumn{6}{c}{$18.55$}\\
\bottomrule
\end{tabular}\label{tab:cantilever}
\end{table}

Because these methods use the same assembly and linear solver processes, we report only the CPU time required to compute the numerical shape functions and condensed stiffness matrices of coarse elements in \autoref{tab:times}. The computation was carried out on a computer with an Intel(R) Core(TM) i7-14700KF CPU and 32 GB of RAM. We observe that although LPIML-l is slightly more expensive than RPIML-l, both methods are substantially faster than MsFEM-l.
\begin{table}[width=.9\linewidth,cols=4,pos=h]
\centering
\caption{CPU time for computing the numerical shape functions and condensed stiffness matrices}
\begin{tabular}{l*{3}{c}}   
\toprule
& \multicolumn{3}{c}{The CPU time(s)}\\
\cmidrule(lr){2-4}
Mesh& MsFEM-l&LPIML-l&RPIML-l\\
\midrule
$20\times2\times2$ &0.57313&0.11004&0.09088\\
$40\times4\times4$ &4.40120&0.46050&0.34397\\
$80\times8\times8$ &34.82777&2.86382&2.56659\\
$160\times16\times16$ &277.83698&24.52482&20.65644\\
\bottomrule
\end{tabular}\label{tab:times}
\end{table}
\FloatBarrier

\section{Conclusion}
In this work, we derive error estimates for problem-independent machine learning (PIML) under the assumption that the generalization error associated with learning the transition matrix is given. We show that the total numerical error consists of two components: a discretization error controlled by the mesh size \(h\) and a learning-induced error controlled by the generalization error \(\varepsilon\). To eliminate the non-vanishing learning-induced term in the error bound, we propose an improved PIML scheme that enforces linear completeness. We prove that the resulting method achieves the same optimal convergence order as the corresponding substructure method. The theoretical results are further validated through a series of two- and three-dimensional numerical experiments. Future work will focus on extending the proposed framework to higher-order PDEs and developing the corresponding error estimates. We also note that the proposed framework can be extended to higher-order multiscale finite element methods \cite{41}, provided that the fine-element shape functions possess the required polynomial completeness and corresponding boundary constraints.
\FloatBarrier
\bibliographystyle{unsrtnat}
\bibliography{reference}

@article{6,
  title     = {Physics-informed neural networks: A deep learning framework for solving forward and inverse problems involving nonlinear partial differential equations},
  author    = {Raissi, Maziar and Perdikaris, Paris and Karniadakis, George E.},
  journal   = {Journal of Computational Physics},
  volume    = {378},
  pages     = {686--707},
  year      = {2019},
  doi       = {10.1016/j.jcp.2018.10.045}
}

@article{7,
  title     = {The Deep Learning Galerkin Method for the General Stokes Equations},
  author    = {Li, Jian and Yue, Jing and Zhang, Wen and Duan, Wansuo},
  journal   = {Journal of Scientific Computing},
  volume    = {93},
  number    = {1},
  pages     = {5},
  year      = {2022},
  doi       = {10.1007/s10915-022-01930-8}
}

@article{8,
  title     = {Deep Ritz - Finite element methods: Neural network methods trained with finite elements},
  author    = {Grekas, Georgios and Makridakis, Charalambos G.},
  journal   = {Computer Methods in Applied Mechanics and Engineering},
  volume    = {437},
  pages     = {117798},
  year      = {2025},
  doi       = {10.1016/j.cma.2025.117798}
}

@article{9,
  title     = {On the Convergence of Physics Informed Neural Networks for Linear Second-Order Elliptic and Parabolic Type {PDEs}},
  author    = {Shin, Yeonjong and Darbon, J{\'e}r{\^o}me and Karniadakis, George Em},
  journal   = {Communications in Computational Physics},
  volume    = {28},
  number    = {5},
  pages     = {2042--2074},
  year      = {2020},
  doi       = {10.4208/cicp.OA-2020-0193}
}

@article{10,
  title     = {On the convergence of {PINNs}},
  author    = {Doum{\`e}che, Nathan and Biau, G{\'e}rard and Boyer, Claire},
  journal   = {Bernoulli},
  volume    = {31},
  number    = {3},
  pages     = {2127--2151},
  year      = {2025},
  doi       = {10.3150/24-BEJ1799}
}

@article{11,
  title     = {Convergence and error control of consistent {PINNs} for elliptic {PDEs}},
  author    = {Bonito, Andrea and DeVore, Ronald and Petrova, Guergana and Siegel, Jonathan W.},
  journal   = {IMA Journal of Numerical Analysis},
  volume    = {46},
  number    = {1},
  pages     = {90--148},
  year      = {2026},
  doi       = {10.1093/imanum/draf008}
}

@article{12,
  title     = {Convergence analysis of a {PINNs}-based approach to the inverse source problem of the heat equation with local measurements},
  author    = {Zhang, Xuezhao and Li, Zhiyuan},
  journal   = {Journal of Computational and Applied Mathematics},
  volume    = {476},
  pages     = {117077},
  year      = {2026},
  doi       = {10.1016/j.cam.2025.117077}
}

@article{13,
  title     = {When and why {PINNs} fail to train: A neural tangent kernel perspective},
  author    = {Wang, Sifan and Yu, Xinling and Perdikaris, Paris},
  journal   = {Journal of Computational Physics},
  volume    = {449},
  pages     = {110768},
  year      = {2022},
  doi       = {10.1016/j.jcp.2021.110768}
}

@article{14,
  title     = {On the eigenvector bias of Fourier feature networks: From regression to solving multi-scale {PDEs} with physics-informed neural networks},
  author    = {Wang, Sifan and Wang, Hanwen and Perdikaris, Paris},
  journal   = {Computer Methods in Applied Mechanics and Engineering},
  volume    = {384},
  pages     = {113938},
  year      = {2021},
  doi       = {10.1016/j.cma.2021.113938}
}

@article{15,
  title     = {hp-{VPINNs}: Variational physics-informed neural networks with domain decomposition},
  author    = {Kharazmi, Ehsan and Zhang, Zhongqiang and Karniadakis, George Em},
  journal   = {Computer Methods in Applied Mechanics and Engineering},
  volume    = {374},
  pages     = {113547},
  year      = {2021},
  doi       = {10.1016/j.cma.2020.113547}
}

@article{16,
  title     = {Stochastic quadrature rules for solving {PDEs} using neural networks},
  author    = {Taylor, Jamie M. and Pardo, David},
  journal   = {Computers \& Mathematics with Applications},
  volume    = {214},
  pages     = {1--26},
  year      = {2026},
  doi       = {10.1016/j.camwa.2026.04.002}
}

@article{17,
  title     = {Variational Physics Informed Neural Networks: the Role of Quadratures and Test Functions},
  author    = {Berrone, Stefano and Canuto, Claudio and Pintore, Moreno},
  journal   = {Journal of Scientific Computing},
  volume    = {92},
  number    = {3},
  pages     = {100},
  year      = {2022},
  doi       = {10.1007/s10915-022-01950-4}
}

@book{18,
  author    = {Ciarlet, Philippe G.},
  title     = {The Finite Element Method for Elliptic Problems},
  publisher = {North-Holland},
  address   = {Amsterdam},
  year      = {1978}
}

@article{19,
  title     = {A new 8-node quadrilateral spline finite element},
  author    = {Li, Chong-Jun and Wang, Ren-Hong},
  journal   = {Journal of Computational and Applied Mathematics},
  volume    = {195},
  number    = {1--2},
  pages     = {54--65},
  year      = {2006},
  doi       = {10.1016/j.cam.2005.07.017}
}

@article{20,
  title     = {Virtual Elements for Linear Elasticity Problems},
  author    = {Beir{\~a}o da Veiga, L. and Brezzi, F. and Marini, L. D.},
  journal   = {SIAM Journal on Numerical Analysis},
  volume    = {51},
  number    = {2},
  pages     = {794--812},
  year      = {2013},
  doi       = {10.1137/120874746}
}

@article{21,
  title     = {The scaled boundary finite-element method---alias consistent infinitesimal finite-element cell method---for elastodynamics},
  author    = {Song, Chongmin and Wolf, John P.},
  journal   = {Computer Methods in Applied Mechanics and Engineering},
  volume    = {147},
  number    = {3--4},
  pages     = {329--355},
  year      = {1997},
  doi       = {10.1016/S0045-7825(97)00021-2}
}

@article{22,
  title     = {Problem-independent machine learning ({PIML})-based topology optimization---A universal approach},
  author    = {Huang, Mengcheng and Du, Zongliang and Liu, Chang and Zheng, Yonggang and Cui, Tianchen and Mei, Yue and Li, Xiao and Zhang, Xiaoyu and Guo, Xu},
  journal   = {Extreme Mechanics Letters},
  volume    = {56},
  pages     = {101887},
  year      = {2022},
  doi       = {10.1016/j.eml.2022.101887}
}

@article{23,
  title     = {A Problem-Independent Machine Learning ({PIML}) enhanced substructure-based approach for large-scale structural analysis and topology optimization of linear elastic structures},
  author    = {Huang, Mengcheng and Cui, Tianchen and Liu, Chang and Du, Zongliang and Zhang, Jiameng and He, Chuhui and Guo, Xu},
  journal   = {Extreme Mechanics Letters},
  volume    = {63},
  pages     = {102041},
  year      = {2023},
  doi       = {10.1016/j.eml.2023.102041}
}

@article{24,
  title     = {Problem-independent machine learning-enhanced structural topology optimization of complex design domains based on isoparametric elements},
  author    = {Zhang, Linfeng and Huang, Mengcheng and Liu, Chang and Du, Zongliang and Cui, Tianchen and Guo, Xu},
  journal   = {Extreme Mechanics Letters},
  volume    = {72},
  pages     = {102237},
  year      = {2024},
  doi       = {10.1016/j.eml.2024.102237}
}

@article{30,
  title     = {A mechanics-based data-free Problem Independent Machine Learning ({PIML}) model for large-scale structural analysis and design optimization},
  author    = {Huang, Mengcheng and Liu, Chang and Guo, Yilin and Zhang, Linfeng and Du, Zongliang and Guo, Xu},
  journal   = {Journal of the Mechanics and Physics of Solids},
  volume    = {193},
  pages     = {105893},
  year      = {2024},
  doi       = {10.1016/j.jmps.2024.105893}
}

@inproceedings{32,
  title     = {Learning from Integral Losses in Physics Informed Neural Networks},
  author    = {Saleh, Ehsan and Ghaffari, Saba and Bretl, Tim and Olson, Luke and West, Matthew},
  booktitle = {Proceedings of the 41st International Conference on Machine Learning},
  pages     = {43077--43111},
  year      = {2024},
  volume    = {235},
  series    = {Proceedings of Machine Learning Research},
  month     = {21--27 Jul},
  publisher = {PMLR},
  url       = {https://proceedings.mlr.press/v235/saleh24a.html}
}

@article{33,
  author  = {Zhang, H. W. and Liu, Y. and Zhang, S. and Tao, J. and Wu, J. K. and Chen, B. S.},
  title   = {Extended multiscale finite element method: its basis and applications for mechanical analysis of heterogeneous materials},
  journal = {Computational Mechanics},
  year    = {2014},
  volume  = {53},
  number  = {4},
  pages   = {659--685},
  doi     = {10.1007/s00466-013-0924-x}
}

@book{34,
  author    = {Ivakhnenko, A. G. and Lapa, V. G.},
  title     = {Cybernetic Predicting Devices},
  year      = {1965},
  publisher = {CCM Information Corporation},
  address   = {New York, NY, USA}
}

@article{35,
  title={Learning nonlinear operators via DeepONet based on the universal approximation theorem of operators},
  author={Lu, Lu and Jin, Pengzhan and Pang, Guofei and Zhang, Zhongqiang and Karniadakis, George Em},
  journal={Nature Machine Intelligence},
  volume={3},
  number={3},
  pages={218--229},
  year={2021},
  doi={10.1038/s42256-021-00302-6}
}

@article{36,
  title={DGM: A deep learning algorithm for solving partial differential equations},
  author={Sirignano, Justin and Spiliopoulos, Konstantinos},
  journal={Journal of Computational Physics},
  volume={375},
  pages={1339--1364},
  year={2018},
  doi={10.1016/j.jcp.2018.08.029}
}

@article{37,
  title={The Deep Ritz method: A deep learning-based numerical algorithm for solving variational problems},
  author={E, Weinan and Yu, Bing},
  journal={Communications in Mathematics and Statistics},
  volume={6},
  number={1},
  pages={1--12},
  year={2018},
  doi={10.1007/s40304-018-0127-z}
}

@inproceedings{38,
  title = 	 {On the Spectral Bias of Neural Networks},
  author =       {Rahaman, Nasim and Baratin, Aristide and Arpit, Devansh and Draxler, Felix and Lin, Min and Hamprecht, Fred and Bengio, Yoshua and Courville, Aaron},
  booktitle = 	 {Proceedings of the 36th International Conference on Machine Learning},
  pages = 	 {5301--5310},
  year = 	 {2019},
  volume = 	 {97},
  series = 	 {Proceedings of Machine Learning Research},
  month = 	 {09--15 Jun},
  publisher =    {PMLR},
  url = 	 {https://proceedings.mlr.press/v97/rahaman19a.html}
}

@article{39,
  title={Physics-informed machine learning},
  author={Karniadakis, George Em and Kevrekidis, Ioannis G and Lu, Lu and Perdikaris, Paris and Wang, Sifan and Yang, Liu},
  journal={Nature Reviews Physics},
  volume={3},
  number={6},
  pages={422--440},
  year={2021},
  doi={10.1038/s42254-021-00314-5}
}

@article{40,
  title={Physics-informed neural networks with hard constraints for inverse design},
  author={Lu, Lu and Pestourie, Rapha{\"e}l and Yao, Wenjie and Wang, Zhicheng and Verdugo, Francesc and Johnson, Steven G.},
  journal={SIAM Journal on Scientific Computing},
  volume={43},
  number={6},
  pages={B1105--B1132},
  year={2021},
  doi={10.1137/21M1397908}
}

@article{41,
  author  = {Wang, Yitong and Yang, Yanfang},
  title   = {Higher-order multiscale finite element methods for elliptic equations with oscillating coefficients},
  journal = {Journal of Computational and Applied Mathematics},
  volume  = {472},
  pages   = {116835},
  year    = {2026},
  doi     = {10.1016/j.cam.2025.116835}
}

@article{42,
  author  = {Wu, Zhetong and Zhang, Hanbo and Ye, Hongfei and Zhang, Hongwu and Zheng, Yonggang and Guo, Xu},
  title   = {{PINN} enhanced extended multiscale finite element method for fast mechanical analysis of heterogeneous materials},
  journal = {Acta Mechanica},
  year    = {2024},
  volume  = {235},
  number  = {7},
  pages   = {4895--4913},
  doi     = {10.1007/s00707-024-03984-1}
}

@article{chen2010spline,
  title     = {A family of spline finite elements},
  author    = {Chen, Juan and Li, Chong-Jun and Chen, Wan-Ji},
  journal   = {Computers \& Structures},
  volume    = {88},
  number    = {11--12},
  pages     = {718--727},
  year      = {2010},
  doi       = {10.1016/j.compstruc.2010.02.011}
}


\end{document}